\documentclass[a4paper,leqno]{amsart}

\usepackage[pagebackref=true]{hyperref}
\renewcommand*{\backref}[1]{}
\renewcommand*{\backrefalt}[4]{{[\tiny\ifcase #1 Not cited.\relax\or Page~#2.\else Pages #2.\fi]}}

\usepackage[svgnames]{xcolor}
\usepackage[hyperpageref]{backref}
\hypersetup{pdftitle={Cyclotomic Carter-Payne homomorphisms},
  pdfauthor={Sinead Lyle and Andrew Mathas},
  colorlinks = true,
  linkcolor =ForestGreen,
  anchorcolor = red,
  citecolor = blue,
  urlcolor = blue
}

\usepackage{booktabs}
\usepackage{amsmath,amssymb,mathrsfs,latexsym,tikz,cite}
\usetikzlibrary{calc,intersections,through,backgrounds}

\newcommand\Comment[2][\relax]{\space\par\medskip\noindent%
   \fbox{\begin{minipage}{\textwidth}\textbf{Comment\ifx\relax#1\else---#1\fi}\newline%
        #2\end{minipage}}\medskip
}

\colorlet{darkgreen}{green!60!black}
\tikzset{dots/.style={very thick,loosely dotted,color=blue!50},
         greendot/.style={fill,circle,color=darkgreen,inner sep=1.5pt,outer sep=0}
}
\def\greendot(#1);{\node[greendot] at(#1){};}

\newenvironment{braid}[1][10]{
  \begin{tikzpicture}[baseline=#1mm,blue,line width=1pt, xscale=0.4,yscale=0.4,
                      draw/.append style={rounded corners},
                      every node/.append style={font=\fontsize{5}{5}\selectfont}]%
  }{\end{tikzpicture}
}
\newcommand\iicrossing[1][i]{%
  \begin{braid}\tikzset{baseline=0mm,scale=0.5}
    \useasboundingbox (-0.6,0) rectangle (1,1.9);
    \draw(0,1)node[left=-1mm]{$#1$}--(1,0);
    \draw(1,1)node[right=-1mm]{$#1$}--(0,0);
  \end{braid}%
}
\newcommand\ijcrossing[1][j]{%
  \begin{braid}\tikzset{baseline=0mm,scale=0.5}
    \useasboundingbox (-0.6,0) rectangle (1,1.9);
    \draw(0,1)node[left=-1mm]{$i$}--(1,0);
    \draw(1,1)node[right=-1mm]{$#1$}--(0,0);
  \end{braid}%
}

\def\BraidHeight{4}
\newcommand\BraidTableau[2][3]{
  \newcount\RowRes\RowRes=0
  \newcount\N\N=1
  \newcount\tabRes
  \foreach \tabRow in {#2} {
    \tabRes=\RowRes
    \ifnum\N>1
      \draw[brown](\N-0.5,\BraidHeight+1) -- +(0,-1);  
      \draw[brown](\N-0.5,0) -- +(0,-1);          
    \fi
    \foreach \tabCol in \tabRow {
      \draw[black](\the\N,\BraidHeight+1.5)node{$\the\N$};
      \ifnum\tabCol>0  
        \draw(\the\N,\BraidHeight)node[above]{$\the\tabRes$}
              --(\tabCol,0)node[below]{$\the\tabRes$};
      \else\ifnum\tabCol<0 
        \draw[red](\the\N,\BraidHeight)node[above]{$\the\tabRes$}
              --(-\tabCol,0)node[below]{$\the\tabRes$};
      \fi\fi
      \global\advance\tabRes by 1
      \ifnum\tabRes=#1\global\tabRes=0\fi
      \global\advance\N by 1
    }
    \global\advance\RowRes by -1
    \ifnum\RowRes<0\global\advance\RowRes by #1\fi
  }
  \draw[brown](\N-0.5,\BraidHeight+1)--(0.5,\BraidHeight+1)--+(0,-1)--(\N-0.5,\BraidHeight)--cycle;
  \draw[brown](\N-0.5,0)--(0.5,0)--(0.5,-1)--(\N-0.5,-1)--cycle;
}

\newcommand\YCrossing[5]{%
  \coordinate (A) at (#2,\BraidHeight);
  \coordinate (B) at (#3,\BraidHeight);
  \coordinate (C) at (#4,0);
  \coordinate (D) at (#5,0);
  \draw[red](A)node[above]{$#1$};
  \draw[red](B)node[above]{$#1$};
  \draw[red](C);
  \draw[red](D);
  \path[name path=AD](A)--(D);
  \path[name path=BC](B)--(C);
  \path[name intersections={of=AD and BC,by=I}];
  \draw[red](A)--($ (A)!0.9!(I) $) -- ($ (I)!0.1!(C) $)--(C);
  \draw[red](B)--($ (B)!0.9!(I) $) -- ($ (I)!0.1!(D) $)--(D);
}

\newcount\tableauRow\newcount\tableauCol
\newcommand\TikTableau[3]{
  \begin{tikzpicture}[scale=#1,draw/.append style={thick,black},baseline=#2mm]
    \tableauRow=0
    \foreach \Row in {#3} {
       \tableauCol=1
       \foreach\k in \Row {
          \draw(\the\tableauCol,\the\tableauRow)+(-.5,-.5)rectangle++(.5,.5);
          \draw(\the\tableauCol,\the\tableauRow)node{\k};
          \global\advance\tableauCol by 1
       }
       \global\advance\tableauRow by -1
    }
  \end{tikzpicture}
}
\newcommand\Tableau[2][4]{\TikTableau{0.5}{#1}{#2}}
\newcommand\tableau[1]{{%
  \tikzset{every node/.append style={font=\fontsize{6}{6}\selectfont}}%
  \TikTableau{0.28}{-2.8}{#1}}%
}

\newdimen\shadedBaseline\shadedBaseline=-4mm
\newcommand\ShadedTableau[2][\relax]{
  \begin{tikzpicture}[scale=0.5,draw/.append style={thick,black},baseline=\shadedBaseline]
    \ifx\relax#1\relax%
    \else 
      \foreach\box in {#1} { \filldraw[blue!20]\box+(-.5,-.5)rectangle++(.5,.5); }
    \fi
    \tableauRow=0
    \foreach \Row in {#2} {
       \tableauCol=1
       \foreach\k in \Row {
          \draw(\the\tableauCol,\the\tableauRow)+(-.5,-.5)rectangle++(.5,.5);
          \draw(\the\tableauCol,\the\tableauRow)node{\k};
          \global\advance\tableauCol by 1
       }
       \global\advance\tableauRow by -1
    }
  \end{tikzpicture}
}

\def\Tritab(#1|#2|#3){\Bigg(%
  \ \Tableau[-4]{#1}\ \Bigg|\ \Tableau[-4]{#2}%
  \ \Bigg|\ \Tableau[-3]{#3}\ \
Bigg)
}

\newcommand\YoungDiagram[2][0]{
  \begin{tikzpicture}[scale=0.5,draw/.append style={thick,black},baseline=#1mm]
    \newcount\tabRow\tabRow=0
    \foreach \col in {#2} {
       \draw(1,\the\tabRow)grid ++(\col,1);
       \global\advance\tabRow by -1
    }
  \end{tikzpicture}
}
\newcommand\TriDiagram[4][0]{\Bigg(%
   \ \YoungDiagram[#1]{#2}\ \Bigg|%
   \ \YoungDiagram[#1]{#3}\ \Bigg|%
   \ \YoungDiagram[#1]{#4}\ \Bigg)
}

\title{Cyclotomic Carter-Payne homomorphisms}

\def\Email#1{\email{\href{mailto:#1}{#1}}}
\author{Sin\'ead Lyle}
\address{School of Mathematics \\ University of East Anglia \\ Norwich NR4 7TJ \\ UK}
\Email{s.lyle@uea.ac.uk}
\author{Andrew Mathas}
\address{School of Mathematics and Statistics \\ University of Sydney \\ NSW 2006 \\ Australia}
\Email{andrew.mathas@sydney.edu.au}

\numberwithin{equation}{section}
\swapnumbers

\makeatletter
\def\subsectionautorefname{\S\@gobble}
\makeatother

\usepackage{aliascnt}
\def\NewTheorem#1{%
  \newaliascnt{#1}{equation}
  \newtheorem{#1}[#1]{#1}
  \aliascntresetthe{#1}
  \expandafter\def\csname #1autorefname\endcsname{#1}
}
\def\equationautorefname~#1\null{(#1)\null}

\newcounter{main}
\theoremstyle{plain}

\newtheorem*{MainTheorem}{Main Theorem}

\NewTheorem{Assumption}
\NewTheorem{Proposition}
\NewTheorem{Theorem}
\NewTheorem{Corollary}
\newtheorem*{COROLLARY}{Corollary}
\NewTheorem{Lemma}
\NewTheorem{Definition}
\theoremstyle{definition}
\theoremstyle{remark}
\NewTheorem{Remark}

\newaliascnt{Example}{equation}
\aliascntresetthe{Example}

\newenvironment{Example}%
  {\refstepcounter{Example}\trivlist
   \item[\hskip\labelsep\theExample.~\textbf{Example}\space]
   \ignorespaces
  }{\unskip\nobreak\hfil%
    \penalty50\hskip2em\hbox{}\nobreak\hfil$\Diamond$%
    \parfillskip=0pt\finalhyphendemerits=0\penalty-100\endtrivlist
    \medskip
}

\newcommand\blam{{\boldsymbol\lambda}}
\newcommand\bmu{{\boldsymbol\mu}}
\newcommand\bnu{{\boldsymbol\nu}}
\newcommand\bsig{{\boldsymbol\sigma}}
\newcommand\brho{{\boldsymbol\rho}}
\newcommand\btau{{\boldsymbol\tau}}
\newcommand\charge{{\boldsymbol\kappa}}

\newcommand{\Sym}{\mathfrak{S}}

\def\D{\mathscr D}

\newcommand\thea[1][\blam]{\theta^{#1}_\bmu}

\let\<=\langle
\let\>=\rangle
\def\({\big(}
\def\){\big)}

\def\bi{\mathbf{i}}
\def\bj{\mathbf{j}}
\newcommand\iarrow[1][i]{\xrightarrow{\ #1\ }}
\newcommand\nuarrow[1][\nu]{\iarrow[\nu]}

\newcommand\tA[1][e]{\mathsf{T}^{(#1)}_A}

\def\s{\mathfrak{s}}
\def\t{\mathfrak{t}}
\def\tmu{\t^{\bmu}}
\def\a{\mathfrak{a}}
\def\u{\mathfrak{u}}
\def\tnulam{{\t^\bnu_{\blam}}}
\def\tnu{\t^\bnu}
\def\tnumu{{\t^\bnu_{\bmu}}}

\newcommand\tnusig[1][\bsig]{{\t^\bnu_{#1}}}
\def\tnutau{{\t^\bnu_{\btau}}}
\def\tlam{\t^\blam}
\def\N{\mathbb{N}}
\def\Z{\mathbb{Z}}
\def\LMult{L^\blam_\bmu}

\def\W{\mathcal{W}_{n+\gamma}}
\def\floor#1/#2{\lfloor\tfrac{#1}{#2}\rfloor}
\def\ceil#1/#2{\lceil\tfrac{#1}{#2}\rceil}

\newcommand\Belt[1][e]{\mathbf{B}^{(#1)}_A}

\let\gedom=\trianglerighteq
\let\gdom=\vartriangleright

\newcommand\Parts[1][n]{\mathscr{P}^\Lambda_{#1}}
\renewcommand\H[1][n]{\mathscr{H}^\Lambda_{#1}}
\newcommand\R[1][n]{\mathscr{R}^\Lambda_{#1}}
\newcommand\UnR[1][n]{\underline{\mathscr{R}}^\Lambda_{#1}}

\def\m{{\mathbf{m}}}

\newcommand\Psin[1][n]{\boldsymbol\Psi_{#1}}
\newcommand\Sp[1][\blam]{\underline{S}^{#1}}

\DeclareMathOperator\Add{Add}
\DeclareMathOperator\Rem{Rem}

\DeclareMathOperator\row{row}
\DeclareMathOperator\res{res}
\DeclareMathOperator\Ext{Ext}
\DeclareMathOperator\Hom{Hom}
\DeclareMathOperator\Mod{Mod}

\DeclareMathOperator\Std{Std}
\newcommand{\IniNu}[1][n]{\Std^e_{#1}(\bnu)}
\newcommand\Stdlm[1][k]{\Std_{\blam}(\bnu,#1)}
\newcommand{\StdCheck}[1][\bsig]{\Std^\vee_{#1}(\bnu)}
\newcommand{\StdHat}[1][\bsig]{\Std^\wedge_{#1}(\bnu)}

\DeclareMathOperator\Shape{Shape}

\DeclareMathOperator\Res{Res}

{\catcode`\|=\active
  \gdef\set#1{\mathinner{\lbrace\,{\mathcode`\|"8000%
  \let|\midvert #1}\,\rbrace}}
}
\def\midvert{\egroup\mid\bgroup}
\def\map#1#2{\,{:}\,#1\!\longrightarrow\!#2}

\newcommand\Set[3][80]{\Big\{\ #2\ \Big| \vcenter{\hsize #1mm\centering#3}\Big\}}
\newcommand\SET[3][80]{\Bigg\{\ #2\ \Bigg| \vcenter{\hsize #1mm\centering#3}\Bigg\}}

{\catcode`\|=\active
  \gdef\set#1{\mathinner{\lbrace\,{\mathcode`\|"8000%
                                   \let|\midvert #1}\,\rbrace}}
}

\begin{document}
\maketitle

\begin{abstract}
  We construct a new family of homomorphisms between (graded) Specht modules of
  the quiver Hecke algebras of type~$A$. These maps have many
  similarities with the homomorphisms constructed by Carter and Payne in the
  special case of the symmetric groups, although the maps that we obtain are
  both more and less general than these.
\end{abstract}

\section{Introduction}
The degenerate and non-degenerate cyclotomic Hecke algebras of type
$G(\ell,1,n)$ are an important class of algebras that arise naturally in the
representation theory of the groups of Lie type~\cite{Broue:conjectures}, the
study of rational Cherednik algebras~\cite{GordonLosev:CycCatOCherednik} and
the categorification of the irreducible highest weight modules of the affine
special linear groups~\cite{Ariki:can,BK:GradedDecomp}.

In the representation theory of non-semisimple algebras the key unsolved
problems revolve around computing decomposition multiplicities and $\Ext^n$-spaces
between important classes of modules. For the cyclotomic Hecke algebras $\H$ of
type $G(\ell,1,n)$ significant progress has been made on the decomposition
number problem through categorification, however, the calculation of
$\Ext^n$-spaces for Specht modules remains an open problem.

The simplest $\Ext^n$-spaces are the hom-spaces and even here very little is known.
For example, it was only recently shown by Dodge~\cite{Dodge:LargeHoms} that the
dimension of the hom-spaces between Specht modules of the symmetric groups can
be arbitrarily large. This was quite surprising because, prior to Dodge's work,
there were no known examples of hom-spaces between Specht modules in odd
characteristic which had dimension greater than one.  Apart from the work in
this paper, the only results on hom-spaces between Specht modules for the
cyclotomic Hecke algebras are those contained in the recent work of
Corlett~\cite{Corlett:PhD,Corlett:OneNode}.

We construct a new family of explicit non-zero homomorphisms between the Specht
modules of cyclotomic Hecke algebras.  Our results are a cyclotomic
generalisation of the famous Carter-Payne Theorem~\cite{CarterPayne} for the
symmetric groups. In the special case of the symmetric group we construct many
of the homomorphisms described by Carter and Payne. We also construct additional
maps which are not Carter-Payne maps.

The main new tool that we use to construct our cyclotomic Carter-Payne
homomorphisms is to work in the graded setting and to construct homomorphisms
for the Specht modules, defined over~$\Z$, for the cyclotomic quiver Hecke
algebras, or cyclotomic KLR algebras,  of type~$A$. These algebras are certain
$\Z$-graded algebras which, over a field, are isomorphic to the cyclotomic Hecke
algebras of type~$G(\ell,1,n)$ by work of Brundan and
Kleshchev~\cite{BK:GradedKL}. We develop a number of new tools for working in
the graded setting which are likely to be of independent interest; see, for
example, \autoref{S:Stubborn}.

To state our main theorem, let $\R$ be the cyclotomic quiver Hecke algebra of
type~$A$ over~$\Z$ which is determined by the quiver with vertex set
$I=\Z/e\Z$, where $e\in\{0,2,3,4,\dots\}$. For each multipartition $\blam$ there
is a graded Specht module~$S^\blam$ which is a $\Z$-free
$\R$-module~\cite{BKW:GradedSpecht}. If $K$ is a field then $S^\blam\otimes_\Z K$
is isomorphic to the graded Specht module constructed
in~\cite{HuMathas:GradedCellular,KMR:UniversalSpecht} and, in turn, this
module is a graded lift of the (ungraded) Specht module $\Sp$ of the
cyclotomic Hecke algebras~$\H$~\cite{DJM:cyc,AMR}.

In \autoref{D:CycCPPair} we give a purely combinatorial, but quite technical,
condition for when a pair $(\blam,\bmu)$ of multipartitions is a
\textbf{cyclotomic Carter-Payne pair}. The most important special case is when
the multipartition $\bmu$ can be obtained from $\blam$ by moving a horizontal
strip of $\gamma$ nodes to an earlier row in the diagram of~$\blam$ without
changing their residues, where $\gamma<|I|$. Our main result is the
following.

\begin{MainTheorem}
\label{MainTheorem}
  Suppose that $(\blam,\bmu)$ is a cyclotomic Carter-Payne pair, where~$\blam$
  and~$\bmu$ are multipartitions of~$n$. Then
  $\Hom_{\R}(S^\blam\<\delta\>,S^\bmu)\ne0$, where $\delta$ is a positive
  integer determined by $\blam$ and $\bmu$.
\end{MainTheorem}

An explicit formula for the integer $\delta$ is given in \autoref{T:main}.
Applying Brundan and Kleshchev's graded isomorphism theorem~\cite{BK:GradedKL}
we obtain the corresponding result for the cyclotomic Hecke algebra~$\H$.

\begin{COROLLARY}
  Suppose that $\H$ is the cyclotomic Hecke algebra defined over a field~$K$ and
  that $(\blam,\bmu)$ is a cyclotomic Carter-Payne pair. Then
  $\Hom_{\H}(\Sp,\Sp[\bmu])\ne0$.
\end{COROLLARY}

The Carter-Payne theorem was proved for symmetric groups in 1980. It proved very
difficult to generalise this result to the Iwahori-Hecke algebras of the
symmetric groups~\cite{Lyle:1CarterPayne,Dixon,LM:CarterPayne}, which is the
level one case of the algebras considered here. We take a very different
approach to this problem here, using use the machinery of the cyclotomic quiver
Hecke algebras to construct homomorphisms $S^\blam\<\delta\>\rightarrow S^\bmu$.

The homomorphisms of our Main Theorem are described explicitly as multiplication
by certain polynomials in the KLR generators $y_{n+1},\dots,y_{n+\gamma}$, which
are nilpotent versions of the classical Jucys-Murphy elements. These polynomials
act on a Specht module for~$\R[n+\gamma]$, which we consider as an $\R$-module
by restriction, and which admits a filtration in which the Specht modules $S^\blam$ and
$S^\bmu$ appear as quotients.
 It is not very difficult to show that the $\R$-endomorphism that we define on $S^\bnu$
induces an $\R$-module homomorphism between these Specht modules,
however, the real challenge is in showing that our maps
are non-zero. Rather than working algebraically, it turns out to be much easier
to use a variation of the diagram calculus introduced by Khovanov and
Lauda~\cite{KhovLaud:diagI}. For the uninitiated, it takes some time to get used
to these diagrams, but it is worth the effort because they reduce long pages of
calculations with relations to a series of diagrams which neatly encode the same
information in a more transparent way.

As a gentle introduction to the use of these diagrams we start by proving our
Main Theorem in \textit{one node case}, which is when $\gamma=1$.  This case
contains most of the features of the general case, however, its proof is
considerably easier both to write down and to understand. The basic idea for
proving \autoref{MainTheorem} when $\gamma>1$ is the same as the one node case,
however, the general argument requires more sophisticated techniques for working
with diagrams and tableaux. For example, in \autoref{S:Stubborn} we introduce
the notion of \textit{stubborn strings} which we think will have wider
application. Using these ideas we prove our Main Theorem as \autoref{T:main}.

As a general rule it is very hard to construct non-trivial homomorphisms between
Specht modules and it is even harder to compose such homomorphisms. One
advantage of our explicit construct is that we are able to compose some of our
cyclotomic Carter-Payne homomorphisms in \autoref{T:composition}. As a second
application we prove a  `row removal' theorem for cyclotomic Carter-Payne
homomorphisms \autoref{T:Slice}.

\section{Cyclotomic Hecke algebras and Specht modules}
The cyclotomic quiver Hecke algebras were introduced for all oriented quivers in a
series of papers by Khovanov-Lauda~\cite{KhovLaud:diagI,KhovLaud:diagII} and
Rouquier~\cite{Rouq:2KM}. These algebras are certain $\Z$-graded algebras which
depend on an oriented quiver.

This section defines the cyclotomic quiver Hecke algebras of type~$A$ and
recalls the results from the literature that we need. Throughout this paper we
work over the ring of integers~$\Z$. Hence, by base change, all of our results
hold over an arbitrary commutative ring.

\subsection{Graded algebras and homomorphisms}
In this paper a \textbf{graded algebra} will mean a $\Z$-graded algebra and a
\textbf{graded module} will be a $\Z$-graded module. If $A$ is a graded algebra
then $\Mod-{A}$ is the category of finitely generated graded (right) $A$-modules
with degree preserving maps. We use the standard notation of graded
representation theory. In particular, if~$M=\bigoplus_{d\in\Z}M_d$ then
$m\in M_d$ is \textbf{homogeneous} of \textbf{degree} $d=\deg m$. If~$n\in\Z$
then $M\<n\>$ is the graded module obtained by shifting the grading on $M$ up
by~$n$ so that $M\<n\>_d=M_{d-n}$.

If $A$ is a graded algebra let $\underline{A}$ be the ungraded algebra obtained
by forgetting the grading on~$A$. Similarly, $\underline{M}$ is the (ungraded)
$\underline{A}$-module obtained by forgetting the grading on the graded
$A$-module $M$. It is well-known that if $M$ and $N$ are graded $A$-modules then
$$\Hom_{\underline{A}}(\underline{M},\underline{N})
    \cong\bigoplus_{d\in\Z}\Hom_A(M\<d\>,N),$$
where $\Hom_A(M\<d\>,N)\cong\Hom_A(M,N\<-d\>)$ is isomorphic to the space of
homogeneous $A$-module homomorphisms $f\map{M}{N}$ of degree~$d$ such that
$f(M_z)\subseteq N_{z+d}$, for $z\in\Z$.

\subsection{Cyclotomic quiver Hecke algebras}
\label{S:grading}
We are now ready to define the cyclotomic quiver Hecke algebras of type~$A$.

Recall from the introduction that we have fixed an integer $e\in\{0,2,3,4,\dots\}$ and that
$I=\Z/e\Z$. Let $\Gamma_e$ be the oriented quiver with vertex set~$I$ and directed edges
$i\longrightarrow i-1$, for $i\in I$.
(The orientation of~$\Gamma_e$ is opposite to that used
in~\cite{BK:GradedKL,HuMathas:GradedCellular}.)
Thus $\Gamma$ is a quiver of type~$A_\infty$ if $e=0$
or of type $A_{e-1}^{(1)}$ if $e > 0$. The corresponding \textbf{Cartan matrix}
$(a_{i,j})_{i, j \in I}$ is given by
\begin{equation}
\label{ECM}
a_{i,j} =
\begin{cases}
2&\text{if $i=j$},\\
0&\text{if $j \neq i, i \pm 1$},\\
-1&\text{if $i \rightarrow j$ or $i \leftarrow j$},\\
-2&\text{if $i \rightleftarrows j$}.
\end{cases}
\end{equation}
(The case $a_{i,j}=-2$ only occurs if $e=2$.) To the quiver $\Gamma_e$ we attach
the standard Lie theoretic data of simple roots $\set{\alpha_i |i \in I}$, fundamental weights $\set{\Lambda_i|i\in I}$,
the positive weight lattice $P^+=\bigoplus_{i\in I}\N\Lambda_i$ and  positive roots
$Q^+=\bigoplus_{i\in I}\N\alpha_i$ and we let $(\;,\;)$ be the bilinear
form determined by
$$(\alpha_i,\alpha_j)=a_{ij}\qquad\text{and}\qquad
          (\Lambda_i,\alpha_j)=\delta_{ij},\qquad\text{for }i,j\in I.$$
More details can be found, for example, in Kac's book~\cite[Chapter~1]{Kac}.

Let $\Sym_n$ be the \textbf{symmetric group} on $n$ letters and let $s_r=(r,r+1)$,
for $1\le r<n$, be the simple transpositions of~$\Sym_n$. Then $\Sym_n$ acts
from the left on elements of the set $ I^n$ by place permutations.

\begin{Definition}[Khovanov-Lauda~\cite{KhovLaud:diagI,KhovLaud:diagII},
                   Rouquier~\cite{Rouq:2KM}]
  \label{D:KLR}
  Suppose that $\Lambda\in P^+$ and $n\ge1$. The \textbf{cyclotomic quiver Hecke
  algebra} $\R=\R(\Z)$ of weight $\Lambda$ and type $\Gamma_e$ is the unital
  associative $\Z$-algebra with generators
  $$\{\psi_1,\dots,\psi_{n-1}\} \cup
  \{ y_1,\dots,y_n \} \cup \set{e(\bi)|\bi\in I^n}$$
  and relations
  {\setlength{\abovedisplayskip}{2pt}
   \setlength{\belowdisplayskip}{1pt}
  \begin{align*}
    y_1^{(\Lambda,\alpha_{i_1})}e(\bi)&=0,
    & e(\bi) e(\bj) &= \delta_{\bi\bj} e(\bi),
    &{\textstyle\sum_{\bi \in I^n}} e(\bi)&= 1,\\
    y_r e(\bi) &= e(\bi) y_r,
    &\psi_r e(\bi)&= e(s_r{\cdot}\bi) \psi_r,
    &y_r y_s &= y_s y_r,
  \end{align*}
  \begin{align}
    \psi_r y_{r+1} e(\bi)&=(y_r\psi_r+\delta_{i_ri_{r+1}})e(\bi),&
    y_{r+1}\psi_re(\bi)&=(\psi_r y_r+\delta_{i_ri_{r+1}})e(\bi),\label{E:ypsi}\\
    \psi_r y_s  &= y_s \psi_r,&&\text{if }s \neq r,r+1,\notag\\
    \psi_r \psi_s &= \psi_s \psi_r,&&\text{if }|r-s|>1,\label{E:ypsiCommute}
  \end{align}
  \begin{align}
    \psi_r^2e(\bi) &= \begin{cases}
         0,&\text{if }i_r = i_{r+1},\\
         (y_{r+1}-y_r)e(\bi),&\text{if  }i_r\rightarrow i_{r+1},\\
         (y_r - y_{r+1})e(\bi),&\text{if }i_r\leftarrow i_{r+1},\\
         (y_{r+1} - y_{r})(y_{r}-y_{r+1}) e(\bi),
              &\text{if }i_r\rightleftarrows i_{r+1}\\
          e(\bi),&\text{otherwise},\\
    \end{cases}
\label{E:psi^2}\\
    \psi_{r}\psi_{r+1} \psi_{r} e(\bi) &= \begin{cases}
      (\psi_{r+1} \psi_{r} \psi_{r+1} +1)e(\bi),\hspace*{7mm}
          &\text{if }i_r=i_{r+2}\rightarrow i_{r+1} ,\\
    (\psi_{r+1} \psi_{r} \psi_{r+1} -1)e(\bi),
        &\text{if }i_r=i_{r+2}\leftarrow i_{r+1},\\
    \rlap{$\big(\psi_{r+1} \psi_{r} \psi_{r+1} +y_r -2y_{r+1}+y_{r+2}\big)e(\bi)$,}\\
             &\text{if }i_r=i_{r+2} \rightleftarrows i_{r+1},\\
    \psi_{r+1} \psi_{r} \psi_{r+1} e(\bi),&\text{otherwise,}
    \end{cases}
\label{E:braid}
  \end{align}
  }
for $\bi,\bj\in I^n$ and all admissible $r$ and $s$. Moreover, $\R$ is naturally
$\Z$-graded with degree function determined by
$$\deg e(\bi)=0,\qquad \deg y_r=2\qquad\text{and}\qquad \deg
  \psi_s e(\bi)=-a_{i_s,i_{s+1}},$$
for $1\le r\le n$, $1\le s<n$ and $\bi\in I^n$.
\end{Definition}

Throughout the paper, we fix a dominant weight $\Lambda\in P^+$ and an ordered
$\ell$-tuple $\charge=(\kappa_1,\dots,\kappa_\ell)\in\Z^\ell$, the
\textbf{multicharge}, such that
$\ell=\sum_{i\in I}(\Lambda,\alpha_i)$ and
$\Lambda=\Lambda_{\kappa_1}+\dots+\Lambda_{\kappa_\ell}$. The algebra~$\R$
depends only on~$\Lambda$ (and on $\Gamma_e$ and $n$), however, some of the
tableau combinatorics which we introduce below will depend on the choice
of~$\charge$.

Suppose that $w\in\Sym_n$. Then $s_{r_1}\dots s_{r_l}$ is a \textbf{reduced
expression} for $w$, and $w$ has \textbf{length} $\ell(w)=l$, if $l$ is minimal
such that $w=s_{r_1}\dots s_{r_l}$, for $1\le r_k<n$. Fix a reduced expression
$w=s_{r_1}\dots s_{r_l}$ for $w$ and define $\psi_w = \psi_{r_1}\ldots
\psi_{r_l}$. In view of~\autoref{E:braid}, $\psi_w$ usually depends on the
chosen reduced expression for $w$.

\subsection{Cyclotomic Hecke algebras}
  Recall that $\Lambda\in P^+$ and that we have fixed an integer
  $e\in\{0,2,3,4,\dots\}$ and multicharge $\charge$. We now define the integral
  cyclotomic Hecke algebras $\H$ of type $G(\ell,1,n)$, where $\ell=\sum_{i\in
  I}(\Lambda,\alpha_i)$ is the \textbf{level} of $\Lambda$.

Fix a field~$K$. If $k\in\Z$ and $t\in K$ define the $t$-quantum integer
$$[k]_t=\begin{cases}
  \phantom{-(}1+t+\dots+t^{k-1},&\text{if }k\ge0,\\                                        -(t^{-1}+t^{-2}+\dots+t^{k}),&\text{if }k<0.
\end{cases}
$$

Now fix a non-element $\xi=\xi(e)\in K$ where $e$ is minimal such that
$[e]_\xi=0$, if $e>1$, or $[k]_\xi\ne0$ for all $k\in\N$, if $e=0$.

Following \cite{HuMathas:SeminormalQuiver}, we make the following definition.

\begin{Definition}
\label{D:HeckeAlgebras}
  Suppose that $K$ is a field containing a non-zero element $\xi=\xi(e)$ as
  above and let $\charge\in\Z^\ell$ be a multicharge such that
  $\Lambda=\Lambda_\charge$.  The \textbf{(integral) cyclotomic Hecke algebra} $\H$ of type
  $G(\ell,1,n)$ with Hecke parameter $\xi$ is the unital associative $K$-algebra
  with generators $L_1,\dots,L_n$, $T_1,\dots,T_{n-1}$ which are subject to the
  relations
  $$ \begin{array}{r@{\ }lr@{\ }l@{\quad}l}
      \prod_{l=1}^\ell(L_1-[\kappa_l])&=0,  &
      (T_r+1)(T_r-\xi )&=0,  \\
      L_rL_t&=L_tL_r, &
    T_rT_s&=T_sT_r &\text{if }|r-s|>1,\\
    T_sT_{s+1}T_s&=T_{s+1}T_sT_{s+1}, &
    T_rL_t&=L_tT_r,&\text{if }t\ne r,r+1,\\
    \multicolumn4c{L_{r+1}T_r-T_rL_r=1+(\xi-1)L_{r+1},}
  \end{array}$$
  where $1\le r<n$, $1\le s<n-1$ and $1\le t\le n$.
\end{Definition}

The connection between the cyclotomic quiver Hecke algebras of type $\Gamma_e$
and the cyclotomic Hecke algebras of type~$G(\ell,1,n)$ is given by the
following remarkable result of Brundan and Kleshchev.

\begin{Theorem}%
  [\protect{Brundan-Kleshchev's graded isomorphism theorem~\cite[Theorem~1.1]{BK:GradedKL}}]%
  \label{T:BKiso}
Suppose $K$ is a field.
Then $\H\cong\UnR$.
\end{Theorem}

In particular, over a field~$K$, the cyclotomic quiver Hecke algebra
$\R(K)=\R(\Z)\otimes_\Z K$ has dimension $\ell^n n!$ and, up to isomorphism,
$\H$ depends only on~$e$ and~$\Lambda$. By
\cite{HuMathas:SeminormalQuiver}, $\R(\Z)$ is free of rank~$\ell^nn!$ if $e=0$
or if~$e$ is a prime integer.

\subsection{Partitions and tableaux}
\label{SSPar}
We now introduce the tableau combinatorics which describes the representation
theory of~$\R$.

A \textbf{partition} of~$m$ is a sequence $\mu=(\mu_1,\mu_2,\dots)$ of
non-negative integers such that $m=|\mu|=\mu_1+\mu_2+\ldots$ and $\mu_1\ge\mu_2\ge\dots$. An
\textbf{$\ell$-multipartition} of~$n$ is an
ordered $\ell$-tuple of partitions $\bmu = (\mu^{(1)}|\dots|\mu^{(\ell)})$ such
that $|\mu^{(1)}|+\dots+|\mu^{(\ell)}|=n$. The partition $\mu^{(l)}$ is the
$l$th \textbf{component} of $\bmu$. Let~$\Parts$ be the set of all
$\ell$-multipartitions of~$n$. The set $\Parts$ depends only on $\ell$, and not on
$\Lambda$, however, soon we will introduce the residues of a multipartition
which will depend on~$\Lambda$ or, more accurately, on~$\charge$.
When $\ell$ is clear from context we will talk of
multipartitions.

We identify the multipartition $\bmu\in\Parts$ with its \textbf{diagram}
$$
\set{(l,r,c)\in\{1,\dots,\ell\}\times\N^2| 1\le c\le \mu_r^{(l)}}.
$$
The elements of (the diagram of) $\bmu$ are nodes. More generally, a
\textbf{node} is any element of $\{1,\dots,\ell\}\times\N^2$. We think of~$\bmu$
as an $\ell$-tuple of arrays of boxes in the plane. For example, the diagram of
$(3,1|1^3|4,2)$ is
$$\TriDiagram{3,1}{1,1,1}{4,2}.$$

The \textbf{length} of a multipartition $\bmu$ is
$\ell(\bmu)=\#\set{(l,r)|\mu^{(l)}_r\ne0}$. That is, $\ell(\bmu)$ is the number
of non-empty rows in the diagram of~$\bmu$.
A node $A\in\bmu$ is a \textbf{removable node (of~$\bmu$)}\, if $\bmu\backslash
\{A\}$ is (the diagram of) a multipartition. A node $B\not\in\bmu$ is an
\textbf{addable node (for~$\bmu$)}\, if $\bmu\cup \{B\}$ is a multipartition.
Let $\Add(\bmu)$ and $\Rem(\bmu)$ be the sets of addable and removable nodes of the multipartition~$\bmu$.

Let $\bmu$ and $\bnu$ be multipartitions. Then $\bmu$ \textbf{dominates} $\bnu$,
 or $\bmu\gedom\bnu$, if
$$
\sum_{l=1}^{m-1}|\mu^{(l)}|+\sum_{r=1}^s\mu_r^{(m)}\ge
\sum_{l=1}^{m-1}|\nu^{(l)}|+\sum_{r=1}^s\nu_r^{(m)}
$$
for all $1\le m\le\ell$ and $s\ge 1$. Therefore, $\bmu\gedom\bnu$ if
$\bnu$ is obtained from $\bmu$ by moving nodes down or to the right in the
diagrams. We write $\bmu\gdom\bnu$ if $\bmu\gedom\bnu$ and $\bmu\ne\bnu$.

Let $\bmu$ be a multipartition. A \textbf{$\bmu$-tableau} $\t$ is a bijection
$\t\map\bmu\{1,2,\dots,n\}$. We think of a $\bmu$-tableau $\t$ as a labelling of
the diagram of~$\bmu$ and we set $\Shape(\t)=\bmu$.  In this way we talk of the
\textbf{components}, \textbf{rows} and \textbf{columns} of tableaux. In
particular, if $\t(l,r,c)=k$ then $k$ lies in row $(l,r)$ of $\t$ and we write
$\row_\t(k) = (l,r)$.  The rows of~$\bmu$ (and the rows of $\t$) are
totally ordered lexicographically, from top left to bottom right, by defining
\begin{equation}
\label{E:RowOrder}
  (l,r)<(l',r') \qquad\text{ if $l<l'$ or if $l=l'$ and $r<r'$.}
\end{equation}

The \textbf{residue} of a node $(l,r,c)$ is $\res(l,r,c)=\kappa_l-r+c\pmod e$,
so that $\res(l,r,c)\in I$. An \textbf{$i$-node} is any node of
residue~$i$, for $i\in I$. Similarly, if $\t$ is a tableau and $\t(l,r,c)=k$ then the
\textbf{residue} of~$k$ in~$\t$ is $\res_\t(k)=\res(l,r,c)$.  The
\textbf{residue sequence} of $\t$ is
$$\res(\t) = (\res_\t(1),\res_\t(2),\dots,\res_\t(n))\in I^n.$$
Note that the residues depend upon the choice of multicharge~$\charge$.

A $\bmu$-tableau~$\t$ is \textbf{row standard} if its entries increase along the
rows in each component. A $\bmu$-tableau~$\t$ is \textbf{standard} if it is
row standard and its entries increase down the columns in each component. Let
$\Std(\bmu)$ be the set of standard $\bmu$-tableaux.

If $\t$ is a standard $\bmu$-tableau and $1\le k\le n$ let $\t_{\downarrow
k}$ be the subtableau of~$\t$ which contains $1,2\dots, k$. Then
$\Shape(\t_{\downarrow k})$ is a multipartition for $1\le k\le n$.
If~$\s,\t\in\Std(\bmu)$ then $\s$ \textbf{dominates}~$\t$, and we write
$\s\gedom\t$, if
$$\Shape(\s_{\downarrow k})\gedom\Shape(\t_{\downarrow k}),
       \qquad\text{for }1\le k\le n.$$
This defines a partial order on the set of row standard $\bmu$-tableaux.

Let $\tmu$ be the unique row standard $\bmu$-tableau such that $\tmu\gedom\t$,
for all $\t\in\Std(\bmu)$. Then $\tmu$ has the numbers $1,2,\dots,n$ entered in
order from left to right along its rows. Set $\bi^\bmu = \res(\tmu)$. For any
$\bmu$-tableau~$\t$ let $d(\t)\in\Sym_n$ be the unique permutation such that
$\t=\tmu d(\t)$. It is a well-known result, commonly attributed to Ehresmann and
James, that $\s\gedom\t$ if and only if~$d(\s)\le d(\t)$, where $\le$ is the
Bruhat order on~$\Sym_n$; see, for example,~\cite[Theorem~3.8]{M:Ulect}.

If $\bmu \in \Parts$ and $A$ is a removable $i$-node of~$\bmu$, let
\begin{align*}
d_A(\bmu)&= \#\{\text{addable $i$-nodes of $\bmu$ strictly below } A\} \\
&\hspace*{10mm}-\#\{\text{removable $i$-nodes of $\bmu$ strictly below } A\}.
\end{align*}
For $\t \in \Std(\bmu)$, we define the \textbf{degree} of $\t$ recursively as follows.  If $n=0$ then $\deg \t=0$.  Otherwise let $A=\t^{-1}(n)$ be the node occupied by $n$ in $\t$ and set
\begin{equation}
\label{E:degrees}
\deg \t = \deg \t\backslash \{A\}+d_A(\bmu).
\end{equation}
This degree function describes the grading on the Specht module~$S^\bmu$.

\subsection{Graded Specht modules} \label{S:SpechtModules}
The graded Specht modules were first defined over a field by Brundan, Kleshchev
and Wang~\cite{BKW:GradedSpecht}. In this paper we follow the approach
of~\cite{KMR:UniversalSpecht} which defines the Specht modules by generators and
relations over an arbitrary ring. The main advantage of
using~\cite{KMR:UniversalSpecht} is that the graded Specht modules are defined
over~$\Z$ and the action of~$\R(\Z)$ on these modules is explicitly described in
terms of the KLR generators. In contrast, using the definitions
in~\cite{BKW:GradedSpecht} it is only possible to compute inside the graded
Specht modules by repeatedly applying \autoref{T:BKiso}.

We remark that~\cite{HuMathas:GradedCellular} gives a third construction of the
graded Specht modules using cellular algebra techniques. This approach equips
the graded Specht modules with a homogeneous bilinear form. We do not need this
result here.

Fix a multipartition~$\bmu$ and a node $A=(l,r,c)\in\bmu$.
If $(l,r+1,c) \in \bmu$ then~$A$ is a (row) \textbf{Garnir node} of $\bmu$.
Suppose that $e\ne0$. Then the \textbf{$(e,A)$-Garnir belt} is the set of nodes
\begin{align*}
\Belt =&\set{(l,r,c)\in\bmu|r\ge c\text{ and } e\ceil{r-c+1}/e\leq\mu^{(l)}_r-c+1}\\
      &\qquad\cup\set{(l,r+1,c)\in\bmu|r\le c\text{ and }c\geq e\ceil{c-r+1}/e}.
\end{align*}

Let $b_A=\#\Belt/e$ and write $b_A=a_A+c_A$ where $ea_A$ is the number of nodes
in~$\Belt$ in row $(l,r)$. Let $\D_A$ be the set of minimal length right coset
representatives of~$\Sym_{a_A}\times\Sym_{c_A}$ in~$\Sym_{b_A}$; see, for
example, \cite[Proposition~3.3]{M:Ulect}. If $e=0$ set $\Belt=\emptyset$,
$b_A=0=a_A=c_A$ and $\D_A=1$.

Suppose $A$ is a Garnir node of $\bmu$.  Let $\tA$ be the $\bmu$-tableau which
agrees with~$\tmu$ for all numbers $k<\tmu(A)=\tmu(l,r,c)$ and $k>\tmu(l,r+1,c)$
and where the remaining entries in~rows $(l,r)$ and $(l,r+1)$ are filled in increasing order from
left to right first along the nodes in row~$(l,r+1)$ which are in the first
$c$~columns but not in~$\Belt$, then along the nodes in row~$(l,r)$ of $\Belt$
followed by the nodes in row~$(l,r+1)$ of~$\Belt$, and then along the remaining
nodes in row~$(l,r)$.

\begin{Example}  Suppose that $e=2$, $\bmu=(3,1|7,5|2,1)$ and $A=(1,3,2)$. Then
  $$\tA[2]=\bigg( \hspace*{1mm}\Tableau[-2]{{1,2,3},{4}}\hspace*{1mm}\bigg|\hspace*{1mm}%
  \begin{tikzpicture}[baseline=-2mm,scale=0.5,draw/.append style={thick,black}]
  \newcount\col
  \foreach\Row/\row in {{5,6,8,9,10,11,14}/0,{7,12,13,15,16}/-1} {
     \col=1
     \foreach\k in \Row {
        \draw(\the\col,\row)+(-.5,-.5)rectangle++(.5,.5);
        \draw(\the\col,\row)node{\k};
        \global\advance\col by 1
      }
   }
   \draw[red,double,very thick](2.5,-0.5)--(2.5,0.5)--(6.5,0.5)--(6.5,-0.5)--(1.5,-0.5)
                     --(1.5,-1.5)--(3.5,-1.5)--(3.5,-0.5);
   \draw[red,double,very thick](4.5,0.5)--(4.5,-0.5);
  \end{tikzpicture}%
  \hspace*{1mm}\bigg|\hspace*{1mm}\Tableau[-2]{{17,18},{19}}\hspace*{1mm}\bigg)$$
  The double lines in~$\tA[2]$ show the $(2,A)$-Garnir belt of~$\bmu$ and they
  show how it decomposes into a disjoint union of ``$e$-bricks''.
  In general, if $e\ne 0$ then~$b_A$ is the number of $e$-bricks in the Garnir
  belt and that~$a_A$ is the number of $e$-bricks in its first row. In this
  case,~$b_A=3$ and~$a_A=2$. Therefore, $\D_A=\{1,s_2,s_2s_1\}$.
\end{Example}

Let $k_A=\tA(A)$ be the number occupying~$A$ in~$\tA$. For $1\le r<b_A$ define
$$w^A_r=\prod_{a=k_A+e(r-1)}^{k_A+re-1}(a,a+e).$$
Although it does not affect  the last definition, because all of the factors commute,
in this paper we use the convention will that all products are read \textbf{from left to
right}. The set $\set{w^A_r|1\le r<b_A}$ generates a subgroup of~$\Sym_n$
isomorphic to $\Sym_{b_A}$, where the isomorphism is determined by the map
$w^A_r\mapsto s_r$, for $1\le r<b_A$. Set $\bi_A=\res(\tA)$. If $d\in\D_A$
choose a reduced expression $d=s_{r_1}\dots s_{r_k}$ for $d$ and define
$$\tau^A_d = e(\bi_A)(\psi_{w^A_{r_1}}+1)\dots(\psi_{w^A_{r_k}}+1)\in\R.$$
Then $\tau_d^A$ is independent of the choice of reduced expression for~$d$
and independent of the choice of reduced expression for each $\psi_{w^A_{r_i}}$ by
\cite[Theorem~5.11]{KMR:UniversalSpecht}. We can now define the Specht modules
of~$\R$.

\begin{Definition}[\protect{\cite[Definition~5.9]{KMR:UniversalSpecht}}]\label{D:SpechtModule}
  Suppose that $\bmu\in\Parts$. The \textbf{Specht module} $S^\bmu$ of~$\R$
  is the $\R$-module generated by the homogeneous element $v_{\tmu}$ of
  degree~$\deg\tmu$ subject to the relations:
  \begin{enumerate}
    \item $v_{\tmu} e(\bi)=\delta_{\bi\bi^\bmu}v_{\tmu}$.
    \item $v_{\tmu} y_s=0$, for $1\le s\le n$.
    \item $v_{\tmu} \psi_r=0$ whenever $r$ and $r+1$ are in the same row of $\tmu$,
    for $1\le r<n$.
    \item $\sum_{d\in\D_A}v_{\tmu}\psi_{\tA}\tau^A_d=0$, for all Garnir nodes $A\in\bmu$.
  \end{enumerate}
\end{Definition}

The relations in part~(d) are the \textbf{homogeneous Garnir relations}. By
\cite[Theorem~6.23]{KMR:UniversalSpecht}, if~$K$ is a field then the graded Specht
module $S^\bmu\otimes_\Z K$ is isomorphic, as a graded $\R$-module, to the
graded Specht module studied previously by Brundan, Kleshchev and
Wang~\cite{BKW:GradedSpecht} and by Hu and Mathas~\cite{HuMathas:GradedCellular}.

For any  $\bmu$-tableau $\t$ define $v_\t=v_{\tmu}\psi_{d(\t)}\in S^\bmu$. These elements
depend on the choice of reduced expression for the permutation $d(\t)$. Nonetheless, by fixing arbitrary choices of these reduced expressions we have the following:

\begin{Theorem}[\!\protect{\cite[Theorem 4.2]{BKW:GradedSpecht},
                        \cite[Theorem~6.21]{KMR:UniversalSpecht}}]
Suppose that $\bmu\in\Parts$. Then $\set{v_\t|\t\in\Std(\bmu)}$ is a homogeneous
basis of~$S^\bmu$ with $\deg v_\t=\deg\t$, for all $\t\in\Std(\bmu)$.
\end{Theorem}

For future use we now state three properties of the homogeneous basis of
$S^\bmu$ which go back to results of Brundan, Kleshchev and
Wang~\cite{BKW:GradedSpecht}. In view of \cite[Theorem
6.23]{KMR:UniversalSpecht} these results are valid over~$\Z$. First observe that
\autoref{D:SpechtModule}~(a) and the relation $\psi_re(\bi)=e(s_r\bi)\psi_r$
implies that if $\t$ is any $\bmu$-tableau then
\begin{equation}
\label{E:tableauxResidue}
  v_\t e(\bi) = \delta_{\bi,\res(\t)} v_\t,
\end{equation}
for all $\bi\in I^n$.

\begin{Lemma}[\protect{\cite[Corollary 4.6 \& Proposition 4.7]{BKW:GradedSpecht}}] \label{L:DomTableaux}
Let $s_{i_1}s_{i_2}\ldots s_{i_l}$ be a reduced expression for a permutation $w \in \Sym_n$ such that $\t = \tmu w$ is standard.  Then
\[v_{\tmu} \psi_{i_1}\psi_{i_2}\ldots \psi_{i_l} = \sum_{\substack{\s \in \Std(\bmu), \\ \s \gedom \t}} a_\s v_\s\]
for some $a_\s \in \Z$.   Furthermore $a_\s \neq 0$ only if $\res(\s)=\res(\t)$.
\end{Lemma}

\begin{Lemma}[~{\cite[Lemma 4.8]{BKW:GradedSpecht}}]
\label{L:YDown}
Suppose that $\t \in \Std(\bnu)$ and $1 \leq r \leq n$.  Then
\[v_\t y_r = \sum_{\substack{\s \in \Std(\bmu), \\ \s \gdom \t}} a_\s v_\s\]
for some $a_\s \in \Z$.  Furthermore if $a_\s \neq 0$ then $\res(\s)=\res(\t)$.
 \end{Lemma}

\subsection{Braid diagrams}
\label{S:braid} This section introduces the diagram
combinatorics of Khovanov and Lauda~\cite[\S2]{KhovLaud:diagI}, adapting it to
our situation. In principle, all of the arguments in this paper could be given
algebraically, however, they become much more transparent, and shorter, when we
use these diagrams.

An \textbf{$n$-braid diagram} is an $I$-labelled graph on
$\{1',\dots,n',1,\dots,n\}$ such that each vertex in $\{1',\dots,n'\}$ is
connected to a unique vertex in $\{1,\dots,n\}$ and each edge is labelled by an
element of~$I$ and is decorated by a finite, non-negative, number of dots.

As we explain below,
the braid diagrams in this paper represent an element of a Specht module
$S^\bmu$, for some multipartition~$\bmu\in\Parts$. All of the
relations in~$\R$ and in~$S^\bmu$ have diagrammatic counterparts which describe
how to rewrite a given diagram as a linear combination of other diagrams.
Considered as a map from $\{1',\dots,n'\}$ to $\{1,\dots,n\}$, each braid
diagram~$B$ determines a unique permutation $\pi_B\in\Sym_n$. We emphasize that
the element of~$S^\bmu$ corresponding to~$B$ is, in general, not determined by
$\pi_B$ and is dependent on the crossings and the dots in~$B$.

We represent braid diagrams as graphs in the plane and we do not distinguish
between under and over crossings. Braid diagrams are drawn with the vertices
labelled by $\{1',\dots,n'\}$ at the top of the diagram and the vertices
labelled by $\{1,\dots,n\}$ at the bottom, with both sets of vertices ordered
from left to right in the obvious way. For $1\le m\le n$ the (unique) edge with
vertex~$m$ is the \textbf{$m$-string} of \textbf{residue}~$i$, where $i\in I$ is
the corresponding edge label. We omit the labels for the vertices
$\{1',\dots,n'\}$ and put the edge labels $\bi=(i_1,\dots,i_n)\in I^n$ at the
top of each string. To improve readability, we colour the strings of the
diagrams, however, these colours have no mathematical meaning (and may not be
distinguishable in black and white!).

Writing $\bi^\bmu=(i_1,\dots,i_n)$, the braid diagram representing $v_{\tmu}\in S^\bmu$ is:
$$v_{\tmu} =
\begin{braid}\tikzset{baseline=8mm}
  \draw(1,4)node[above]{$i_1$}--(1,0)node[below]{$1$};
    \draw(2,4)node[above]{$i_2$}--(2,0)node[below]{$2$};
    \draw[dots](2.5,4)--(4.5,4);
    \draw[dots](2.5,0)--(4.5,0);
    \draw(5,4)node[above]{$i_r$}--(5,0)node[below]{$r$};
    \draw[dots](5.5,4)--(7.5,4);
    \draw[dots](5.5,0)--(7.5,0);
    \draw(8,4)node[above]{$i_n$}--(8,0)node[below]{$n$};
  \end{braid}
  =
  \begin{braid}\tikzset{baseline=9mm}
    \draw(1,4)node[above]{$i_1$}--(1,0)node[below]{$1$};
    \draw[dots](1.5,0)--(2.5,0);
    \draw(3,4)node[above]{$i_a$}--(3,0)node[below]{$a$};
    \draw[dots](3.7,4)--(5.3,4);
    \draw[dots](3.5,0)--(5.5,0);
    \draw[brown](0.5,5)--(3.5,5)--(3.5,4)--(0.5,4)--cycle;
    \draw(6,4)node[above]{$i_{m+1}$}--(6,0)node[below]{$m+1$};
    \draw[dots](6.5,0)--(7.5,0);
    \draw(8,4)node[above]{$i_n$}--(8,0)node[below]{$n$};
    \draw[brown](5.0,5)--(8.5,5)--(8.5,4)--(5.0,4)--cycle;
    \draw(2,4.8)node[above]{$\bmu^{(1)}_1$};
    \draw(7,4.8)node[above]{$\bmu^{(\ell)}_z$};
  \end{braid}
$$
As in the right-hand diagram (where $a=\mu^{(1)}_1$, $m=n-\mu^{(\ell)}_z<n$ and
$\mu^{(\ell)}_{z+1}=0$), we will sometimes (but not always) indicate the rows
of~$\bmu$ by drawing a circle around the corresponding residues at the top of
the diagram.

For $1\le r<n$, $\psi_r$ acts on a braid diagram by crossing the $r$-string
and the $(r+1)$-string at the bottom of the diagram.  If $1\le s\le n$ then $y_s$ acts by putting a dot at the
bottom of the
$s$-string. If $\bj\in I^n$ then $e(\bj)$ acts as the identity on a diagram $B$ if
$\bj = \pi_B^{-1} \cdot \bi$ , and otherwise it acts as
zero.

All of the relations in \autoref{D:KLR} and \autoref{D:SpechtModule} have
corresponding diagrammatic versions. For example, the relations in parts~(b) and~(c) of \autoref{D:SpechtModule} become:
\begin{equation}
\label{E:InsideRow}
  \begin{braid}[7]
    \draw(1,4)node[above]{$i_1$}--(1,0)node[below]{$1$};
    \draw(2,4)node[above]{$i_2$}--(2,0)node[below]{$2$};
    \draw[dots](2.5,4)--(4.5,4);
    \draw[dots](2.5,0)--(4.5,0);
    \draw(5,4)node[above]{$i_s$}--(5,0)node[below]{$s$};
    \greendot(5,2);
    \draw[dots](5.5,4)--(7.5,4);
    \draw[dots](5.5,0)--(7.5,0);
    \draw(8,4)node[above]{$i_n$}--(8,0)node[below]{$n$};
  \end{braid}=0
  \quad\text{and}\quad
  \begin{braid}[7]
    \draw(1,4)node[above]{$i_1$}--(1,0)node[below]{$1$};
    \draw[dots](1.5,0)--(2.5,0);
    \draw(3,4)node[above]{$i_a$}--(3,0)node[below]{$a$};
    \draw[dots](3.7,4)--(4.7,4);
    \draw[dots](3.5,0)--(5.5,0);
    \draw[brown](0.5,5)--(3.5,5)--(3.5,4)--(0.5,4)--cycle;
    \draw(6,4)node[above]{$i_r$}--(6.5,0);
    \draw(7,4.5)node{$i_{r+1}$};
    \draw(7,-0.4)node{$r+1$};
    \draw(6.5,4)--(6,0)node[below]{$r$};
    \draw[brown](4.7,5)--(8.3,5)--(8.3,4)--(4.7,4)--cycle;
    \draw[dots](8.8,4)--(9.9,4);
    \draw[dots](8.3,0)--(9.9,0);
    \draw(10.5,4)node[above]{$i_m$}--(10.5,0)node[below]{$m$};
    \draw(12,4)node[above]{$i_n$}--(12,0)node[below]{$n$};
    \draw[brown](10,5)--(12.5,5)--(12.5,4)--(10,4)--cycle;
  \end{braid}=0
  \end{equation}
  respectively. Similarly when $e \neq 2$, relations~\autoref{E:ypsiCommute}
  and~\autoref{E:braid} yield the diagrams:
  \begin{equation}
\label{E:iiCrossing}
    \begin{braid}[7]
    \draw(1,4) node[above]{$i$}--(3,0);
    \draw(3,4) node[above]{$j$}--(1,0);
    \greendot(2.5,1);
  \end{braid}
  =
  \begin{braid}[7]
    \draw(1,4) node[above]{$i$}--(3,0);
    \draw(3,4) node[above]{$j$}--(1,0);
    \greendot(1.5,3);
  \end{braid}
  +\delta_{ij}\begin{braid}[7]
    \YCrossing i 1 3 1 3
  \end{braid}
  \quad\text{and}\quad
  \begin{braid}[7]
    \draw(1,4)node[above]{$i$}--(3,0);
    \draw(2,4)node[above]{$i{\pm}1$}--(1,2)--(2,0);
    \draw(3,4)node[above]{$i$}--(1,0);
  \end{braid}
  =
  \begin{braid}[7]
    \draw(1,4)node[above]{$i$}--(3,0);
    \draw[red](2,4)node[above]{$i{\pm}1$}--(3,2)--(2,0);
    \draw(3,4)node[above]{$i$}--(1,0);
  \end{braid}
  \pm
  \begin{braid}[7]
    \draw(1,4)node[above]{$i$}--(1,0);
    \draw(2,4)node[above]{$i{\pm}1$}--(2,0);
    \draw(3,4)node[above]{$i$}--(3,0);
  \end{braid}
  \end{equation}
  both of which should be viewed as `local' relations for strings inside a
  braid diagram. The homogeneous Garnir relations of
  \autoref{D:SpechtModule}~(d) are harder to represent diagrammatically, however,
  in the special case when $A=(l,r,c)$ and either $e=0$ or $\delta=\bmu^{(l)}_r-c+1<e$
  this relation simplifies to $v_{\tmu}\psi_{\tA}=0$, or
  \begin{equation} \label{E:SmallGarnir}
  \begin{braid}[7]
    \draw[dots](-.5,4.5)--(0.5,4.5);\draw[dots](10,4.5)--(11,4.5);
    \draw[dots](-.5,0.1)--(0.5,0.1);\draw[dots](10,0.1)--(11,0.1);
    \draw[brown](1,5)--(5.4,5)--(5.4,4)--(1,4)--cycle;
    \draw(3.5,4.8)node[above]{$\bmu^{(l)}_r$};
    \draw(7.5,4.8)node[above]{$\bmu^{(l)}_{r+1}$};
    \draw(3,4)--(3,0);
    \draw[dots](1.7,0.1)--(3.0,0.1);
    \draw(3.5,4)--(3.5,0);
    \draw(1.5,4)--(1.5,0);
    \draw(8.5,4)--(8.5,0);
    \draw(9,4)--(9,0);
    \draw(4,4)--(7,0);
    \draw(5,4)--(8,0);
    \draw(7.5,-0.4)node{$\underbrace{\hspace*{6mm}}_{\delta<e}$};
    \draw(6,-1.4)node{$\underbrace{\hspace*{18mm}}_{\mu^{(l)}_r+1}$};
   \draw(6,4)--(4,0);
    \draw[dots](4.3,.1)--(5.9,0.1);
    \draw(8,4)--(6,0);
    \draw[brown](5.6,5)--(9.5,5)--(9.5,4)--(5.6,4)--cycle;
  \end{braid}=0.
\end{equation}

We will use these diagrams extensively in this paper. The reader who has not
used these diagrams before may find, initially, that they are cumbersome to work
with. The reader should persevere, however, because these diagrams allow us to
replace long calculations with the KLR relations with one or two diagrams which
make the arguments more transparent. As these diagrams take some getting use to
we give a toy example of how they are used which also illustrates why they are
useful.

\begin{Example}
\label{Ex:diagrams}
Take $e=2$, let~$\mu$ be the partition $(5,4)$ and consider
$v_\t=v_{\tmu}\psi_5\psi_6\psi_7\psi_8y_9$, where
$\t=\tableau{{1,2,3,4,9},{5,6,7,8}}$. Let $i=\res_{\tmu}(5)=0$.
This product corresponds to:
$$v_\t=v_{\tmu}\psi_5\psi_6\psi_7\psi_8y_9
   =\begin{braid}[8]\tikzset{xscale=1.2}
      \draw[brown](2,5)--(5.18,5)--(5.18,4)--(2,4)--cycle;
      \draw[brown](5.32,5)--(7.3,5)--(7.3,4)--(5.32,4)--cycle;
      \foreach \k in {2.5,3,...,4.5} {
         \draw(\k,4)--(\k,0);
      }
      \foreach \k in {5.5,6,6.5,7}  {
          \draw(\k,4)-- +(-0.5,-4);
      }
      \draw(5,4)node[above]{$i$}--(7,0);
      \greendot(6.85,0.3);
    \end{braid}.
  $$
The dot at the bottom of the diagram represents~$y_9$ because it is on the
$9$-string, which is the ninth string reading from left to right at the bottom
of the diagram. The $9$-string has residue~$0$ because its top vertex is $5'$,
the fifth string reading left to right at the top of the diagram, and
$\res_{\tmu}(m)=i$. The $9$-string crosses the four
$k$-strings, for $k=5,6,7,8$. Reading these crossing in order, from top to
bottom, corresponds to multiplying by $\psi_5,\dots,\psi_8$, respectively.

Now it is, of course, easy to compute this example using the relations.
Diagrammatically, it is even easier because \autoref{E:iiCrossing} says that we
slide the dot all the way up to the top of the $9$-string, where it comes zero,
except that whenever we cross another string of residue~$i$ then we must add the
a new diagram created by removing the dot and \textit{cutting} the
\iicrossing-crossing. The $8$-string and the $6$-string have residue~$i$, so
$$v_\t y_9=\begin{braid}[8]\tikzset{xscale=1.2}
      \draw[brown](2,5)--(5.18,5)--(5.18,4)--(2,4)--cycle;
      \draw[brown](5.32,5)--(7.3,5)--(7.3,4)--(5.32,4)--cycle;
      \foreach \k in {2.5,3,...,4.5} {
         \draw(\k,4)--(\k,0);
      }
      \foreach \k in {5.5,6,6.5,7}  {
          \draw(\k,4)-- +(-0.5,-4);
      }
      \draw(5,4)node[above]{$i$}--(7,0);
      \greendot(6.85,0.3);
    \end{braid}
    =\begin{braid}[8]\tikzset{xscale=1.2}
      \draw[brown](2,5)--(5.18,5)--(5.18,4)--(2,4)--cycle;
      \draw[brown](5.32,5)--(7.3,5)--(7.3,4)--(5.32,4)--cycle;
      \foreach \k in {2.5,3,...,4.5} {
         \draw(\k,4)--(\k,0);
      }
      \foreach \k in {5.5,6,6.5}  {
          \draw(\k,4)-- +(-0.5,-4);
      }
      \draw[red](5,4)node[above]{$i$}--($ (5,4)!0.75!(7,0) $)--(6.5,0);
      \draw[red](7,4)node[above]{$i$}--($ (7,4)!0.65!(6.5,0) $)--(7,0);
    \end{braid}
    +
    \begin{braid}[8]\tikzset{xscale=1.2}
      \draw[brown](2,5)--(5.18,5)--(5.18,4)--(2,4)--cycle;
      \draw[brown](5.32,5)--(7.3,5)--(7.3,4)--(5.32,4)--cycle;
      \foreach \k in {2.5,3,...,4.5} {
         \draw(\k,4)--(\k,0);
      }
      \foreach \k in {5.5,6.5,7}  {
          \draw(\k,4)-- +(-0.5,-4);
      }
      \draw[red](5,4)node[above]{$i$}--($ (5,4)!0.35!(7,0) $)--(5.5,0);
      \draw[red](6,4)node[above]{$i$}--($ (6,4)!0.35!(5.5,0) $)--(7,0);
    \end{braid}.
  $$
  The second diagram is zero by \autoref{E:InsideRow}, so
  $v_\t y_9=v_\s$, where
  $\s=\tableau{{1,2,3,4,8},{5,6,7,9}}$.

  If the reader if not yet convinced that these diagrams simplify these
  calculations, then exactly the same argument shows that if $e \in \{0,2,3,\ldots\}$ and
  $\mu=(m,n-m)$, where $m\le n\le 2m$, then
 $$v_{\tmu} \psi_{m} \psi_{m+1} \ldots \psi_{n-1} y_n =
\begin{cases} v_\s, & \res_{\tmu}(m)=\res_{\tmu}(n), \\
 0, & \text{otherwise},
 \end{cases} $$
  where $\s=\tmu s_m\dots s_{n-2}$.
\end{Example}

\subsection{Restricting Specht modules}
\label{S:Restriction}
We now describe the framework that we use to prove our Main Theorem. These
ideas are a natural extension of work of Ellers and
Murray~\cite{EllersMurray:branching}.

Let $\bnu$ be a multipartition of $n+1$. By~\cite[Theorem~4.11]{BKW:GradedSpecht} and
\cite[Theorem~6.23]{KMR:UniversalSpecht}, the restriction $\Res S^\bnu$
of~$S^\bnu$ to $\R[n]$ has a graded Specht module filtration
\begin{equation}
\label{E:SpechtFiltration}
     S^\bnu =S^\bnu_0\supset S^\bnu_1\supset\dots \supset S^\bnu_{z+1}=0
\end{equation}
where $S^\bnu_k/S^\bnu_{k+1}\cong S^{\bnu\backslash\{A_{k+1}\}}\<\deg A_{k+1}\>$, for $0\le k\leq z$,
and $A_1<\dots<A_{z+1}$ are the removable nodes of~$\bnu$, where the nodes $A_1,\dots,A_{z+1}$
are ordered lexicographically so that $\row(A_1)<\dots<\row(A_{z+1})$ according
to~\autoref{E:RowOrder}. In particular, $S^\bnu_k$ is the $\R$-submodule of
$S^\bnu$ with basis $\set{v_\t|\t\in\Std(\bnu)\text{ and }\Shape(\t_{\downarrow
n})\gedom\bnu\backslash\{A_{k+1}\}}.$

More generally, suppose that $\bnu \in \Parts[n+\gamma]$ where $\gamma >0$.  The next result then follows by considering the restriction of $S^\bnu$ to $\R$ and repeatedly applying~\autoref{E:SpechtFiltration}.

\begin{Lemma} \label{L:Filtration}
Suppose that $\t \in \Std(\bnu)$ and that $h \in \R$.  Then
$$v_\t h = \sum_{\s \in \Std(\bnu)} a_\s v_\s,$$
where $a_\s \neq 0$ only if $\row_\s(n+g)\geq \row_\t(n+g)$ for $1 \leq g \leq \gamma$.
\end{Lemma}

If $\t \in \Std(\bnu)$, let
$\Std_\t^\gamma(\bnu) =\set{\s\in\Std(\bnu)|\s^{-1}(n+g)=\t^{-1}(n+g)\text{ for }1\le g\le\gamma}$.
Now suppose that $\bsig \in \Parts$ and that $\bsig \subset \bnu$.
Let
\begin{align*}
\StdCheck &=\set{\s\in\Std(\bnu)|\Shape(\s_{\downarrow n}) \gdom \bsig}, \\
\StdHat & =\set{\s\in\Std(\bnu)|\bsig \not\gedom\Shape(\s_{\downarrow n})}.
\end{align*}
Then $\StdCheck\subseteq\StdHat$
where this inclusion is usually strict if $\gamma>1$.
If~$\t_{\bsig} \in \Std(\bnu)$ and $(\t_{\bsig})_{\downarrow n}=\t^{\bsig}$
then define
\begin{align*}
  \check{S}^\bnu_{\t_{\bsig}} &= \<v_\s \mid \s\in\Std_{\t_{\bsig}}^\gamma(\bnu)\cup\StdCheck[\bsig]\>_\Z, &
\check{S}^\bnu_\bsig &= \<v_\s \mid \s\in\StdCheck[\bsig]\>_\Z, \\
\hat{S}^{\bnu}_{\t_{\bsig}} &= \<v_\s \mid \s\in\Std_{\t_{\bsig}}^\gamma(\bnu)\cup\StdHat[\bsig]\>_\Z, &
\hat{S}^{\bnu}_{\bsig} &= \<v_\s \mid \s\in\StdHat[\bsig]\>_\Z.
\end{align*}

Fix multipartitions $\blam$ and $\bmu$ of~$n$ such that $\bmu\gdom\blam$ and
$\blam,\bmu\subset\bnu$ and fix standard tableaux $\t_{\blam},\t_{\bmu}\in\Std(\bnu)$
such that $(\t_{\blam})_{\downarrow n}=\tlam$, $(\t_{\bmu})_{\downarrow n}=\tmu$
and $\t_{\bmu}\gdom\t_{\blam}$.
It follows from~\autoref{L:Filtration} that $\check{S}^\bnu_{\t_\blam}$,
$\check{S}^\bnu_{\blam}$, $\hat{S}^{\bnu}_{\t_\bmu}$ and $\hat{S}^{\bnu}_{\bmu}$ are all
$\R$-submodules of $S^\bnu$ and, moreover, that
$S^\bnu$ has an $\R$-module filtration
\begin{equation}
\label{E:Snulam}
  S^\bnu\supseteq \check{S}^\bnu_{\t_{\blam}}\supset \check{S}^\bnu_\blam
  \supseteq \hat{S}^{\bnu}_{\t_{\bmu}}\supset \hat{S}^{\bnu}_\bmu\supset0.
\end{equation}
Furthermore,
$S^\blam\<\deg\t_\blam-\deg\t^\blam\>\cong \check{S}^\bnu_{\t_\blam}/\check{S}^\bnu_\blam$ and
$S^\bmu\<\deg\t_\bmu-\deg\t^\bmu\>\cong \hat{S}^{\bnu}_{\t_\bmu}/\hat{S}^{\bnu}_\bmu$ as
graded $\R$-modules where the isomorphisms are given by
$v_\s+\check{S}^\bnu_\blam\mapsto v_{\s_{\downarrow n}}$ and
$v_\t+\hat{S}^{\bnu}_\bmu\mapsto v_{\t_{\downarrow n}}$, for
$\s\in\Std_{\t_\blam}^\gamma(\bnu)$ and $\t\in\Std_{\t_\bmu}(\bnu)$,
respectively. In particular, as $\R$-modules, $\check{S}^\bnu_{\t_\blam}/\check{S}^\bnu_\blam$
is generated by $v_{\t_\blam}+\check{S}^\bnu_\blam$ and, similarly,
$\hat{S}^{\bnu}_{\t_\bmu}/\hat{S}^{\bnu}_\bmu$ is generated by
$v_{\t_\bmu}+\hat{S}^{\bnu}_\bmu$.

The following result is a graded analogue of \cite[Corollary~2.3]{LM:CarterPayne}.

\begin{Lemma}
\label{L:EllersMurray}
  For $\bnu\in\Parts[n+\gamma]$ and $\blam,\bmu\in\Parts$  with $\blam,\bmu\subset\bnu$, let $\t_{\blam},\t_{\bmu}\in\Std(\bnu)$
  be standard tableaux such that $(\t_{\blam})_{\downarrow n}=\tlam$,
  $(\t_{\bmu})_{\downarrow n}=\tmu$ and $\t_{\bmu}\gdom\t_{\blam}$. Suppose that there
  exists a homogeneous element $L\in\R[n+\gamma]$ such that $L$ commutes with $\R$ and
 $$\check{S}^\bnu_{\t_{\blam}}L\subseteq \hat{S}^{\bnu}_{\t_{\bmu}}, \quad
   \check{S}^\bnu_\blam L\subseteq \hat{S}^{\bnu}_\bmu\quad\text{and}\quad
   v_{\t_{\blam}} L\notin \hat{S}^{\bnu}_\bmu.$$
 Then $\Hom_{\R}(S^\blam\<d\>,S^\bmu)\ne0$, where
 $d=\deg\t_{\blam}-\deg\tlam-\deg\t_{\bmu}+\deg\tmu+\deg L$.
\end{Lemma}

\begin{proof}Since $L$ commutes with $\R$, our assumptions
  ensure that there is a well-defined, non-zero, degree preserving
  $\R$-module homomorphism such that
  $$\theta_L\map{\check{S}^\bnu_{\t_{\blam}}/\check{S}^{\bnu}_\blam\<\deg L\>}{\hat{S}^{\bnu}_{\t_{\bmu}}/\hat{S}^{\bnu}_\bmu}; \qquad
                     x+\check{S}^{\bnu}_\blam\mapsto xL+\hat{S}^{\bnu}_\bmu,$$
  for $x\in \check{S}^\bnu_{\t_{\blam}}$. By construction,
  $\check{S}^\bnu_{\t_{\blam}}/\check{S}^\bnu_\blam\cong S^\blam\<\deg\t_{\blam}-\deg\tlam\>$
  and
  $\hat{S}^{\bnu}_{\tmu}/\hat{S}^{\bnu}_\bmu\cong S^\bmu\<\deg\t_{\bmu}-\deg\tmu\>$ as
  $\R$-modules, so the result follows.
\end{proof}

\section{Cyclotomic Carter-Payne pairs}
In this section we prove our Main Theorem by using the framework of
\autoref{L:EllersMurray} to explicitly construct a non-zero homomorphism between
two Specht modules $S^\blam$ and $S^\bmu$ when $(\blam,\bmu)$ is a cyclotomic
Carter-Payne pair. The first section treats the simplest case, when $\blam$ and $\bmu$
differ by one node. This case is easy to understand and it illustrates
some of the main ideas of the general case. Doing the one node case first also
has the advantage that it allows us to assume that $e\ne2$ in later sections.
Once the one node case is settled we introduce the machinery that we need to
tackle the general case and then prove our Main Theorem.

\subsection{One node cyclotomic Carter-Payne homomorphisms}
\label{S:OneNode} We
are now ready to prove our main theorem in the special case when $\blam$ and
$\bmu$ differ by one node.

Throughout this section, we fix $\bnu \in \Parts[n+1]$ and $i \in I$. If
$\bsig\in\Parts$ write $\bsig\iarrow\bnu$ if $\bsig\subset\bnu$ and
$\res(\bnu\backslash\bsig)=i$, that is, $\bnu\backslash\bsig$ is a removable
$i$-node of~$\bnu$.

\begin{Theorem}
\label{T:OneNode}
Suppose that $\blam, \bmu \iarrow \bnu$ with $\bmu \gdom \blam$. Then
$$\Hom_{\R}(S^\blam\<a+b\>,S^\bmu)\ne0,$$
where $a=\#\set{\alpha\in\Add(\bnu)|\res(\alpha)=i\text{ and }
         \row(\bnu\backslash\blam)<\row(\alpha) \leq \row(\bnu\backslash\bmu)}$
and $b=\#\set{\alpha\in\Rem(\bnu)|\res(\alpha)=i\text{ and }
        \row(\bnu\backslash\blam)<\row(\alpha) \leq \row(\bnu\backslash\bmu)}$.
In particular, $a+b>0$.
\end{Theorem}

It is possible, with some persistence, to prove this result algebraically, however, as a limbering
up exercise for the next sections we prove \autoref{T:OneNode} using the braid
diagrams introduced in \autoref{S:braid}.

If $\bsig \iarrow \bnu$, let $\tnusig$ be the unique standard
$\bnu$-tableau such that $(\tnusig)_{\downarrow n}=\t^\bsig$.  Then
\[v_{\tnusig} = v_{\tnu} \psi_m \psi_{m+1} \ldots \psi_n\]
where $m$ is in the same position of $\tnu$ as $n+1$ is in $\tnusig$.  If $\t$
is any other tableau with $\Shape(\t_{\downarrow n})=\Shape(\tnusig_{\downarrow
n})$ then there exists $w \in \Sym_n$ such that $\tnusig w = \t$ and we define
$v_\t \in S^\bnu$ by $v_\t =v_{\tnusig}\psi_w$. Note that $\psi_w$ may depend on
the chosen reduced expression for $w$.

Suppose that $\bsig,\btau \iarrow\bnu$ and $\btau\gdom\bsig$.
Define $w_{(\btau,\bsig)}\in\Sym_n$ to be the permutation $(m,m-1,\dots,l)$, where $m$ is
the number in position~$\bnu\backslash\btau$ of~$\tnusig$ and~$l$ is the number in position
$\bnu\backslash\bsig$ of~$\tnutau$.

\begin{Lemma}
\label{L:OneNodePush}
  Suppose that $\bsig\iarrow\bnu$. Then
  $$v_{\tnusig}y_{n+1}
      =\sum_{\substack{\btau\iarrow\bnu\\\btau\gdom\bsig}}v_{\tnu_{\btau}} \psi_{w_{(\btau,\bsig)}}.$$
\end{Lemma}

\begin{proof}We want to compute $v_{\tnusig}y_{n+1}$ using diagrams. By design
  this calculation is almost the same as \autoref{Ex:diagrams}. Nonetheless, we
  give a detailed explanation of the notation to ensure that the reader
  understands the diagrams because this language is indispensable later. In
  terms of diagrams we want to compute:
  $$
    v_{\tnusig}y_{n+1}
    =\begin{braid}
      \draw[dots](0,4.5)--(1.5,4.5);\draw[dots](9.8,4.5)--(11,4.5);
      \draw[dots](0,0)--(1.5,0);\draw[dots](8.8,0)--(10,0);
      \draw[brown](2,5)--(5.4,5)--(5.4,4)--(2,4)--cycle;
      \draw[brown](5.6,5)--(9.4,5)--(9.4,4)--(5.6,4)--cycle;
      \draw[brown](11.6,5)--(14.4,5)--(14.4,4)--(11.6,4)--cycle;
      \foreach \k in {2.5,3,...,4.5} {
         \draw(\k,4)--(\k,0);
      }
      \draw(5,4)node[above]{$i$}--(14,0);
      \foreach \k in {6,6.5,...,9}  {
          \draw(\k,4)-- +(-1,-4);
      }
      \foreach \k in {12,12.5,...,14}  {
          \draw(\k,4)-- +(-1,-4);
      }
      \greendot(13.5,0.2);
    \end{braid}
  $$
  The dot represents $y_{n+1}$ because it is on the $(n+1)$-string, which is the
  $(n+1)$th string when we read from left to right along the bottom of the diagram.
  Reading from top to bottom, the $(n+1)$-string crosses the
  $k$-string for $k=m,\dots,n$. Each crossing corresponds to multiplying by
  some~$\psi_r$, so $v_{\tnu}\psi_m\dots\psi_n y_{n+1}=v_{\tnusig}y_{n+1}$ as
  claimed.

  To simplify the diagram above we use~\autoref{E:ypsiCommute} to move the dot past
  the \ijcrossing[j]-crossings above it --- here $i$ and $j$ are the residues of
  the corresponding strings in the diagram. If $j\ne i$ then the dot jumps over
  the crossing without penalty. When $j=i$, however, we have an
  \iicrossing-crossing and by \autoref{E:iiCrossing}
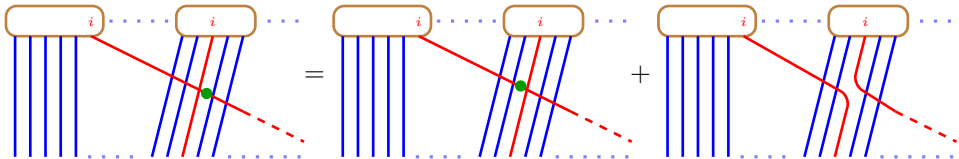
  $$
  \begin{braid}
      \draw[dots](4.6,4.5)--(6.8,4.5);
      \draw[dots](9.8,4.5)--(11,4.5);
      \draw[dots](3.9,0)--(5.7,0);
      \draw[dots](8.5,0)--(11,0);
      \draw[brown](1.2,5)--(4.4,5)--(4.4,4)--(1.2,4)--cycle;
      \draw[brown](6.8,5)--(9.4,5)--(9.4,4)--(6.8,4)--cycle;
      \foreach \k in {1.5,2,...,3.5} {
         \draw(\k,4)--(\k,0);
      }
      \draw(7,4)--(6,0);
      \draw(7.5,4)--(6.5,0);
      \draw(8.5,4)--(7.5,0);
      \draw(9,4)--(8,0);
      \draw[red](4,4)node[above]{$i$}--(9,1.5);\draw[red,dashed](9,1.5)--(11,0.5);
      \draw[red](8,4)node[above]{$i$}--(7,0);
      \greendot(7.8,2.1);
    \end{braid}
    \!\!=
    \begin{braid}%
      \draw[dots](4.6,4.5)--(6.8,4.5);
      \draw[dots](9.8,4.5)--(11,4.5);
      \draw[dots](3.9,0)--(5.7,0);\draw[dots](8.5,0)--(11,0);
      \draw[brown](1.2,5)--(4.4,5)--(4.4,4)--(1.2,4)--cycle;
      \draw[brown](6.8,5)--(9.4,5)--(9.4,4)--(6.8,4)--cycle;
      \foreach \k in {1.5,2,...,3.5} {
         \draw(\k,4)--(\k,0);
      }
      \draw(7,4)--(6,0);
      \draw(7.5,4)--(6.5,0);
      \draw(8.5,4)--(7.5,0);
      \draw(9,4)--(8,0);
      \draw[red](4,4)node[above]{$i$}--(9,1.5);\draw[red,dashed](9,1.5)--(11,0.5);
      \draw[red](8,4)node[above]{$i$}--(7,0);
      \greendot(7.35,2.35);
    \end{braid}
    \!\!+
    \begin{braid}
      \draw[dots](4.6,4.5)--(6.8,4.5);\draw[dots](9.8,4.5)--(11,4.5);
      \draw[dots](3.9,0)--(5.7,0);\draw[dots](8.5,0)--(11,0);
      \draw[brown](1.2,5)--(4.4,5)--(4.4,4)--(1.2,4)--cycle;
      \draw[brown](6.8,5)--(9.4,5)--(9.4,4)--(6.8,4)--cycle;
      \foreach \k in {1.5,2,...,3.5} {
         \draw(\k,4)--(\k,0);
      }
      \draw(7,4)--(6,0);
      \draw(7.5,4)--(6.5,0);
      \draw(8.5,4)--(7.5,0);
      \draw(9,4)--(8,0);
      \draw[red](4,4)node[above]{$i$}--(7.5,2)--(7,0);
      \draw[red](8,4)node[above]{$i$}--(7.6,2.4)--(9,1.5);
      \draw[red,dashed](9,1.5)--(11,0.5);
    \end{braid}
    $$
With the first diagram on the right-hand side, we continue moving the dot to the left.
With the second diagram there are three possibilities. First, if the rightmost
labelled string of residue~$i$ does not correspond to a node
at the right-hand end of a row of the diagram $\bnu$, as drawn in the diagram
above, then the diagram is zero by the row relation~\autoref{E:InsideRow}. If the
node is at the end of row~$(l,r)$ and $\nu^{(l)}_r=\nu^{(l)}_{r+1}$,
where $1\le l\le\ell$, then the diagram is again zero by the Garnir
relation~\autoref{E:SmallGarnir}. Otherwise, the right-hand string of residue~$i$
moves $n+1$ to a removable node in~$\bnu$. Let $\btau$ be the multipartition of
$n$ obtained by removing this node. Then $\btau\iarrow\bnu$, $\btau\gdom\bsig$
and  the diagram is equal to $v_\tnutau\psi_{w_{(\btau,\bsig)}}$.

The lemma follows by repeating this argument and noting that once the dot gets
all the way to the top of the diagram when we get zero by~\autoref{E:InsideRow}.
\end{proof}

We now generalise the definitions used in \autoref{L:OneNodePush}. Suppose that
$\bsig_0,\dots,\bsig_{z+1}$ are multipartitions of~$n$ such that $\bsig_k\iarrow\bnu$, for
$0\le k\le z+1$, ordered so that
$\bsig_{z+1}\gdom\bsig_{z}\gdom\dots\gdom\bsig_0$. Define permutations
$$w_{(\bsig_k,\bsig_{k-1}, \dots,\bsig_0)}=w_{(\bsig_1,\bsig_0)} w_{(\bsig_2,\bsig_1)}
             \dots w_{(\bsig_{k},\bsig_{k-1})},$$
for $0\le k\le z+1$. Observe that
$w_{(\bsig_{j+1},\bsig_{j})} w_{(\bsig_{k+1},\bsig_{k})} =  w_{(\bsig_{k+1},\bsig_{k})} w_{(\bsig_{j+1},\bsig_{j})}$, whenever $0 \leq j<k \leq z$.
For $0\le k\le z+1$ define $\bnu$-tableaux by
$$\t_{(\bsig_k,\bsig_{k-1}, \dots,\bsig_0)} = \tnu_{\bsig_k}
                    w_{(\bsig_k,\bsig_{k-1}, \dots,\bsig_0)}.$$
Then $\t_{(\bsig_k,\bsig_{k-1}, \dots,\bsig_0)}$ is standard.
Note that
$\psi_{w_{(\bsig_k,\bsig_{k-1},\ldots,\bsig_{0})}}$ is independent of the choice
of reduced expression for $w_{(\bsig_k,\bsig_{k-1},\ldots,\bsig_{0})}$, so that
$$v_{\t_{(\bsig_k,\bsig_{k-1}, \dots,\bsig_0)}} = v_{\tnu_{\bsig_k}}
            \psi_{w_{(\bsig_k,\bsig_{k-1},\ldots,\bsig_{0})}}.$$

\begin{Example}
\label{E:1node}
Suppose that $e=3$ and $\kappa=(0,2,0,1)$. Take $n+1=20$, $\bnu=(4,3,1|2|4,3|2,1)$ and $i=0$.  The removable $i$-nodes are those shaded below
$$\Bigg(\ \ShadedTableau[(4,0)]{{0,1,2,0},{2,0,1},{1}}\ \Bigg|\
          \ShadedTableau[((2,0)]{{2,0}}\ \Bigg|\
          \ShadedTableau[(4,0)]{{0,1,2,0},{2,0,1}}\ \Bigg|\
          \ShadedTableau[(1,-1)]{{1,2},{0}}\
\Bigg),$$
and the partitions $\btau$ such that $\btau\iarrow \bnu$ are
\begin{align*}
\bsig_0 & = (3,3,1|2|4,3|2,1), &
\bsig_1 & = (4,3,1|1|4,3|2,1), \\
\bsig_2 & = (4,3,1|2|3,3|2,1), &
\bsig_3 & = (4,3,1|2|4,3|2).
\end{align*}
Then
$w_{(\bsig_3,\bsig_2,\bsig_1,\bsig_0)} = (9,8,7,6,5,4)(13,12,11,10)(19,18,17,16,15,14)$
and
\[\t_{(\bsig_3,\bsig_2,\bsig_1,\bsig_0)} = \Bigg(\ \ShadedTableau[]{{1,2,3,9},{4,5,6},{7}}\ \Bigg|\
          \ShadedTableau[]{{8,13}}\ \Bigg|\
          \ShadedTableau[]{{10,11,12,19},{14,15,16}}\ \Bigg|\
          \ShadedTableau[]{{17,18},{20}}\
\Bigg).\]
\end{Example}

By definition, $w_{(\bsig_k,\bsig_{k-1},\ldots,\bsig_{0})}\in\Sym_n$ so $y_{n+1}$ commutes with
$\psi_{w_{(\bsig_k,\bsig_{k-1},\ldots,\bsig_{0})}}$, for $0\le k\le z+1$. Hence,
\autoref{L:OneNodePush} immediately implies the following.

\begin{Corollary}
\label{C:MultiPush}
  Suppose that $\bsig\iarrow\bnu$ and that $k\ge1$. Then
  $$ v_{\tnusig}y_{n+1}^k
=  \sum_{\substack{\bsig_k,\dots,\bsig_1\iarrow\bnu\\\bsig_k\gdom\dots\gdom\bsig_1\gdom\bsig}}
  v_{\tnu_{\bsig_k}} \psi_{(\bsig_k,\dots,\bsig_1,\bsig)} =
  \sum_{\substack{\bsig_k,\dots,\bsig_1\iarrow\bnu\\\bsig_k\gdom\dots\gdom\bsig_1\gdom\bsig}}
  v_{\t_{(\bsig_k,\dots,\bsig_1,\bsig)}}.$$
\end{Corollary}

We can now complete the proof of our Main Theorem in the case when
$\gamma=1$.

\begin{proof}[Proof of \autoref{T:OneNode}]
  We remind the reader that the rows of $\bnu$ are indexed by pairs
  $r=(\text{comp},\text{row})$, which determine a row in a component of~$\bnu$.
  These row indices are totally ordered lexicographically by~\autoref{E:RowOrder}.

By definition, $\blam,\bmu \iarrow \bnu$ are multipartitions such that $\bmu \gdom
\blam$ where $\bnu\backslash\blam$ is in row~$r$ of $\bnu$ and $\bnu\backslash
\bmu$ lies in row~$s$ for some $r<s$.

Let $\alpha_0,\dots,\alpha_{z+1}$ be the removable $i$-nodes in~$\bnu$ in
rows~$r$ to~$s$ inclusive, ordered by their row indices. Define multipartitions
$\blam=\bsig_0,\bsig_1,\dots,\bsig_{z+1}=\bmu$ of~$n$ by setting
$\bsig_k=\bnu\backslash\lbrace\alpha_k\rbrace$, for $0\le k\le z+1$. Then
$\bmu=\bsig_{z+1}\gdom\bsig_z\gdom\dots\gdom\bsig_0=\blam$ and $\bsig_k\iarrow\bnu$,
for $0\le k\le z+1$.
Recalling the notation from~\autoref{S:Restriction}, define $S^\bnu_{\bsig_k}:=\hat{S}^\bnu_{\bsig_k}=\check{S}^\bnu_{\bsig_k}$ and $S^\bnu_{\tnu_{\bsig_k}}:=\hat{S}^\bnu_{\tnu_{\bsig_k}}=\check{S}^\bnu_{\tnu_{\bsig_k}}$, so that $S^\bnu_{\bsig_k}$ has basis $\set{v_\t |\row_\t(n+1)>\row_{\tnu_{\bsig_k}}(n+1)}$,
for $0\le k\le z+1$.
It then follows from~\autoref{E:Snulam} that
the Specht module $S^\bnu$ has an $\R$-module filtration
$$S^\bnu\supseteq S^\bnu_{\tnulam}\supset S^\bnu_\blam \supseteq
S^{\bnu}_{\tnumu}\supset S^{\bnu}_\bmu\supseteq0.$$

Recalling the framework of~\autoref{L:EllersMurray}, define an idempotent $e_n(i)\in\R[n+1]$ by
$$e_n(i) = \sum_{(i_1,\ldots,i_n) \in I^n}e(i_1,\ldots,i_n,i)$$
and let $L = e_n(i)y_{n+1}^{z+1}$. Then $L$ commutes with $\R$ and $\deg L = 2(z+1)$.
We claim that multiplication by~$L$ induces a non-zero homomorphism from
$S^\blam$ to $S^\bmu$.

Suppose that $x \in S^\bnu_\blam$.  Then $x$ is a $\Z$-linear combination of elements $v_\t$ indexed by standard tableaux $\t$ such that $\row_\t(n+1) > r$. If $\res_\t(n+1) \ne i$ then $v_\t e_n(i)=0$ by~\autoref{E:tableauxResidue}, so that $v_\t L=0$.  So
suppose that $\res_\t(n+1)=i$.  Then $v_\t e_n(i) = v_\t$.
If $\row_\t(n+1)>s$ then $\bmu \not \gedom \t_{\downarrow n}$ and $v_\t y^{z+1}_{n+1} \in S^\bnu_\bmu$ by
\autoref{L:YDown}.
Otherwise $\Shape(\t_{\downarrow n}) = \bsig_k$,
for some $1 \leq k \leq z+1$.  Now $v_\t = v_{\tnu_{\bsig_k}} \psi_w$, for $w \in \Sym_n$ and
\[v_\t y_{n+1}^{z+1} = v_{\tnu_{\bsig_k}} \psi_w y_{n+1}^{z+1}
= v_{\tnu_{\bsig_k}} y_{n+1}^{z+1} \psi_w \in S^\bnu_{\bsig_{z+1}}=S^\bnu_\bmu,\]
by~\autoref{C:MultiPush}. Therefore, if $x \in S^\bnu_\blam$ then
$xL \in  S^{\bnu}_\bmu$.

Now consider $v_\tnulam L$.  Using~\autoref{C:MultiPush} again,
\begin{align*}
  v_\tnulam L & = v_{\tnulam} y_{n+1}^{z+1}
  \equiv v_{\t_{(\bmu,\bsig_z,\dots,\bsig_1,\blam)}}
  \pmod{S^\bnu_{\bmu}}.
\end{align*}
By construction, $v_{\t_{(\bmu,\bsig_z,\dots,\bsig_1,\blam)}} \in S^{\bnu}_{\tnumu}$ because
$\row_{\t_{(\bmu,\bsig_z,\ldots,\bsig_1,\blam)}}(n+1) = s$ which also implies
that $v_{\t_{(\bsig_z,\dots,\bsig_1,\blam)}} \notin S^\bnu_\bmu$.
Therefore, in the notation of \autoref{L:EllersMurray}, the map~$\theta_L$
induces a non-zero $\R$-module homomorphism from $S^\blam\<\delta\>$ to~$S^\bmu$,
where $\delta=(\deg\tnulam-\deg\tlam)-(\deg\tnumu-\deg\tmu)+\deg L$.

To complete the proof of \autoref{T:OneNode} it remains to show that $\delta$, as
defined in the last paragraph, is equal to the number of addable or removable $i$-nodes in rows
$r+1$ to $s$.
Applying~\autoref{E:degrees}, $\deg\tnulam-\deg\tlam$ is equal to the number of
addable $i$-nodes minus the number of removable $i$-nodes in~$\bnu$ which are
strictly below row~$r$ while $\deg\tnumu-\deg\tmu$ is equal to the number of
addable $i$-nodes minus the number of removable $i$-nodes in~$\bnu$ which are
strictly below row~$s$.
Therefore,
$(\deg\tnulam-\deg\tlam)-(\deg\tnumu-\deg\tmu)$ is the number of addable
$i$-nodes minus the number of removable $i$-nodes in~$\bnu$ in rows~$r+1$
to~$s$. On the other hand, by definition, $\deg L=2(z+1)$ and~$\bnu$ has exactly
$z+1$ removable $i$-nodes in rows~$r+1$ to~$s$. This completes the proof.
\end{proof}

\subsection{Stubborn strings} \label{S:Stubborn} The proof of
\autoref{T:OneNode} is quite straightforward but extending this result to the
general case requires several new ideas. We start by introducing
\textit{stubborn strings} which will allow us to prove a significant
generalisation of \autoref{L:OneNodePush}. The results and techniques introduced
in this section are likely to be of independent interest.

Fix integers $n\ge0$ and $\gamma>0$ and a multipartition
$\bnu\in\Parts[n+\gamma]$. As we will see, when $e=2$ our Main Theorem reduces
to \autoref{T:OneNode}, so we will assume that $e\ne2$. This is advantageous
because it simplifies relations~\autoref{E:psi^2} and~\autoref{E:braid}.

The basis $\set{v_\u|\u\in\Std(\bnu)}$ of the Specht module $S^\bnu$ defined in
\autoref{S:SpechtModules} is determined by a choice of reduced expressions and
the arguments which follow are easier for certain choices of these reduced
expressions.  The following lemma, which we leave as an easy
exercise, will help us make these choices.

\begin{Lemma}
\label{L:ReducedExps}
  Suppose that $\u\in\Std(\bnu)$ and let
  $\bsig=\Shape(\u_{\downarrow n})$. Then there exists a unique tableau
  $\s\in\Std(\bnu)$ such that $\s_{\downarrow n}=\t^\bsig$ and
  $\s^{-1}(n+h)=\u^{-1}(n+h)$, for $1\le h\le\gamma$. Moreover,
  $\u=\s w$,  for some $w\in\Sym_n$ and $\ell(d(\u)) = \ell(d(\s))+\ell(w)$.
\end{Lemma}

If $\u\in\Std(\bnu)$ then write $\u=\s w$ as in \autoref{L:ReducedExps}. Then
$d(\u)=d(\s)w$ with $\ell(d(\u))=\ell(d(\s))+\ell(w)$. In the sequel we always assume that $v_\u=v_\s\psi_w$,
for some fixed (and for the moment arbitrary) reduced expressions for~$d(\s)$
and~$w$.

Let $\W$ be the set of words in
$\{\psi_1,\ldots,\psi_{n+\gamma-1},y_1,\ldots,y_{n+\gamma}\}$. We identify a
word $\varpi\in\W$ with a braid diagram $B_\varpi$ in~$S^\bnu$ by concatenating
the diagrams for~$v_{\tnu}$ and the letters in~$\varpi$ the natural way. By
definition, $B_\varpi$ is a braid diagram which represents an element
of~$S^\bnu$. We have not yet attempted to apply the relations to rewrite this
element in terms of the basis $\set{v_\t|\t\in\Std(\bnu)}$ of~$S^\bnu$.

If $\varpi \in \W$ define $\pi_\varpi$ to be the underlying permutation in
$\Sym_{n+\gamma}$ determined by ignoring any dots on the strings in~$B_\varpi$.
Define $\res(\varpi) = \res(\tnu \pi_\varpi)$.

Finally, given $\varpi\in\W$ let $\approx_\varpi$ be the equivalence relation on
$\{1,2,\dots,n+\gamma\}$ generated by $l\approx_\varpi m$ whenever the
$l$-string and $m$-string in~$B_\varpi$ intersect and have the same residue.

\begin{Proposition} \label{P:Rewrite}
  Suppose that $\varpi\in\W$. Then there exist integers $a_\s\in\Z$ such that
  \[v_{\tnu} \varpi = \sum_{\s\in\Std(\bnu)}a_\s v_\s\]
  where $a_\s\ne0$ only if $\res(\s)=\res(\varpi)$ and
  $(\pi_{\varpi}d(\s)^{-1})(m)\approx_\varpi m$, for $1 \le m \le n+\gamma$.
\end{Proposition}

\begin{proof} Since $\set{v_\s|s\in\Std(\bnu)}$ is a basis of $S^\bnu$ we can
  certainly write $v_{\tnu}\varpi=\sum_\s a_s v_\s$, for some $a_\s\in\Z$. If
  $\bi=\res(\varpi)$ then $v_{\tnu}\varpi=v_{\tnu}\varpi e(\bi)$. So $a_\s\ne0$ only
  if $\res(\s)=\bi$ by~\autoref{E:tableauxResidue}.

  To prove the second claim in the lemma observe that if $i\in I$ then all of
  the KLR relations in \autoref{D:KLR}, interpreted in terms of braid diagrams,
  either preserve or break the $\iicrossing$-crossings in a braid diagrams. In
  particular, the relations in~$\R$ never create new $\iicrossing$-crossings in
  the sense that two strings of the same residue intersect in the diagrams
  obtained after applying the relations only if they intersected before the
  relations were applied. Similarly,  relations~(a)--(c) in
  \autoref{D:SpechtModule} do not create new $\iicrossing$-crossings.

  Similarly, we claim that the homogeneous Garnir relations in
  \autoref{D:SpechtModule}~(d) do not generate new $\iicrossing$-crossings for
  any $i\in I$. When $e=0$ then there is nothing to prove so suppose that $e>0$.
  Fix a node $A\in\bnu$. Then the terms in the corresponding homogeneous Garnir
  relation are indexed by the tableaux $\set{\tA d|d\in\D_A}$. All of these
  tableaux are standard except for the Garnir tableau ~$\tA d_0$, where $d_0$ is
  the unique element of~$\D_A$ of maximal length. Inside the Garnir belt, the
  braid diagram for $d_0$ takes the form:
  $$
  \begin{braid}\tikzset{yscale=1.2}
    \foreach \k in {2,2.7,...,6.2} {
        \draw[blue!50](\k,4)--(\k+4.9,1);
    }
    \foreach \k in {6.9,7.6,...,10.4} {
        \draw[blue!50](\k,4)--(\k-4.9,1);
    }
    \draw(1.5,0.4)node[right]{$\underbrace{\hspace*{13mm}}_{ec_A\text{ nodes}}$};
    \draw(6.4,0.4)node[right]{$\underbrace{\hspace*{18mm}}_{ea_A\text{ nodes}}$};
  \end{braid}
  $$
  where $a_A$ and $c_A$ are the number of $e$-bricks in the first and second
  rows of the Garnir belt, respectively. It follows that if $d\in\D_a$ then two
  strings form an $\iicrossing$-crossing in the diagram for~$\tA d$ only if they
  form an $\iicrossing$-crossing in the diagram for~$\tA d_0$. Consequently, the
  homogeneous Garnir relations do not create new $\iicrossing$-crossings, as
  claimed.

  We have shown that the none of the relations in $S^\bnu$ generate new
  $\iicrossing$-crossings. Therefore, if $a_\s\ne0$ and the $l$-string and the $m$-string
  in~$d(\s)$ intersect and have the same residue then $l\approx_\varpi m$. Note
  that $d(\s)^{-1}(m)$ is the number at the top of the $m$-string in the diagram
  for~$\s$. It follows that $(\pi_{\varpi}d(\s)^{-1})(m)\approx_\varpi m$ as
  required.
\end{proof}

\begin{Definition}
  Suppose that $\varpi \in \W$ and that $1 \le m \le n+\gamma$. Then
  the $m$-string is \textbf{stubborn} in~$B_{\varpi}$ if $\pi^{-1}_{\varpi}(m) \leq \pi^{-1}_{\varpi}(l)$ for all $l$ such that $l \approx_\varpi m$.
\end{Definition}

An immediate consequence of the definition and \autoref{P:Rewrite} is the following.

\begin{Corollary} \label{C:Stubborn}
  Suppose that $\varpi\in\W$ and that the $m$-string is stubborn in $B_\varpi$.  Then
  \[v_{\tnu}\varpi = \sum_{\s\in\Std(\bnu)}a_\s v_\s\]
  where $a_\s\ne0$ only if $\res(\s)=\res(\tnu\pi_\varpi)$ and
  $d(\s)^{-1}(m) \ge \pi_\varpi^{-1}(m)$.
\end{Corollary}

Define
$\Std_n(\bnu)=\set{\t\in\Std(\bnu)|\t_{\downarrow n}=\t^\bsig
                 \text{ where }\bsig=\Shape(\t_{\downarrow n})}$.
These tableaux have already appeared in \autoref{L:ReducedExps}. We are
particularly interested in the following subset of $\Std_n(\bnu)$: set
$$
\IniNu = \SET{\t\in\Std_n(\bnu)}%
    {$\res_\t(n+g)\ne\res_{\t}(n+h)$, for $1\le h<g\le\gamma$, and\\
    $\row_\t(n+h)=\row_\t(n+f)$ whenever there exist $f\le h\le g$ with $\row_\t(n+g)=\row_\t(n+f)$}
$$

The next result should be compared with \autoref{L:OneNodePush}.

\begin{Lemma} \label{L:Pushing}
  Suppose that $\t\in\IniNu$. Fix an integer~$g$ with $1\le g\le\gamma$ and
  suppose $f$ is minimal such that $1\le f\le g$ and  $\row_\t(n+f)=\row_\t(n+g)$.
  Then
  $$v_\t y_{n+g}=\sum_{s\in\Std(\bnu)}a_\s v_\s$$
  for some integers $a_\s$ such that $a_\s\ne0$ only if $\res(\s)=\res(\t)$ and
  $$\row_\s(n+h)\ge\row_\t(n+h),\qquad\text{ for }1\le h\le\gamma,$$
  and this inequality is strict whenever $f\le h\le g$.
\end{Lemma}

\begin{proof}
  In terms of diagrams, to prove the lemma we need to compute the following:
  $$
    v_\t y_{n+g}
    =\begin{braid}
      \draw[brown](4.5,5)--(9.4,5)--(9.4,4)--(4.5,4)--cycle;
      \draw[brown](1.6,5)--(3.4,5)--(3.4,4)--(1.6,4)--cycle;
      \draw[brown](10.6,5)--(13.4,5)--(13.4,4)--(10.6,4)--cycle;
      \draw[brown](16.6,5)--(19.4,5)--(19.4,4)--(16.6,4)--cycle;
      \foreach \k in {2,2.5} {
        \draw[blue!50](\k,4)--(\k-1,0)--(\k-1,-1);
      }
      \foreach \k in {5,5.5}{
        \draw[blue!50](\k,4)-- +(-1,-4)-- (\k-1,-1);
      }
      \foreach \k in {11,11.5,...,12.5}{
        \draw[blue!50](\k,4)-- +(-4,-4)-- (\k-4,-1);
      }
      \foreach \k in {17,17.5,...,19}{
        \draw[blue!50](\k,4)-- (\k-6,0)-- (\k-6,-1);
      }
      \draw(7.4 ,4)node[above]{$i_{f}$}--(15,0)--(15,-1)node[below]{$n_{\!f}$};
      \draw(8.2,4)node[above]{$i_g$}--(15.8,0)--(15.8,-1)node[below]{$n_{\!g}$};
      \draw[densely dashdotted,darkgreen!60]%
                (1,4)--(13.5,0)--(13.5,-1)node[below]{$n_1$};
      \draw[densely dashdotted,darkgreen!60](9,4)--(16.6,0)--(16.6,-1);
      \draw[densely dashdotted,,darkgreen!60](3,4)--(15,-0.2)--(18,-1);
      \draw[densely dashdotted,darkgreen!60](12.7,4)--(18.5,0)--(18.5,-1);
      \draw[densely dashdotted,darkgreen!60]%
          (13.2,4)--(19,0)--(19,-1)node[below]{$n_\gamma$};
      \greendot(15.8,-0.8);
    \end{braid}
 $$
 where for convenience we write $n_h$ instead of $n+h$ and $i_h = \res_\t(n+h)$
 for $1 \leq h \leq \gamma$. Thus, the $n_h$-string is the string of
 residue~$i_h$ that ends at $n+h$ in the bottom row of the diagram. As in the
 diagram, $n_f=n+f$ and $n_g=n+g$ are in the same row of~$\t$ and if $m\le n_f$
 is in a row of higher index in~$\t$ than~$n_f$ then the $m$-string must cross the
 $n_h$-string whenever $f\le h\le g$. The dashed lines indicate that we do not know
 (or care) where the $n_h$-strings are in the diagram if $h<f$ or if $h>g$. As
 we will see we can may ignore these lines because $i_h\ne i_g$, whenever
 $h\ne g$ and $1\le h\le\gamma$.

 By assumption the $(n+h)$-strings, for $1\le h\le\gamma$, have distinct
 residues so the fact that some of the $(n+h)$-strings might cross the
 $n+f,\dots,n+g$ strings will not cause us any difficulties. Exactly as in the
 proof of \autoref{L:OneNodePush} we can slide the dot up the $(n+g)$-string
 until we reach a string of residue~$i_g=\res_\t(n+g)$.
 Applying~\autoref{E:iiCrossing}, and the relation $v_{\tnu}y_{r}=0$ for all $1
 \le r \le n+\gamma$, shows that~$v_\t y_{n+g}$ is equal to the sum of diagrams
 of the form
  $$
  \begin{braid}
      \draw[brown](2,5)--(7.5,5)--(7.5,4)--(2,4);
      \draw[dots](5.4,4.5)--(7.2,4.5);
      \draw[brown,dashed](1,5)--(2,5); \draw[brown,dashed](1,4)--(2,4);
      \draw[dots](2.35,4.5)--(4.6,4.5);
      \draw(2,4)node[above]{$i_{f}$}--(11,0)--(11,-1)node[below]{$n_{f}$};
      \draw(4,4)--(13,0)--(13,-1);
      \draw(7,4)--(16,0)--(16,-1);
      \draw[densely dashdotted,darkgreen!60](-3,4)--(10.5,-0.2)--(17,-1);
      \draw[densely dashdotted,darkgreen!60](13,4)--(17.5,0)--(17.5,-1);
      \draw[densely dashdotted,darkgreen!60](13.5,4)--(18,0)--(18,-1)%
                   node[below]{$n_\gamma$};
      \draw[densely dashdotted,darkgreen!60](-5,4)--(9,0)--(9,-1)node[below]{$n_1$};
      \draw[blue!50](10,4)--(4,0)--(4,-1);
      \draw[blue!50](12,4)--(6,0)--(6,-1);
      \draw[red](11,4)node[above]{$i_g$}--($ (11,4)!0.4!(5,0) $)
                    --(14,0)--(14,-1)node[below]{$n_{\!g}$};
      \draw[red](5,4)node[above]{$i_g$}--($ (5,4)!0.4!(14,0) $)--(5,0)--(5,-1);
   \end{braid}
 $$
 Let  $\varpi$ be the word corresponding to one of these diagrams. Noting that
 the $n+1,\ldots,n+\gamma$ strings all have different residues, we see that for
 $1 \le h \le \gamma$, the $(n+h)$-string is stubborn.  Applying
 \autoref{C:Stubborn}, $$v_{\tnu}\varpi=\sum_{\s\in\Std(\bnu)}b_\s v_\s,$$ for
 $b_\s\in\Z$ such that $b_\s\ne0$ only if $\res(\s)=\res(\t)$ and
 $d(\s)^{-1}(n+h)\ge \pi_\varpi^{-1}(n+h)\ge d(\t)^{-1}(n+h)$, for $1\le
 h\le\gamma$. In particular, if $b_\s\ne0$ then
 $d(\s)^{-1}(n+g)>d(\t)^{-1}(n+g)$.  But this implies that
 $\row_\s(n+g)>\row_\t(n+g)$ because there is no position in $\row_\t(n+g)$
 which is to the right of the position containing $n+g$ in $\t$ and which has
 the same residue.
 Similarly, if $f\le h< g$ then $\row_\t(n+h)=\row_\t(n+g)$.  Then
 $\row_\s(n+h)>\row_\t(n+h)$ because otherwise
 $\row_\s(n+h)=\row_\t(n+h)$ which implies that the entry in~$\s$ which is in
 the same position as~$n+g$ in~$\t$ is both larger than $n+h$ and has
 residue~$\res_\t(n+g)$, since $\res(\s)=\res(\t)$, which is impossible. Hence,
 $\row_\s(n+h)\ge\row_\t(n+h)$, for $1\le h\le\gamma$, and this inequality is
 strict if $f\le h\le g$. This completes the proof.
\end{proof}

We also need the following more specialised version of \autoref{L:Pushing}.

\begin{Lemma} \label{L:GeneralSigma}
  Suppose that $\t \in \IniNu$.  Choose $1 \leq g \leq \gamma$ and suppose $1 \leq f \leq g$ is minimal such that $\row_\t(n+f)=\row_\t(n+g)$.
 Set $n_h = n+h$ and $i_h=\res_\t(n+h)$, for $1 \leq h \leq \gamma$.
Then $v_{\t} y_{n+g}$ is the sum of all diagrams of the form
$$
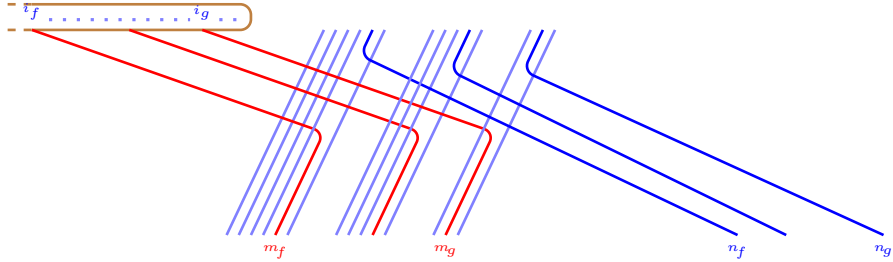
\begin{braid}\tikzset{xscale=0.8,yscale=1.7}
  \draw[brown](0,4.5)--(9.0,4.5)--(9.0,4)--(0,4);
  \draw[dots](7.7,4.2)--(8.7,4.2);
  \draw[brown,dashed](-1,4.5)--(0,4.5);
  \draw[brown,dashed](-1,4)--(0,4);
  \foreach \k in {12,12.5,13} {
    \draw[blue!50](\k,4)--+(-4,-4);
  }
  \draw[dots](0.7,4.2)--(6.5,4.2);

  \foreach \k in {0,4,7} {
    \draw[red](\k,4)--(\k+12,2)--(\k+10,0);
    \draw[blue!50](\k+9.5,0)--+(4,4);
    \draw[blue!50](\k+10.5,0)--+(4,4);
  }
  \draw[red](10,0)node[below]{$m_{\!f}$};
  \draw[red](17,0)node[below]{$m_{\!g}$};
  \draw(14,4)--(13.5,3.5)--(29,0)node[below]{$n_{\!f}$};
  \draw(18,4)--(17.2,3.2)--(31,0);
  \draw(21,4)--(20.2,3.2)--(35,0)node[below]{$n_{\!g}$};
  \draw(0,4)node[above]{$i_f$};
  \draw(7,4)node[above]{$i_g$};
  \draw[blue!50](13,0)--+(4,4);
  \draw[blue!50](12.5,0)--+(4,4);
\end{braid}
$$
where the sum is over all integers $1\leq m_f<\dots<m_g\le n$ such that
$\res_{\t}(m_h)=\res_{\t}(n+h)$ and $\row_\t(m_h)>\row_\t(n+h)$ for $f\le h\le
g$, and where all other strings appear in the same positions as they appear
in~$\t$.
\end{Lemma}

\textit{Warning!}  The notation~$m_h$ for the integers $m_f,\dots,m_g$ in
\autoref{L:GeneralSigma} is convenient because $m_h$ depends on~$n_h$, however,
the reader should not confuse $m_h$ and $m+h$ (since $m$ is not defined) even
though $n_h=n+h$, for $1\le h\le\gamma$.

\begin{proof}As in the proof of \autoref{L:Pushing}, in terms of diagrams
  $$
    v_{\t}y_{n+g}
    =\begin{braid}\tikzset{yscale=1.2}
      \draw[brown](1.5,5)--(6.4,5)--(6.4,4)--(1.5,4)--cycle;
      \draw[brown](7.6,5)--(9.4,5)--(9.4,4)--(7.6,4)--cycle;
      \draw[brown](10.6,5)--(13.4,5)--(13.4,4)--(10.6,4)--cycle;
      \draw[brown](16.6,5)--(19.4,5)--(19.4,4)--(16.6,4)--cycle;
      \foreach \k in {8,8.5}{
          \draw[blue!50](\k,4)-- (\k-4,0)--(\k-4,-1);
      }
      \foreach \k in {11,11.5,...,12.5}{
          \draw[blue!50](\k,4)-- (\k-4,0)--(\k-4,-1);
      }
      \foreach \k in {17,18.5,19}{
          \draw[blue!50](\k,4)-- (\k-5,0)--(\k-5,-1);
      }
      \draw(4,4)node[above]{$i_f$}--(15.1,0)--(15.1,-1)node[below]{$n_{\!f}$\ };
      \draw(5,4)node[above]{$i_g$}--(15.8,0)--(15.8,-1)node[below]{\ $n_{\!g}$};
      \draw[densely dashdotted,darkgreen!60](6,4)--(16.5,0)--(16.5,-1);
      \draw[densely dashdotted,darkgreen!60](9,4)--(17.5,0)--(17.5,-1);
      \draw[densely dashdotted,darkgreen!60](13,4)--(18,0)--(18,-1);
      \draw[densely dashdotted,darkgreen!60](0,4)--(13.5,-0.2)--(19,-1)%
            node[below]{$n_\gamma$};
      \greendot(15.8,-0.8);
    \end{braid}
 $$
 We saw in the proof of \autoref{L:Pushing} that the dashed strings do not move, so
 for clarity we will now ignore them. As in the proof of \autoref{L:Pushing},
 sliding the dot up the $(n+g)$-string shows that $v_\t y_{n+g}$ is the sum
 of all diagrams of the form
$$
  \begin{braid}\tikzset{xscale=0.9,yscale= 1.2}
   \draw[brown](2,5)--(7.5,5)--(7.5,4)--(2,4);
    \draw[dots](5.4,4.6)--(7.2,4.6);
   \draw[brown,dashed](1,5)--(2,5); \draw[brown,dashed](1,4)--(2,4);
    \draw[dots](2.35,4.6)--(4.6,4.6);
    \draw(2,4)node[above]{$i_f$}--(11,0)node[below]{$n_{\!f}$};
    \draw(4,4)--(13,0);
    \draw[red](14,0)node[below]{$n_{\!g}$};
    \draw(7,4)--(16,0);
    \draw[blue!50](10,4)--(4,0);
    \draw[blue!50](12,4)--(6,0);
    \YCrossing{i_g}{5}{11}{5}{14}
    \draw[red](5,0)node[below]{$m_{\!g}$};
 \end{braid}
$$
If $f=g$ then this proves the lemma, so we may assume that $f<g$.

Each of these diagrams above corresponds to a tableau $\s$ where $\s=\t
(m_g,n_g)$ and $\res_\t(m_g)=i_g=\res_\t(n_g)$. Note that $m_g\le n$ because the
residues of the $n_h$-strings are all distinct, for $f\le h\le g$.
Unfortunately, since $f<g$, the tableau~$\s$ is not standard because $n_{g-1}$
and $m_g$ both appear in the same row of~$\s$, with $n_{g-1}$ directly to the
left of~$m_g$, and $m_g<n_{g-1}$. To prove the lemma we show that this diagram
can be written as a linear combination of diagrams which correspond to tableaux
of the form $\t w$, where $w=(m_f,n_f)\dots (m_g,n_g)\in\Sym_{n+g}$, with
$n_f,\dots,n_g$ all appearing in later rows of~$\t w$ than $m_f,\dots,m_g$. As
we will see, the key observation is that each $n_h$-string crosses an
$m_h$ string of the same residue, for $f\le h\le g$.

Fix $h$ with $f\le h\le g$. We claim that the
diagram displayed above is equal to the sum of all diagrams of the form
$$
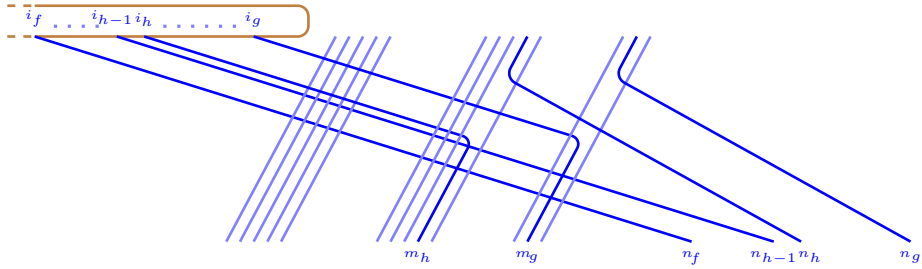
\begin{braid}\tikzset{xscale=0.9,yscale=1.7}
  \draw[brown](0,4.6)--(10.0,4.6)--(10.0,4)--(0,4);
  \draw[brown,dashed](-1,4.6)--(0,4.6);
  \draw[brown,dashed](-1,4)--(0,4);
  \draw(0,4)node[above]{$i_f$}--(24,0)node[below]{$n_{\!f}$};
  \draw(3,4)node[above]{\!\!\!$i_{h-1}$}--(27,0)node[below]{$n_{h-1}$};
  \foreach \k in {11,11.5,...,13} {
    \draw[blue!50](\k,4)--+(-4,-4);
  }
  \draw[dots](0.7,4.2)--(2.6,4.2);
  \draw[dots](4.7,4.2)--(7.6,4.2);
  \foreach \k in {4,8} {
    \draw(\k,4)--(\k+12,2)--(\k+10,0);
    \draw[blue!50](\k+9.5,0)--+(4,4);
    \draw[blue!50](\k+10.5,0)--+(4,4);
  }
  \draw(14,0)node[below]{$m_{h}$};
  \draw(18,0)node[below]{$m_{\!g}$};
  \draw(18,4)--(17.2,3.2)--(28,0)node[below]{\quad$n_{h}$};
  \draw(22,4)--(21.2,3.2)--(32,0)node[below]{$n_{g}$};
  \draw(4,4)node[above]{$i_{h}$};
  \draw(8,4)node[above]{$i_g$};
  \draw[blue!50](13,0)--+(4,4);
  \draw[blue!50](12.5,0)--+(4,4);
\end{braid}
$$
where all other strings are in exactly the same position as they are in the
diagram for~$\t$. That is, we have replaced the tableau~$\s=\t(m_g,n_g)$ in the
paragraph above with $\t(m_h,n_h)\dots (m_g,n_g)$, for some integers
$m_h<m_{h+1}<\dots<m_g$ such that $\res_\t(m_k)=\res_\t(n_k)$ and
$\row_\t(m_h)>\row_\t(n_h)$, for $h\le k\le g$. The sum is over all choices of
$m_h,\dots,m_g$ subject to these constraints. The lemma is exactly our claim
when $h=f$.

We prove our claim by downwards induction on~$h$. If $h=g$ then the claim is
automatically true. By induction we can assume that $f+1\le h\le g$ and that the
claim holds for $h$. For the inductive step we show that the last diagram, where
$m_h,\dots,m_g$ are fixed, is equal to the sum of all diagrams of the required
form with $m_{h-1}<m_h$. To do this we want to pull the $m_h$-string past the
string of residue~$i_{h-1}$, which is immediately to the left of this string.
Applying relation~\autoref{E:iiCrossing}, we can keep on doing this until we
reach an $i_{h-1}$-crossing by which time we will have a
diagram of the form
$$
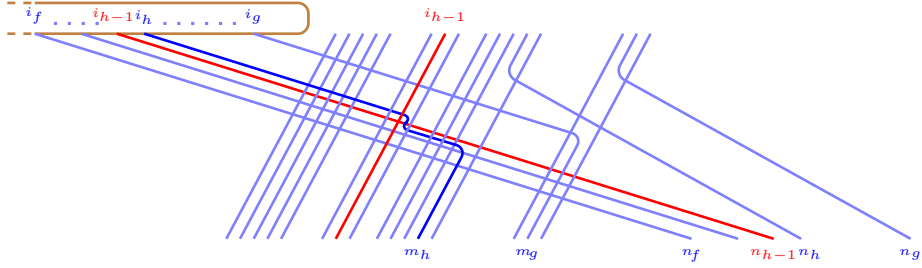
\begin{braid}\tikzset{xscale=0.9,yscale=1.7}
  \draw[brown](0,4.6)--(10,4.6)--(10,4)--(0,4);
  \draw[brown,dashed](-1,4.6)--(0,4.6);
  \draw[brown,dashed](-1,4)--(0,4);
  \draw[blue!50](0,4)node[above,blue]{$i_f$}--(24,0)node[below,blue]{$n_{\!f}$};
  \draw[blue!50](1.7,4)--(25.7,0);
  \draw[red](3,4)node[above]{\!\!$i_{h-1}$}--(27,0)node[below]{$n_{h-1}$};
  \foreach \k in {11,11.5,...,13} {
    \draw[blue!50](\k,4)--+(-4,-4);
  }
  \draw[dots](0.7,4.2)--(2.6,4.2);
  \draw[dots](4.7,4.2)--(7.6,4.2);
  \foreach \k in {4,8} {
    \draw[blue!50](\k+9.5,0)--+(4,4);
    \draw[blue!50](\k+10.5,0)--+(4,4);
  }
  \draw[blue!50](8,4)--(20,2)--(18,0);
  \coordinate (I) at ($ (4,4)!0.82!(16,2) $);
  \draw(4,4)node[above]{$i_{h}$}--(I)--($ (I)+(-0.55,-0.19) $)
                                   --($ (16,2)!0.115!(14,0) $)--(14,0);
  \draw(14,0)node[below]{$m_{h}$};
  \draw(18,0)node[below]{$m_{\!g}$};
  \draw[blue!50](14.5,4)--(10.5,0);
  \draw[red](15,4)node[above]{$i_{h-1}$}--(11,0);
  \draw[blue!50](15.5,4)--(11.5,0);
  \draw[blue!50](18,4)--(17.2,3.2)--(28,0)node[below,blue]{\quad$n_{h}$};
  \draw[blue!50](22,4)--(21.2,3.2)--(32,0)node[below,blue]{$n_{g}$};
  \draw[blue!50](8,4)node[above,blue]{$i_g$};
  \draw[blue!50](13,0)--+(4,4);
  \draw[blue!50](12.5,0)--+(4,4);
\end{braid}
$$
The two strings of residue $i_{h-1}$ and the string of residue $i_h$ are now in
exactly configuration where \autoref{E:iiCrossing} gives a non-trivial braid
relation. Therefore, when we pull the $m_h$-string past the $i_{h-1}$-crossing
we obtain two diagrams. The first diagram is obtained by applying the standard
three string braid relation, so that we pull the $m_h$-string past the
$i_{h-1}$-crossing. Therefore, in the first diagram we can keep pulling the
$m_h$-string to the left until we reach the next $i_{h-1}$-crossing, or until we
reach the top of the diagram in which case we get zero by~\autoref{E:InsideRow}.

The second diagram that arises when we apply the braid relation in
\autoref{E:iiCrossing} is obtained by cutting the $i_{h-1}$-crossing,
vertically, so it
has the form
$$
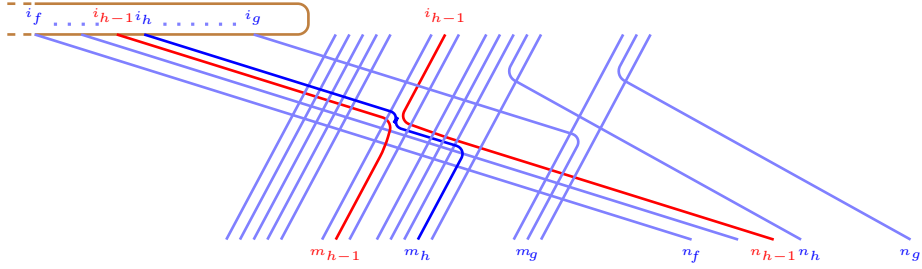
\begin{braid}\tikzset{xscale=0.9,yscale=1.7}
  \draw[brown](0,4.6)--(10,4.6)--(10,4)--(0,4);
  \draw[brown,dashed](-1,4.6)--(0,4.6);
  \draw[brown,dashed](-1,4)--(0,4);
  \draw[blue!50](0,4)node[above,blue]{$i_f$}--(24,0)node[below,blue]{$n_{\!f}$};
  \draw[blue!50](1.7,4)--(25.7,0);
  \draw[red](3,4)node[above]{\!\!$i_{h-1}$}--( $(3,4)!0.42!(27,0)$ )--( $(15,4)!0.55!(11,0)$ )--(11,0)node[below]{$m_{h-1}$};
  \draw[red](15,4)node[above]{$i_{h-1}$}--($(15,4)!0.42!(11,0)$)--($(3,4)!0.5!(27,0)$)
              --(27,0)node[below]{$n_{h-1}$};
  \foreach \k in {11,11.5,...,13} {
    \draw[blue!50](\k,4)--+(-4,-4);
  }
  \draw[dots](0.7,4.2)--(2.6,4.2);
  \draw[dots](4.7,4.2)--(7.6,4.2);
  \foreach \k in {4,8} {
    \draw[blue!50](\k+9.5,0)--+(4,4);
    \draw[blue!50](\k+10.5,0)--+(4,4);
  }
  \draw[blue!50](8,4)--(20,2)--(18,0);
  \coordinate (I) at ($ (4,4)!0.78!(16,2) $);
  \draw(14,0)node[below]{$m_{h}$};
  \draw(18,0)node[below]{$m_{\!g}$};
  \draw[blue!50](14.5,4)--(10.5,0);
  \draw(4,4)node[above]{$i_{h}$}--(I)--($ (I)+(-0.33,-0.22) $)
                                   --($ (16,2)!0.115!(14,0) $)--(14,0);
  \draw[blue!50](15.5,4)--(11.5,0);
  \draw[blue!50](18,4)--(17.2,3.2)--(28,0)node[below,blue]{\quad $n_{h}$};
  \draw[blue!50](22,4)--(21.2,3.2)--(32,0)node[below,blue]{$n_{g}$};
  \draw[blue!50](8,4)node[above,blue]{$i_g$};
  \draw[blue!50](13,0)--+(4,4);
  \draw[blue!50](12.5,0)--+(4,4);
\end{braid}
$$
Observe that this diagram is not of the required form because the
$n_{h-1}$-string crosses the $m_k$-string twice, for $h<k\le g$. If $h<k\le g$
then $i_{h-1}\ne i_k\pm1,i_k$ because the residues
$i_{h+1},\dots,i_g$ are all consecutive and, by assumption, $g\le\gamma<e$.
Therefore, we can move the $n_{h-1}$-string past the $m_k$-strings using the
relation $\psi_r^2e(\bi)=e(\bi)$, for $h<k\le g$. Hence, the last diagram is
equal to
$$
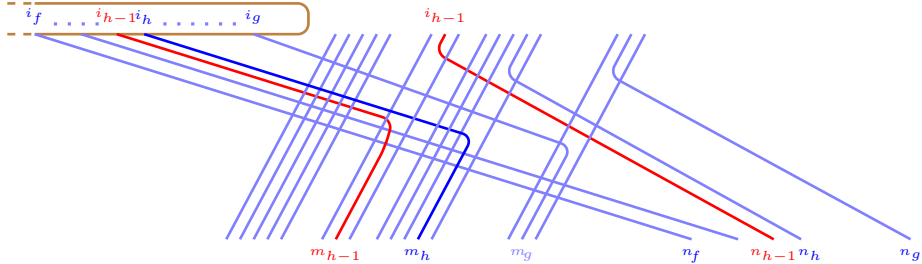
\begin{braid}\tikzset{xscale=0.9,yscale=1.7}
  \draw[brown](0,4.6)--(10,4.6)--(10,4)--(0,4);
  \draw[brown,dashed](-1,4.6)--(0,4.6);
  \draw[brown,dashed](-1,4)--(0,4);
  \draw[blue!50](0,4)node[above,blue]{$i_f$}--(24,0)node[below,blue]{$n_{\!f}$};
  \draw[blue!50](1.7,4)--(25.7,0);
  \draw[red](3,4)node[above]{$i_{h-1}$}--( $(3,4)!0.42!(27,0)$ )
                                   --( $(15,4)!0.55!(11,0)$ )--(11,0)node[below]{$m_{h-1}$};
  \draw[red](15,4)node[above]{$i_{h-1}$}--($(15,4)!0.09!(11,0)$)--(27,0)node[below]{$n_{h-1}$};
  \foreach \k in {11,11.5,...,13} {
    \draw[blue!50](\k,4)--+(-4,-4);
  }
  \draw[dots](0.7,4.2)--(2.6,4.2);
  \draw[dots](4.7,4.2)--(7.6,4.2);
  \foreach \k in {4,7.8} {
    \draw[blue!50](\k+9.5,0)--+(4,4);
    \draw[blue!50](\k+10.5,0)--+(4,4);
  }
  \draw(14,0)node[below]{$m_{h}$};
  \draw[blue!50](8,4)--(19.6,1.8)--(17.8,0)node[below]{$m_{\!g}$};
  \draw[blue!50](14.5,4)--(10.5,0);
  \draw(4,4)node[above]{$i_{h}$}--(16,2)--(14,0);
  \draw[blue!50](15.5,4)--(11.5,0);
  \draw[blue!50](18,4)--(17.2,3.2)--(28,0)node[below,blue]{\quad$n_{h}$};
  \draw[blue!50](21.8,4)--(21,3.2)--(32,0)node[below,blue]{$n_{g}$};
  \draw[blue!50](8,4)node[above,blue]{$i_g$};
  \draw[blue!50](13,0)--+(4,4);
  \draw[blue!50](12.5,0)--+(4,4);
\end{braid}
$$
This completes the proof of the inductive step and hence of our claim. Notice
that the endpoints of the $n+h,\dots,n+h$ strings and the strings which
start at $i_{h},\dots,i_g$ have not changed during the proof of the claim.
Taking $h=f$ we have proved the lemma.
\end{proof}

\subsection{One-row cyclotomic Carter-Payne homomorphisms} \label{S:OneRowCP}
We now prove the version of our Main Theorem which is closest in spirit to the
original Carter-Payne theorem. Even though we are again considering only a
special case of our Main Theorem the arguments in this section extend, almost
verbatim, to the general case. The proof of our Main Theorem requires
considerably more notation, however, so to make the argument clearer we consider
the  ``one-row case'' first and then highlight the changes needed for the general case.

We continue to assume that $\bnu$ is a multipartition of $n+\gamma$
where $n \geq 0$ and $\gamma>0$.

\begin{Definition}
  Suppose that $j\in I$. A \textbf{removable $j$-strip} of
  length~$\delta$ in $\bnu$ is a set of nodes
  $C=\{(l,r,c),(l,r,c+1),\dots,(l,r,c+\delta-1)\} \subseteq \bnu$
  such that $j=\res(l,r,c)$ and $\bnu\backslash C$ is a multipartition of~$n$.
  Set $\row(C)=(l,r)$.
\end{Definition}

That is, a removable $j$-strip in $\bnu$ is a set of nodes $C$ which are
contained in one row of~$\bnu$ such that $\bnu \backslash C$ is a multipartition
and the residue of the leftmost node in~$C$ has residue~$j$.

If $C$ is a removable $j$-strip in $\bnu$ let $\res(C)=\set{\res(\alpha)|\alpha\in C}$.
If $J\subseteq I$ then a removable $J$-node is a removable $i$-node for some
$i\in J$ and an addable $J$-node is an addable $i$-node  for some $i\in J$.

The main result of this section is the following:

\begin{Theorem}
\label{T:OneRow}
  Let $\blam$ and $\bmu$ be multipartitions of~$n$ and suppose that
  $\bmu\gdom\blam$ and that $\bnu=\blam\cup\bmu$ is a multipartition of
  $n+\gamma$ such that $\bnu\backslash\blam$ and $\bnu\backslash\bmu$ are both
  removable $j$-strips of length~$\gamma$, where $j\in I$ and
  $\gamma<|I|$.  Let $J=\res(\bnu\backslash\blam)=\res(\bnu\backslash\bmu)$.
  Then
  $$\Hom_{\R}(S^\blam\<a-b+2c\>,S^\bmu)\ne0,$$
  where
\begin{align*}
   a&=\#\set{\alpha \in \Add(\bnu)|\res(\alpha) \in J \text{ and } \row(\bnu\backslash\blam) < \row(\alpha) \leq \row(\bnu\backslash\bmu)},\\
   b&=\#\set{\alpha \in \Rem(\bnu)|\res(\alpha) \in J \text{ and } \row(\bnu\backslash\blam) <\row(\alpha) \leq \row(\bnu\backslash\bmu)},\\
   c&=\#\Set[70]{C}{$C$ is a removable $j$-strip of length at most $\gamma$ such
   that $\row(\bnu\backslash\blam)< \row(C) \leq \row(\bnu\backslash\bmu)$}.
  \end{align*}
  Moreover, $a-b+2c>0$.
\end{Theorem}

When $\gamma=1$ then it is not hard to see that $b=c$ so that in this case the
degree shift in \autoref{T:OneRow} is $a-b+2c=a+b$, which is strictly positive.
More generally, it is easy to see that if $\alpha$ is a removable $i$-node which
contributes to the number~$b$ in \autoref{T:OneRow} then either $\alpha$ is
contained in a removable $j$-strip of length at most $\gamma$ or the row below
$\row(\alpha)$ has an addable $i$-node for some $i\in J$. As $c>0$ it
follows that $a-b+2c>0$.

In view of the remarks in the last paragraph, \autoref{T:OneNode} is the special
case of \autoref{T:OneRow} corresponding to $\gamma=1$. In particular, if $e=2$
then \autoref{T:OneRow} is equivalent to \autoref{T:OneNode}. Consequently, as
we did in the last section, we can assume that $e\ne2$ when proving
\autoref{T:OneRow}.

For the rest of this section we assume that $\blam,\bmu\subset\bnu$ are
multipartitions satisfying the assumptions of \autoref{T:OneRow}. Once again,
the basic strategy for proving \autoref{T:OneRow} is to apply
\autoref{L:EllersMurray}. The first step is to set up the notation that we need
in order to define a non-zero map from~$S^\blam\<a-b+2c\>$ to $S^\bmu$.

Recall the definition of $\IniNu$ from just before \autoref{L:Pushing}.
Let $\tnulam$ be the unique standard $\bnu$-tableau in $\IniNu$ such that
$\Shape((\tnulam)_{\downarrow n}) = \blam$. Then $S^\blam\cong
S^\bnu_{\tnulam}/S^\bnu_\blam\<\deg\tlam-\deg\tnulam\>$ by \autoref{E:Snulam}.
Extending the definition of $\IniNu$, define
$$\IniNu[\blam]=\set{\t\in\IniNu|\res_\t(n+g)=\res_{\tnulam}(n+g),
               \text{ for }1\le g\le\gamma}.$$
Note that the residues of the nodes in $\bnu\backslash\blam$ are pairwise
distinct. If $\bsig\in\Parts$, $\bsig\subset\bnu$ and
$\res(\bnu\backslash\bsig)=\res(\bnu\backslash\blam)$ then
define $\tnu_\bsig=\tnu_\bsig(\blam)$ to be the unique standard $\bnu$-tableau
in~$\IniNu[\blam]$ such that $(\tnu_{\bsig})_{\downarrow n}=\t^\bsig$.

By assumption, $\bnu\backslash\blam$ and $\bnu\backslash\bmu$ are both removable
$j$-strips where $j=\res_{\tnulam}(n+1)$. Let
$C_0=\bnu\backslash\blam,C_1,\dots,C_{z+1}=\bnu\backslash\bmu$ be the removable
$j$-strips of length at most~$\gamma$ in~$\bnu$ ordered so that
$r_0<r_1<\dots<r_{z+1}$, where $r_k=\row(C_k)$ for $0\le k\le z+1$.  Define
standard tableaux $\t_0,\t_1,\dots,\t_{z+1}$ recursively by setting
$\t_0=\tnu_\blam$ and
$$\t_{k+1}=\t_k\prod_{h=1}^{g_k}(\eta^{(k)}_h,n+h),$$
where $g_k=\#C_k$ and $\eta^{(k)}_h$ is minimal such that
$\res_{\t_k}(\eta^{(k)}_h)=\res_{\t_k}(n+h)$ and
$$\t_k^{-1}(\eta^{(k)}_h)\in C_{k+1}\cup C_{k+2}\cup\dots\cup C_{z+1},$$
for $1\le h\le g_k$. Then, $1\le\eta^{(k)}_h\le n$ for $1\le h\le g_k$.

Before we give an example, we explain the meaning of these definitions.
As $j=\res_{\tnulam}(n+1)$, the removable $j$-strips
$C_0,C_1,\dots,C_{z+1}$ are in precisely those rows between $r_0$ and~$r_{z+1}$
where $n+1$ can appear in a standard tableau with the same residue sequence
as~$\tnulam$.  In such a tableau, if $n+1$ occurs in~$C_k$ then so does $n+h$
whenever $1\le h\le g_k=\#C_k$. Finally, the tableau~$\t_{k+1}$ is obtained from
$\t_k$ by swapping $n+h$ and $\eta^{(k)}_h$, for $1\le h\le g_k$, where
$\eta^{(k)}_h$ is the highest number in~$\t_k$ which is below $n+h$, has residue
$i_h=\res_{\tnulam}(n+h)$ and which sits in a removable $j$-strip. That is, the
tableaux $\t_k$ are constructed by successively moving $n+1,\dots,n+\gamma$ into
later rows as ``slowly'' as possible.

\begin{Example}
\label{Example:tk}
  Suppose that $e=3$ and $\kappa=(0,0,1)$. Let $\blam=(3,3,2,1|2|3,2)$ and $\bmu=(5,3,2,1|2|3)$  so that $n=16$,  $\gamma=2$ and $\bnu=(5,3,2,1|2|3,2)$. The
  residues in~$\bnu$ are
  $$\Bigg(\ \ShadedTableau[(4,0),(5,0),(1,-3)]{{0,1,2,0,1},{2,0,1},{1,2},{0}}\ \Bigg|\
          \ShadedTableau[(1,0),(2,0)]{{0,1}}\ \Bigg|\
          \ShadedTableau[(3,0),(1,-1),(2,-1)]{{1,2,0},{0,1}}\ \Bigg),$$
  where we have shaded all of the nodes which are contained in a removable
  $0$-strip of length at most 2.
  Therefore, $z=3$ and our definitions in this case give:
  \begin{center}
  \begin{tabular}{*6c}\toprule
    $k$ & $\t_k$ & $r_k$ &$g_k$ & $\eta^{(k)}_1$ & $\eta^{(k)}_2$ \\\midrule
    0&\ \Bigg(\ \ShadedTableau{{1,2,3,17,18},{4,5,6},{7,8},{9}}
    \ \Bigg|\ \ShadedTableau{{10,11}}
    \ \Bigg|\ \ShadedTableau{{12,13,14},{15,16}}\ \ \Bigg)
              &(1,1) & 2 & 9 & 11 \\[40pt]
    1&\Bigg(\ \ShadedTableau[(4,0),(5,0),(1,-3)]{{1,2,3,9,11},{4,5,6},{7,8},{17}}
    \ \Bigg|\ \ShadedTableau[(2,0)]{{10,18}}\ \Bigg|\
              \ShadedTableau{{12,13,14},{15,16}}\ \Bigg)
              & (1,4) &1 & 10 & - \\[40pt]
    2&\Bigg(\ \ShadedTableau[(1,-3)]{{1,2,3,9,11},{4,5,6},{7,8},{10}}%
    \ \Bigg|\ \ShadedTableau[(1,0)]{{17,18}}\ \Bigg|\
              \ShadedTableau{{12,13,14},{15,16}}\ \Bigg)
              &(2,1) & 2 & 14 & 16 \\[40pt]
    3&\Bigg(\ \ShadedTableau{{1,2,3,9,11},{4,5,6},{7,8},{10}}%
    \ \Bigg|\ \ShadedTableau[(1,0),(2,0)]{{14,16}}\ \Bigg|\
              \ShadedTableau[(3,0),(2,-1)]{{12,13,17},{15,18}}\ \Bigg)
              &(3,1)& 1 & 15 & - \\[40pt]
    4&\Bigg(\ \ShadedTableau{{1,2,3,9,11},{4,5,6},{7,8},{10}}%
    \ \Bigg|\ \ShadedTableau{{14,16}}\ \Bigg|\
              \ShadedTableau[(3,0),(1,-1)]{{12,13,15},{17,18}}\ \Bigg)
              &(3,2)& - & - & -
     \\\bottomrule
   \end{tabular}
   \end{center}
   where we have shaded the numbers which have moved at each step.
   Observe that $\bmu=\Shape((\t_4)_{\downarrow n})$.
\end{Example}

For $0\le k\le z+1$ set $\bsig_k=\Shape((\t_k)_{\downarrow n})$. Then
$\bmu=\bsig_{z+1}\gdom\bsig_z\gdom\dots\gdom\bsig_0=\blam$.
By definition, $\bsig_k\in\Parts$, $\bsig_k\subset\bnu$ and
$\res(\bnu\backslash\bsig_k)=\res(\bnu\backslash\blam)$. So,
$\tnu_{\bsig_k}\in\IniNu[\blam]$ and we can define permutations $w_k\in\Sym_n$ by
$\t_k=\tnu_{\bsig_k}w_k$, for $0\le k\le z+1$.

The next result is a straightforward application of the definitions. Note that
part~(b) follows by induction on~$k$ because if $1\le h\le \gamma$ then, by
definition, $n+h$ is contained in a removable $j$-strip. Therefore,
$r_k=\row_{\t_k}(n+1)\le\row_{\t_k}(n+h)$ because~$\t_k$ is standard, for
$0\le k\le z+1$.

\begin{Lemma}
\label{L:MildlyUsefulI}
Suppose that $0 \le k \le z+1$.  Then:
\begin{enumerate}
\item $\row_{\t_k}(n+h)=r_k$ for $1 \le h \le g_k$.
\item $\row_{\t_k}(n+h)>r_k$ for $g_k<h\le \gamma$.
\end{enumerate}
Moreover, $\t_{z+1}\gdom\t_z\gdom\dots\gdom\t_1\gdom\t_0=\tnulam$.
\end{Lemma}

We next choose a reduced expression for $w_k$ which will be used to fix the basis
element~$v_{\t_k}\in S^\bnu$, further specifying the choice made in \autoref{L:ReducedExps}.
If $0 \leq k \leq z$ and $1 \leq h \leq g_k$ then set
\begin{equation}
\label{E:lm}
l^{(k)}_h  = \tnu(\t_k^{-1}(n+h))\quad\text{and}\quad
m^{(k)}_h  = \tnu_{\sigma_k}(\t_k^{-1}(\eta^{(k)}_{h})).
\end{equation}
If $1\le l\le m\le n+\gamma$ let
$s(l,m)=(l,l+1)(l+1,l+2)\ldots (m-1,m)$ and set
$\psi(l,m)=\psi_l\dots\psi_{m-1}$. Note that $s(m,m)=1$ and $\psi(m,m)=1$.

\begin{Example}
\label{Example:tkII}
  Continuing \autoref{Example:tk}, the permutations $w_k$ are:
  $$\begin{array}{*6c}\\\toprule
      k & w_k & l^{(k)}_1 & l^{(k)}_2 & m^{(k)}_1 & m^{(k)}_2\\\midrule
       0 & 1  & 4 & 5 & 9 & 11\\
       1 & s(5,11)s(4,9) & 11& - & 11 & -\\
       2 & s(5,11)s(4,9) & 12& 13& 14&16\\
       3 & s(13,16)s(12,14)s(5,11)s(4,9) & 16 & - & 16&-\\
       4 & s(13,16)s(12,14)s(5,11)s(4,9) &-&-&-&-\\
       \bottomrule
   \end{array}
  $$
\end{Example}

It might help the reader to keep the last two examples in mind during the proof
of the following lemma.  The equation in \autoref{L:GoodReduced} is not as
complicated it looks: $g_k-h+1$ appears only because we want the right-hand
product in the formula to be written in terms of increasing values of~$h$.

\begin{Lemma} \label{L:GoodReduced}
Suppose that $0 \leq k \le z$.  Then
\[w_{k+1} = \prod_{h=1}^{g_k} s(l^{(k)}_{g_k-h+1}, m^{(k)}_{g_k-h+1}) \cdot w_k.\]
Furthermore, $\ell(w_{k+1}) =\ell(w_k)+\displaystyle\sum_{h'=1}^{g_k} (m^{(k)}_{h'}-l^{(k)}_{h'})$.
\end{Lemma}

\begin{proof}
  Directly from the definitions, if $0\le k\le z$ then
\begin{align*}
\tnu_{\sigma_{k+1}} w_{k+1} &= \t_{k+1}
                   = \t_k \prod_{h=1}^{g_k}(\eta_h^{(k)},n+h)
                   =\tnu_{\sigma_k} w_k \prod_{h=1}^{g_k}(\eta^{(k)}_h,n+h).
\end{align*}
On the other hand, the definition of $l^{(k)}_j$ and $m^{(k)}_h$ gives that
\[\tnu_{\sigma_k} w_k \prod_{h=1}^{g_k} (\eta^{(k)}_h,n+h)
= \tnu_{\sigma_{k+1}}\prod_{h=1}^{g_k} s(l^{(k)}_{g_k-h+1},m^{(k)}_{g_k-h+1} ) \cdot w_k. \]
Combining these two equations shows that $w_{k+1} = \prod_{h=1}^{g_k}
s(l^{(k)}_{g_k-h+1}, m^{(k)}_{g_k-h+1}) \cdot w_k$ as claimed.
The formula for $\ell(w_{k+1})$ follows by noting that $\ell(s(l,m))=m-l$ if $l \leq m$, that
$\ell(w(l,l+1))=\ell(w)+1$ if and only if $w(l)<w(l+1)$ and that
\[l^{(k)}_{g_1} > \ldots > l^{(k)}_{1} > l^{(k-1)}_{g_{k-1}}>\ldots >l^{(k-1)}_{1}
> \ldots > l^{(1)}_{g_1}> \ldots > l^{(1)}_{1}.\]
We leave the details to the reader.
\end{proof}

Following~\autoref{L:GoodReduced}, we fix reduced expressions for $w_0,\dots,w_{z+1}$ by setting
$\psi_{w_0}=1$ and
\[\psi_{w_{k+1}} =  \prod_{h=1}^{g_k} \psi_{s(l^{(k)}_{g_k-h+1}, m^{(k)}_{g_k-h+1})} \psi_{w_k},\]
for $0 \leq k \le z$. Using this choice of reduced expressions, define
\begin{equation}
\label{E:vtk}
v_{\t_k} = v_{\tnu_{\sigma_k}} \psi_{w_k},
\end{equation}
for $0\le k\le z+1$. Notice that the definition of $v_{\t_k}$  is compatible with
the choice of reduced expression that we made after \autoref{L:ReducedExps}.

Having set up our notation, we are now ready to start the proof of
\autoref{T:OneRow}.
Define $\Psin=\<\psi_1,\ldots,\psi_{n-1}\> \subset \R$ to be the (non-unital) subalgebra of
$\R$ generated by $\psi_1,\dots,\psi_{n-1}$.

\begin{Lemma}
\label{L:LittleLemmaI}
Suppose that $\t\in\IniNu[\blam]$ and that $\row_\t(n+1)=r_k$ for some $0 \leq k \leq z$.
 Let $k \le l\le z$. Then
  \[v_{\t}y_{n+g_k} y_{n+g_{k+1}} \ldots y_{n+g_l}
          = \sum_{\s \in \IniNu[\blam]} a_\s v_\s \psi_\s,\]
where $a_\s \in \Z$ and $\psi_\s \in \Psin$.  Moreover, if $a_\s \neq 0$ then
$\row_{\s}(n+1) \ge r_{l+1}$
and $\row_{\s}(n+g) \ge \row_{\t}(n+g)$,
for $1 \le g \le \gamma$.
\end{Lemma}

\begin{proof}By \autoref{L:ReducedExps}, if $\u\in\Std(\bnu)$ and $\res(\u)=\res(\tnulam)$ then
  $v_\u=v_\s\psi_w$ for a unique $\s\in\IniNu[\blam]$ and $w\in\Sym_n$, so that
  $v_\u y_{n+g}=v_\s y_{n+g}\psi_w$ where $\psi_w\in\Psin$.
  Observe that $y_{n+g}$ commutes with $\Psin$, for $1\le g\le\gamma$.
  Therefore, the lemma follows by applying \autoref{L:Pushing} $(l-k+1)$ times.
\end{proof}

Recall from \autoref{S:Restriction} that
$\hat{S}^\bnu_\bmu=\<v_\t\mid\t\in\Std(\bnu)\text{ and }\bmu \not\gedom \Shape(\t_{\downarrow n})\>$
is an $\R$-submodule of $S^\bnu$. For the next result we note that if
$\t\in\IniNu[\blam]$ and $r_0\le\row_\t(n+1)\le r_{z+1}$ then it is not
necessarily true that $\t=\tnu_{\bsig_k}$, for some $k$ with $0\le k\le z+1$.

\begin{Lemma} \label{L:PushI}
Suppose $\t \in \Std(\bnu)$, with $\res(\t)=\res(\tnulam)$, and that there exist
$0 \leq l \leq z+1$  and $1 \leq g \leq \gamma$ such that $\row_\t(n+1) \geq r_l$ and
$\row_\t(n+g) > \row_{\t_l}(n+g)$.  Then
$$v_\t \prod_{j=l}^{z} y_{n+g_j} \in \hat{S}^\bnu_\bmu,$$
where the product is assumed to be equal to 1 if $l=z+1$.
\end{Lemma}

\begin{proof}
First note that if $\res(\t)=\res(\tnulam)$ then by~\autoref{L:ReducedExps},
there exist $\s \in \Std^e_\blam(\bnu)$ and $w \in \Sym_n$ such that $v_\t = v_\s \psi_w$
where $\psi_w$ commutes with all the terms in the product.
Since $\hat{S}^\bnu_\bmu$ is an~$\R$-module and $\psi_w \in \R$, in order to prove
the lemma it is sufficient to consider only the case that $\t \in \Std^e_\blam(\bnu)$.

Suppose that there exist $0 \leq l \leq z+1$  and $1 \leq g \leq \gamma$ such that $\row_\t(n+1) \geq r_l$ and
$\row_\t(n+g) > \row_{\t_l}(n+g)$.
If $\row_\t(n+g)>r_{z+1}$
then by~\autoref{L:Pushing}, $v_\t \prod_{j=l}^{z} y_{n+g_j}$
is a linear combination of basis elements $v_\a$ where $\row_\a(n+g) >r_{z+1}$.
Then $\bmu \not \gedom \Shape(\a_{\downarrow n})$ and so
$v_\a \in \hat{S}^\bnu_{\bmu}$.
So we may assume that $0 \leq l \leq z$ and that $\row_\t(n+g)\leq r_{z+1}$, that is,
$\row_\t(n+g)= r_k$ for some $0 \leq k \leq z+1$.

Suppose first that $g=1$, so that $\row_\t(n+1)
=r_k$ for some $l<k\leq z$.
 By \autoref{L:LittleLemmaI} we can write $v_\t  y_{n+g_k} \ldots y_{n+g_{z}}$ as
a linear combination of terms $v_\u \psi_\u$, where $\u\in\IniNu[\blam]$,
$\psi_{\u}\in\Psin$ and $\row_\u(n+1)\ge r_{z+1}$. Then
by~\autoref{L:Pushing} and~\autoref{L:Filtration},
$v_\u \psi_\u y_{n+g_l} \ldots y_{n+g_{k-1}}$ is a linear combination of terms
$v_\a$ where $\a \in \Std(\bnu)$ and $\row_\a(n+1)>r_{z+1}$,
so that $v_\a \in \hat{S}^\bnu_\bmu$.

Now suppose that $2 \leq g \leq \gamma$.
By the comments above, we may assume that $\row_\t(n+1)=r_l$.
As $\t\in\IniNu[\blam]$ and $\row_\t(n+1)=r_l$ we have
$\row_{\t_l}(n+g)>r_l$ because $\t\in\IniNu[\blam]$ and $n+g$ is the only number
larger than $n+1$ which has residue $\res_{\t_l}(n+g)$.  So $\row_{\t_l}(n+g)=r_k$ for some $l<k\leq z+1$.
By
\autoref{L:LittleLemmaI}, $v_\t y_{n+g_l}\ldots y_{n+g_{n+k-1}}$ is equal to a
linear combination of terms $v_\u \psi_\u$, where $\u\in\IniNu[\blam]$,
$\psi_\u\in\Psin$, $\row_\u(n+1)\ge r_k$ and $\row_\u(n+g) \geq \row_\t(n+g)>\row_{\t_l}(n+g) = r_k$.
However, since $\row_\u(n+g)>r_k$ and $\u\in\IniNu[\blam]$ this forces
$\row_\u(n+1)>r_k$ because, again, $n+g$ is the only number larger than $n+1$ with
residue $\res_\t(n+g)$.
Therefore, $v_\u y_{n+g_k}\dots y_{n+g_z}\in
\hat{S}^\bnu_\bmu$ by the argument above.  It follows that
$v_\t y_{n+g_l}\dots y_{n+g_z}\in \hat{S}^\bnu_\bmu$, as required.
\end{proof}

The next result, which generalises \autoref{L:OneNodePush}, will show that the
map we construct is non-zero.

\begin{Proposition}
\label{P:yngI}
  Suppose that $0\le k\le z$. Then there exist integers $a_\s$ such that
  $$v_{\t_k}y_{n+g_k}=v_{\t_{k+1}}+\sum_{\s\in\Std(\bnu)}a_\s v_\s,$$
  where $a_\s\ne0$ only if $\res(\s)=\res(\tnulam)$, $\row_\s(n+1)\ge r_{k+1}$ and  either
  there exists $1 \leq g \leq g_k$ with $\row_\s(n+g)>\row_{\t_{k+1}}(n+g)$ or
  $\t_{z+1}\not\gedom\s$ $($or both$)$.
\end{Proposition}

\begin{proof}
It follows from~\autoref{L:YDown} that
$v_{\t_k} y_{n+g_k} = \sum_{\u \in \Std(\bnu)} b_\u v_\u$
where $b_\u \neq 0$ only if~$\res(\u)=\res(\t_k)=\res(\tnulam)$.
By construction, $v_{\t_k}y_{n+g_k}=v_{\tnu_{\bsig_k}} \psi_{w_k} y_{n+g_k} = v_{\tnu_{\bsig_k}}y_{n+g_k}\psi_{w_k}$,
which, by~\autoref{L:Pushing} and~\autoref{L:Filtration}, is a linear sum of
basis elements $v_\u$ where $\u \in \Std(\nu)$ is such that
$\row_\u(n+1)\geq r_{k+1}$.

It remains to show that $b_{\t_{k+1}} =1$ and that if $\u\ne t_{k+1}$ and $b_\u
\neq 0$ then either there exists $1 \leq g \leq g_k$ such that
$\row_\u(n+g)>\row_{\t_{k+1}}(n+g)$ or $\t_{z+1} \not \gedom \u$.  Note that the
last two cases are not mutually exclusive, however, $\t_{k+1}$ does not satisfy
either of these conditions.

By \autoref{L:GeneralSigma}, $v_{\tnu_{\bsig_k}}y_{n+g_k}$ is equal to
the  sum of all diagrams of the form
$$
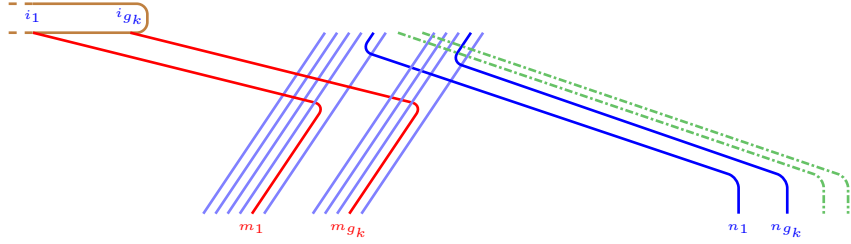
\begin{braid}\tikzset{xscale=0.8,yscale=1.2}
  \draw[brown](0,4.8)--(4.7,4.8)--(4.7,4)--(0,4);
  \draw[brown,dashed](-1,4.8)--(0,4.8);
  \draw[brown,dashed](-1,4)--(0,4);
  \foreach \k in {12,12.5,13} {
    \draw[blue!50](\k,4)--+(-5,-5);
  }
  \foreach \k in {0,4} {
    \draw[red](\k,4)--(\k+12,2)--(\k+9,-1);
    \draw[blue!50](\k+8.5,-1)--+(5,5);
    \draw[blue!50](\k+9.5,-1)--+(5,5);
  }
  \draw[red](9,-1)node[below]{$m_1$};
  \draw[red](13,-1)node[below]{$m_{g_k}$};
  \draw(14,4)--(13.5,3.5)--(29,0)--(29,-1)node[below]{$n_1$};
  \draw(18,4)--(17.2,3.2)--(31,0)--(31,-1)node[below]{$n_{g_k}$};
  \draw(0,4)node[above]{$i_1$};
  \draw(4,3.95)node[above]{$i_{g_k}$};
  \draw[blue!50](12,-1)--+(5,5);
  \draw[blue!50](11.5,-1)--+(5,5);
  \draw[densely dashdotted,darkgreen!60](15,4)--(32.5,0)--(32.5,-1);
  \draw[densely dashdotted,darkgreen!60](16,4)--(33.5,0)--(33.5,-1);
\end{braid}
$$
where the sum is over all integers $1\leq m_1<\dots<m_{g_k}\le n$ such that
$\res_{\tnu_{\bsig_k}}(m_h)=\res_{\tnu_{\bsig_k}}(n+h)$ and
$\row_{\tnu_{\bsig_k}}(m_h)>\row_{\tnu_{\bsig_k}}(n+h)$, for $1\le h\le g_k$,
and where all other strings appear in the same positions that they do
in~$\tnu_{\bsig_k}$.

Fix $m_1<\dots<m_{g_k}$, as in the last paragraph, and let $\m
=(m_1,\ldots,m_{g_k})$.  Let $\varpi_\m\in\W$ be the word in~$\W$
determined by the braid diagram corresponding to~$\m$, as drawn above,
and let $\pi_{\varpi_\m}$ be the corresponding permutation.
 Then
we need to compute
$v_{\tnu} \varpi_\m \psi_{w_k}$, for all $\m$.

First, suppose that $m_h=m^{(k)}_h$, for $1\le h\le g_k$. Then
 by~\autoref{L:GoodReduced} and the discussion following it,
\[v_{\tnu} \varpi_\m \psi_{w_k} = v_{\tnu_{\bsig_{k+1}}} \prod_{h=1}^{g_k} \psi_{s(l^{(k)}_{g_k-h+1},m^{(k)}_{g_k-h+1})} \psi_{w_k} = v_{\t_{k+1}}.\]
This implies that $v_{\t_{k+1}}$ appears in $v_{\t_k}y_{n+g_k}$ with
coefficient~$1$ because, as we will see,~$v_{\t_{k+1}}$ only appears in $v_{\tnu} \varpi_\m
\psi_{w_k}$ when $\m =(m_1^{(k)},\ldots,m_{g_k}^{(k)})$.

Now suppose that there exists $1 \leq g \leq g_k$ with $m_g > m^{(k)}_g$.  Note
that the $(n+g)$-string is stubborn in $\varpi_\m$.
By~\autoref{C:Stubborn}, $v_{\tnu} \varpi_\m$ is a linear combination of
basis elements $v_\a$ indexed by standard $\bnu$-tableaux $\a$ such that
$d(\a)^{-1}(n+g) \geq \pi_{\varpi_\m}^{-1}(n+g)$ and
$\res_\a(n+g)=\res_{\tnulam}(n+g)$, implying that
$\row_\a(n+g)>\row_{\t_{k+1}}(n+g)$.  By~\autoref{L:Filtration}, $v_\a
\psi_{w_k}$ is a linear combination of terms $v_\u$ where $\u \in \Std(\bnu)$ is
such that $\row_\u(n+g)>\row_{\t_{k+1}}(n+g)$.

Finally suppose that $m_h \leq h^{(k)}$, for $1 \leq h \leq g_k$, and that there
exists $1 \leq g \leq g_k$ with $m_g<m_g^{(k)}$. Without loss of
generality, we may assume $g$ is minimal with this property.  In particular,
note that $m_{g_k} \leq m_{g_k}^{(k)}$ so that none of the $n+1,\ldots,n+\gamma$
strings cross.
As in the proof of \autoref{L:GoodReduced},
\begin{equation}
\label{Eq:Swap}
\tnu \pi_{\varpi_\m} w_k = \t_k \prod_{h=1}^{g_k} (\eta_h,n+h),
\end{equation}
 where
$\eta_{h}=\t_k((\t_{\bsig_k}^\bnu)^{-1}(m_h))$, and
$\ell(\pi_{\varpi_\m})+\ell(w_k) = \ell(\pi_{\varpi_\m} w_k)$.  Let
$\t_\m = \tnu \pi_{\varpi_\m} w_k$.  Suppose first that $\t_\m$ is not row standard.
By~\autoref{Eq:Swap}, there exists $1 \leq h \leq g_k$ such that~$n+h$ is immediately to the left of an entry~$s$ in $\t_m$, where $1 \leq s \leq n$.  Then  $v_{\tnu}
\varpi_\m$ is zero by~\autoref{D:SpechtModule}~(c).  Now suppose that $\t_\m$ is row standard but not column standard.
Note that $\row_{\t_m}(n+f) \leq \row_{\t_\m}(n+h)$ for $1 \leq f\leq h \leq \gamma$.  By~\autoref{Eq:Swap}, there exists $h$ with $1 \leq h \leq \gamma$ such that $n+h$ lies above an entry $s$ in $\t_\m$, for some $1 \leq s \leq n$.
But then $v_{\tnu} \varpi_\m = 0$ by \autoref{E:SmallGarnir}, since either $e=0$
or $\gamma<e$.

So we may assume $\t_\m$ is  standard.
Let $x =\t_\m(\t_k^{-1}(n+g))$, so that $x$ occupies the same position in $\t_\m$ that $n+g$ occupies in $\t_k$. Note that $d(\t_\m)(x')=d(\t_{z+1})(x')$ for $1 \leq x' <x$.  Thus $\Shape((\t_{z+1})_{\downarrow x}) \not \gedom \Shape((\t_\m)_{\downarrow x})$
and so $\t_{z+1} \not \gedom \t_\m$.
Applying~\autoref{L:DomTableaux}, $v_{\tnu} \varpi_\m \psi_{w_k}$ is a linear combination of basis elements $v_\u$ where $\u\in \Std(\bnu)$ and $\u \gedom \t_\m$, so that $\t_{z+1} \not\gedom \u$.
\end{proof}

We are now ready to start proving \autoref{T:OneRow}. Recalling the framework of
\autoref{L:EllersMurray}, define
$$\LMult=y_{n+g_0}y_{n+g_1}\dots y_{n+g_z}e^\blam_\bmu,\qquad
\text{where}\quad e^\blam_\bmu= \sum_{\bj\in I^n}e(j_1,\dots,j_n,i_1,\dots,i_\gamma)$$
and $i_h=\res_{\tnulam}(n+h)$, for $1\le h\le\gamma$. Then $\LMult$ commutes
with $\R$, so the map
\begin{equation}
\label{E:thetaOneRow}
    \thea\map{S^\bnu}S^\bnu; \qquad x\mapsto x\LMult,\qquad\text{for }x\in S^\bnu,
\end{equation}
is an $\R$-module endomorphism of $S^\bnu$.

\begin{Example}
  Let $\blam$ and $\bmu$ be as in \autoref{Example:tk} and \autoref{Example:tkII}.
  Then
  $$\thea(x) =xy_{17}^2y_{18}^2 \sum_{\bj\in I^n} e(j_1,\dots,j_n,0,1),$$
  for all $x\in S^\bnu$.

In the notation of \autoref{T:OneRow}, $a=3$, $b=5$ and $c=4$, so that
$a-b+2c=6$ and, consequently, there exists a non-zero homomorphism
$S^\blam\<6\>\to S^\bmu$ by \autoref{T:OneRow}. We will show that $\thea$
induces such a map.  \end{Example}

\begin{proof}[Proof of \autoref{T:OneRow}]
  We follow the same basic strategy that we used to prove \autoref{T:OneNode}, that
 is, we engineer the circumstances that we need to apply \autoref{L:EllersMurray}.
  By \autoref{E:Snulam}, $S^\bnu$ has an $\R$-module filtration
  \begin{equation*}
    S^\bnu\supseteq \check{S}^\bnu_{\tnulam}\supset \check{S}^\bnu_\blam
    \supseteq \hat{S}^{\bnu}_{\tnumu}\supset \hat{S}^{\bnu}_\bmu\supset0
  \end{equation*}
  where $\check{S}^\bnu_{\tnulam}/\check{S}^\bnu_\blam
             \cong S^\blam\<\deg \tnulam -\deg\tlam\>$ and
  $\hat{S}^\bnu_{\tnumu} / \hat{S}^\bnu_\bmu \cong S^\bmu \< \deg\tnumu-\deg\tmu \>$.
  By \autoref{L:EllersMurray} in order to show that $\thea$ induces a map
  $S^\blam\<\delta\>\rightarrow S^\bmu$, for some $\delta \in\Z$, it is enough
  to prove the following three statements:
  $$\text{(A)}\quad \thea(\check{S}^\bnu_\blam)\subseteq \hat{S}^\bnu_\bmu,
  \qquad\text{(B)}\quad\thea(\check{S}^\bnu_{\tnulam})\subseteq \hat{S}^\bnu_{\tnumu}\qquad\text{and}
  \qquad\text{(C)}\quad \thea(v_{\tnulam})\notin \hat{S}^\bnu_\bmu.
  $$

  First consider~(A). Since $\thea$ is an $\R$-homomorphism it is enough to show
  that $\thea(v_\t)\in \hat{S}^\bnu_{\tnumu}$ whenever $\t\in\Std_n(\bnu)$ and
  $\Shape(\t_{\downarrow n}) \gdom \blam$. If $\t\notin\IniNu[\blam]$ then
  $\thea(v_\t)=0$ because $v_\t e^\blam_\bmu=0$ by~\autoref{E:tableauxResidue}.
  Therefore, we may assume that $\t\in\IniNu[\blam]$.  Since
  $\Shape(\t_{\downarrow n}) \gdom \blam$, we have $\row_\t(n+h) \geq r_0$ for
  $1 \leq h \leq \gamma$ with at least one of these equalities being strict, so
  since $\res(\t)=\res(\tnulam)$ we have $\row_\t(n+1)>r_0$.
  Setting $g=1$ and $l=0$ in~\autoref{L:PushI}, $\thea(v_\t)\in \hat{S}^\bnu_\bmu$
  so property~(A) holds.


  By definition, $\thea(v_{\tnulam})=v_{\tnulam}y_{n+g_0}\dots y_{n+g_z}$. To
  prove (B) and (C) we would like to argue by induction on~$k$ by considering
  the elements $v_{\tnulam}y_{n+g_0}\dots y_{n+g_k}$, where $0\le k\le z$. The
  difficulty is that when we write
  $v_{\tnulam}y_{n+g_0}\dots y_{n+g_k}=\sum_\t a_\t v_\t$ then
  it can happen that $a_\t\ne0$ when $\row_\t(n+h)<\row_{\t_z}(n+1)$ for some $h$
  with $1\le h\le\gamma$. To cater for this we have to argue slightly
  circuitously. For $0\le k\le z+1$ define
  $$\Stdlm=\set{\u\in\Std(\bnu)|\res(\u)=\res(\tnulam),
               \row_\u(n+1)=r_k\text{ and }\t_{z+1}\not\gedom\u}.$$
  To complete the proof we claim that if $0\le k\le z+1$ then
  \begin{equation}
\label{E:dagger}
  v_{\tnulam}y_{n+g_0}\dots y_{n+g_{z}} \equiv \Big( v_{\t_k}
  + \sum_{\u \in \Stdlm} a_\u v_{\u} \Big) y_{n+g_k}\dots y_{n+g_{z}}
       \pmod{\hat{S}^\bnu_\bmu},
  \end{equation}
  for some $a_\u \in \Z$. By definition, $v_{\t_{z+1}} \in \hat{S}^\bnu_{\tnumu}$
  and $v_{\t_{z+1}} \notin \hat{S}^\bnu_\bmu$. Moreover,
  $v_\u \in \hat{S}^\bnu_{\tnumu}$ for all $\u \in \Stdlm[z+1]$.  Therefore,
  if $k=z+1$ then \autoref{E:dagger} is the statement
  $$\thea(v_{\tnulam})\equiv v_{\t_{z+1}}\pmod{\hat{S}^\bnu_{\tnumu}}.$$
  So, taking $k=z+1$ in \autoref{E:dagger} establishes (B) and (C) and
  completes the proof.

  To prove that \autoref{E:dagger} holds we argue by induction on~$k$. If $k=0$
  then there is nothing to prove because we can set $a_\u=0$ for all
  $\u\in\Stdlm[0]$. Hence, by induction it is enough to show that
  if \autoref{E:dagger} holds for $k+1$ whenever it holds for~$k$.  First,
  suppose that $\u\in\Stdlm$ and that $a_\u\ne0$ in \autoref{E:dagger}.
  Using~\autoref{L:ReducedExps}, we can write
  $v_\u = v_\s \psi_w$ where $\s \in \Std^e_n(\bnu)$ and $w \in \Sym_n$.
  By~\autoref{L:Pushing} and~\autoref{L:Filtration}, $v_\u y_{n+g_k} = v_\s \psi_w
  y_{n+g_k} = v_\s y_{n+g_k} \psi_w$ is a $\Z$-linear combination of terms
  $v_{\t}$ where $\t\in\Std(\bnu)$ and $\row_{\t}(n+1)>r_k$. Therefore,
  $\row_\t(n+1)\ge r_{k+1}$ since $n+1$ can only appear in a removable $j$-strip
  of length at most~$\gamma$.  If $\row_\t(n+1)>r_{k+1}$ then $v_\t
  y_{n+g_{k+1}}\dots y_{n+g_z}\in \hat{S}^\bnu_\bmu$ by
  \autoref{L:PushI}. Alternatively, if $\row_\t(n+1)=r_{k+1}$ then
  $\t_{z+1}\not\gedom\t$, because $\t_{z+1}\not\gedom\u$ and $\t\gdom\u$ by
  \autoref{L:YDown}. Therefore, $\t\in\Stdlm[k+1]$.

  Now consider the leading term $(v_{\t_k}y_{n+g_k})y_{n+g_{k+1}}\dots y_{n+g_z}$
  appearing in~\autoref{E:dagger}. By \autoref{P:yngI}, $v_{\t_k}y_{n+g_k}$ is
  equal to $v_{\t_{k+1}}$ plus a $\Z$-linear combination of terms $v_\t$ such
  that $\row_\t(n+1)\ge r_{k+1}$ and either $\t_{z+1}\not\gedom\t$, or there
  exists a $g$ such that $\row_\t(n+g)>\row_{\t_{k+1}}(n+g)$ and
  $1\le g\le g_k$. In the latter case, $v_\t y_{n+g_{k+1}} \ldots y_{n+g_z} \in
  \hat{S}^\bnu_\bmu$ by \autoref{L:PushI}. On the other hand, if
  $\row_\t(n+1)=r_{k+1}$ and $\t_{z+1}\not\gedom\t$ then $\t\in\Stdlm[k+1]$.

  We have now shown that \autoref{E:dagger} holds for $k+1$ and so proved our claim.
  Therefore, by \autoref{L:EllersMurray}, $\thea$ induces a non-zero
  homogeneous map $S^\blam\<\delta\>\to S^\bmu$, where
  $\delta=(\deg\tnulam-\deg\tlam)-(\deg\tnumu-\deg\tmu)+\deg L^\blam_\bmu$.
  Arguing as in the last paragraph of proof of \autoref{T:OneNode} shows that
  $\delta=a-b+2c$, where $a$,~$b$ and~$c$ are as in the statement of
  \autoref{T:OneRow}. We leave these details to the reader.
\end{proof}

\subsection{Cyclotomic Carter-Payne pairs} \label{S:ProofMainThm}
We now turn to the statement and proof of our Main Theorem.  Special cases of
our main result include \autoref{T:OneNode} and \autoref{T:OneRow}, both of
which motivate our main result and simplify its proof.

If $\blam$ and $\bnu$ are multipartitions and $\blam\subset\bnu$ let
$\res(\bnu\backslash\blam)=\set{\res(\alpha)|\alpha\in\bnu\backslash\blam}$.
If $J \subseteq I$ and $i,j\in J$ write $j\le_J i$ if
$\{j,j+1,\ldots,i\} \subseteq J$.

\begin{Definition}
\label{D:CycCPPair}
  Suppose that $\blam$ and $\bmu$ are multipartitions of $n\ge0$
  and let
  $\bnu=\blam\cup\bmu$ and $J=\res(\bnu\backslash\blam)$. Then
  $(\blam,\bmu)$ is a \textbf{cyclotomic Carter-Payne pair} if
  $J=\res(\bnu\backslash\bmu)$, $|J|=|\bnu\backslash\blam|$ and $|J|<|I|$,
  and if $i,j\in J$ and $j\le_J i$ then
  $$
    \row\alpha^\bnu_\blam(j)\le\row\alpha^\bnu_\blam(i),\quad
    \row\alpha^\bnu_\bmu(j)=\row\alpha^\bnu_\bmu(i) \quad\text{and}\quad
    \row\alpha^\bnu_\blam(i)\le\row\alpha^\bnu_\bmu(i),
  $$
  and $\alpha^\bnu_\blam(i)$ is in a removable $j$-strip.
  Here, and below, if $k\in J$ then $\alpha^\bnu_\blam(k)\in\bnu\backslash\blam$ and
  $\alpha^\bnu_\bmu(k)\in\bnu\backslash\bmu$ are the unique
  nodes in~$\bnu$ such that $\res\alpha^\bnu_\blam(k) = k = \alpha^\bnu_\bmu(k)$.
\end{Definition}

Note that in \autoref{D:CycCPPair} the residues in the set $J$ are necessarily
distinct because $|J|=|\bnu\backslash\blam|$. It is easy to see that if $\blam$
and $\bmu$ are multipartitions satisfying the assumptions of \autoref{T:OneRow}
then $(\blam,\bmu)$ is a cyclotomic Carter-Payne pair.

\shadedBaseline=-1mm
\begin{Example} \label{Ex:General}
Suppose that $e=0$ and let $\kappa=(0,4,0,0,4)$.  Let $\blam =
(\emptyset|\emptyset|2|3|1)$ and $\bmu=(2|1|3|\emptyset|\emptyset)$ so that
$\bnu=\blam\cup\bmu=(2|1|3|3|1)$.  In the first diagram below, we have shaded the nodes
$\alpha^\bnu_{\blam}(i)$ and in the second we have shaded the nodes
$\alpha^\bnu_{\bmu}(i)$, for $i \in J = \{0,1,2,4\}$.
\begin{align*}
  \alpha^\bnu_\blam(i):&&
  \Big(\ \ShadedTableau[(1,0),(2,0)]{{0,1}}\ \Big|\
          \ShadedTableau[(1,0)]{{4}}\ \Big|\
          \ShadedTableau[(3,0)]{{0,1,2}}\ \Big|\
          \ShadedTableau[]{{0,1,2}}\ \Big|\
          \ShadedTableau[]{{4}}\
\Big),\\
  \alpha^\bnu_\bmu(i):&&
\Big(\ \ShadedTableau[]{{0,1}}\ \Big|\
          \ShadedTableau[]{{4}}\ \Big|\
          \ShadedTableau[]{{0,1,2}}\ \Big|\
          \ShadedTableau[(1,0),(2,0),(3,0)]{{0,1,2}}\ \Big|\
          \ShadedTableau[(1,0)]{{4}}\
\Big).
\end{align*}
Hence,  $(\blam,\bmu)$ is a cyclotomic Carter-Payne pair.
\end{Example}

The aim of this section is to prove the following result, which is the most
general statement of our Main Theorem from the introduction. To state this
result we continue to use the notation from \autoref{D:CycCPPair}.

\begin{Theorem}
\label{T:main}
  Suppose that $(\blam,\bmu)$ is a cyclotomic Carter-Payne pair and let
  $J=\res(\bnu\backslash\blam)$ and $J^*=\set{j\in J|j-1\notin J}$. Then
  $$\Hom_{\R}(S^\blam\<a-b+2c-d\>,S^\bmu)\ne0,$$
  where for $j\in J^*$ we set $\gamma_j=\#\set{i\in J|j\le_J i}$ and
  \begin{align*}
    a&=\#\set{\alpha\in\Add(\bnu)|\res(\alpha)=i\in J\text{ and }
        \row\alpha^\bnu_\blam(i)<\row\alpha\leq\row\alpha^\bnu_\bmu(i)},\\
    b&=\#\set{\alpha\in\Rem(\bnu)|\res(\alpha)=i\in J\text{ and }
        \row\alpha^\bnu_\blam(i)<\row\alpha\leq \row\alpha^\bnu_\bmu(i)},\\
    c&=\#\Set[80]{C}{$C$ is a removable $j$-strip of length at most $\gamma_j$ such
                 that $\row\alpha^\bnu_\blam(j)<\row C\leq \row\alpha^\bnu_\bmu(j)$,
                 for $j\in J^*$},\\
    d&=\#\set{\row\alpha^\bnu_\blam(i)|i\in J}
              -\#\set{\row\alpha^\bnu_\bmu(i)|i\in J}.
  \end{align*}
  Moreover, $a-b+2c-d>0$.
\end{Theorem}

When $e=2$,~\autoref{T:main} reduces to \autoref{T:OneNode} (and
\autoref{T:OneRow}). When $e=3$,~\autoref{T:main} is slightly stronger than
\autoref{T:OneRow} because the nodes in~$\bnu\backslash\blam$ are not
necessarily in the same row.

By \autoref{D:CycCPPair}, if $j\le_J i$ then
$\row\alpha^\bnu_\bmu(i)=\row\alpha^\bnu_\bmu(j)$ but $\alpha^\bnu_\blam(i)$ and
$\alpha^\bnu_\blam(j)$ can be in different rows.
Therefore, the integer~$d$ in \autoref{T:main} is the sum over $j\in J^*$ of the
differences between the number of rows occupied by the nodes
$\set{\alpha^\bnu_\blam(i)|j\le_J i}$ and the nodes
$\set{\alpha^\bnu_\bmu(i)|j\le_J i}$.
Every row that contributes to this sum has, for some $j \in J^*$, a removable $j$-strip $C$ of length at most~$\gamma_j$ where $\row\alpha^\bnu_\blam(j)<\row C\leq \row\alpha^\bnu_\bmu(j)$, that is it also contributes towards the integer $c$ (but not the integer $b$).  The argument given in the proof of \autoref{T:OneRow} then shows that $a-b+2c-d>0$.

By the last paragraph, if $\bnu\backslash\blam$ is
contained in one row  then $d=0$ and \autoref{T:main} reduces to
\autoref{T:OneRow}. More generally, as we will see, the proof of
\autoref{T:main} is almost identical to the proof of \autoref{T:OneRow}. For
this reason, we only sketch the proof of \autoref{T:main}, highlighting
those places where changes to the previous argument are needed.

We now fix a cyclotomic Carter-Payne pair $(\blam,\bmu)$ and turn to the proof
of \autoref{T:main}.  Note that if $e=2$ then $(\blam,\bmu)$ is a cyclotomic
Carter-Payne pair only if $\gamma=1$, in which case~\autoref{T:main} reduces to
\autoref{T:OneNode}. We can therefore assume for the rest of this section
that $e\neq 2$.

To prove \autoref{T:main} we generalise the construction of the tableaux $\t_k$
and the multipartitions $\bsig_k$, for $0\le k\le z+1$, which were used in the
proof of \autoref{T:OneRow}. We invite the reader to check that these two
notations agree in the special case when the pair $(\blam,\bmu)$ satisfies the
assumptions of \autoref{T:OneRow}.

As in \autoref{T:main}, set $\bnu=\blam\cup\bmu$,
$J=\res(\bnu\backslash\blam)$ and $J^*=\set{j\in J|j-1\notin J}$. If $j\in J^*$
define $\gamma_j=\#\set{i\in J|j\le _J i}$ and let $\gamma=\sum_{j\in J^*}\gamma_j$.
Then $\gamma=|\bnu\backslash\blam|$.

Suppose that $\tnulam$ is a standard $\bnu$-tableau such that
$(\tnulam)_{\downarrow n} = \tlam$ and if $i,i+1 \in J$
and $\row \alpha^\bnu_{\blam}(i) = \row \alpha^\bnu_\blam(i+1)$
then
$\tnulam(\alpha^\bnu_\blam(i+1)) = \tnulam(\alpha^\bnu_\blam(i))+1$, and
let
$$\IniNu[\blam]=\set{\t\in\IniNu|\res_\t(n+g)=\res_{\tnulam}(n+g),
               \text{ for }1\le g\le\gamma}.$$
The choice of $\tnulam$ does affect our definition of the $\R$-modules
$\check{S}^\bnu_{\tnulam}, \hat{S}^\bnu_{\tnumu} \subset S^\bnu$, and
consequently the definition of the $\R$-endomorphism $\thea:S^\bnu \rightarrow
S^\bnu$, defined in \autoref{E:thetaMain} below, but it does not affect the
statement of the Main Theorem.

As before, if $\bsig\in\Parts$ is a multipartition such that
$J=\res(\bnu\backslash\bsig)$ and $\bsig\subset\bnu$
define $\tnu_\bsig=\tnu_\bsig(\blam)$ to be the unique standard $\bnu$-tableau
in~$\IniNu[\blam]$ such that $(\tnu_{\bsig})_{\downarrow n}=\t^\bsig$.

Mirroring the definitions in \autoref{S:OneRowCP}, let $r$ be the smallest row
index, with respect to the lexicographic order~\autoref{E:RowOrder}, such that the
corresponding rows of $\bnu$ and $\blam$ are different and let $s$ be the
largest row index such that $\bnu$ and $\bmu$ differ.  Let
$C_0,C_1,\dots,C_{z+1}$ be the complete list of the removable $j$-strips of
length at most~$\gamma_j$, for some $j\in J^*$, lying between rows $r$ and $s$
in the diagram of~$\bnu$ such that $r=r_0<r_1<\dots<r_{z+1}=s$, where
$r_k=\row(C_k)$ for $0\le k\le z+1$. For $0\le k\le z+1$ define $j_k\in J^*$
to be the residue such that $C_k$ is a
removable $j_k$-strip and let $f_k$ be the unique integer such that
$\res_{\tnulam}(n+f_k)=j_k$ and $1\le f_k\le\gamma$. Almost exactly
as before, define standard tableaux
$\t_0,\t_1,\dots,\t_{z+1}$ recursively by letting $\t_0=\tnu_\blam$ and setting
$$\t_{k+1}=\t_k\prod_{h=f_k}^{g_k}(\eta^{(k)}_h,n+h),$$
where $f_k\le g_k\le\gamma$ and $g_k$ is maximal such that
$\row_{\t_k}(n+g_k)=r_k$ and the
integer~$\eta^{(k)}_h$ is minimal such that
$\res_{\t_k}(\eta^{(k)}_h)=\res_{\t_k}(n+h)$ and
$$\t_k^{-1}(\eta^{(k)}_h)\in C_{k+1}\cup C_{k+2}\cup\dots\cup C_{z+1},$$
for $f_k\le h\le g_k$. By construction, $\row_{\t_k}(n+f_k)=r_k$, for
$0\le k\le z+1$.

Generalising \autoref{E:thetaOneRow}, define $\thea$ to be the $\R$-endomorphism
of $S^\bnu$ given by
\begin{equation}
\label{E:thetaMain}
  \thea(x)=xy_{n+g_0}\dots y_{n+g_z}e^\blam_\mu, \qquad\text{ for }x\in S^\bnu,
\end{equation}
where $e^\blam_\bmu=\sum_{\bj\in I^n}e(j_1,\dots,j_n,i_1,\dots,i_\gamma)$ and
$i_h=\res_{\tnulam}(n+h)$ for $1\le h\le\gamma$.
By definition, $\thea$ is an $\R$-module endomorphism of~$S^\bnu$.

\begin{Lemma}
\label{L:ThetaIndependence}
  Suppose that $(\blam,\bmu)$ is a cyclotomic Carter-Payne pair. Then $\thea$
  depends only on $\blam$, $\bmu$ and $\tnulam$.
\end{Lemma}

\begin{proof}The tableau $\tnulam$ determines the idempotent $e^\blam_\bmu$. For
  $1\le g\le\gamma$ let $i_g=\res_{\tnulam}(n+g)$. If $x\in S^\bnu$ then it is
  easy to see that
  $$\thea(x)=x\prod_{g=1}^\gamma y_{n+g}^{c_g}e^\blam_\bmu,$$
  where $c_g$ is equal to the number of removable $j$-strips $C$ of length at
  most $\gamma_j$ such that $j=i_g-|C|+1\in J^*$ and
  $\row_{\tnulam}(n+g)=\row\alpha^\bnu_\blam(i_g)<\row
  C\le\row\alpha^\bnu_\bmu(i_g)$. (The condition $i_g=j+|C|-1$ simply says that
  the rightmost node in $C$ has residue~$i_g$.) This implies the result.
\end{proof}

To show that $\thea$ induces a non-zero homomorphism
from $S^\blam\<a-b+2c-d\>$ to $S^\bmu$ we need some more
notation.  For $0 \leq k \leq z+1$ define multipartitions
$\bsig_k=\Shape((\t_k)_{\downarrow n})$ and define permutations $w_k\in\Sym_n$ by
$(\t_k)_{\downarrow n}=(\tnu_{\bsig_k})_{\downarrow n}w_k$. Then
$\t_{z+1}\gdom\t_z\gdom\dots\gdom\t_1\gdom\t_0=\tnulam$.
Exactly as in~\autoref{E:lm}, define integers
$l_h^{(k)}$ and $m_h^{(k)}$, for $f_k\le h\le g_k$ and $0\le k\le z$ and hence
define $\psi_{w_k}$ as in~\autoref{E:vtk}.

\shadedBaseline=-1mm
\begin{Example}
\label{Ex:GeneralCont}
As in~\autoref{Ex:General}, let $e=0$ and let $\kappa=(0,4,0,0,4)$.  Consider
$\blam = (\emptyset|\emptyset|2|3|1)$ and $\bmu=(2|1|3|\emptyset|\emptyset)$.  The
residues in $\bnu=\blam\cup\bmu=(2|1|3|3|1)$ are
$$\Big(\ \ShadedTableau[]{{0,1}}\ \Big|\
          \ShadedTableau[]{{4}}\ \Big|\
          \ShadedTableau[]{{0,1,2}}\ \Big|\
          \ShadedTableau[]{{0,1,2}}\ \Big|\
          \ShadedTableau[]{{4}}\
\Big).$$
By definition, $\gamma=|\bnu\backslash\blam|=4$ and $z=2$. The tableaux $\t_k$,
for $0\le k\le z+1$ are
\begin{align*}
\tnulam=\t_0&=\Big(\ \ShadedTableau[]{{7,8}}\ \Big|\
          \ShadedTableau[]{{10}}\ \Big|\
         \ShadedTableau[]{{1,2,9}}\ \Big|\
         \ShadedTableau[]{{3,4,5}}\ \Big|\
          \ShadedTableau[]{{6}}\
\Big), \\
\t_1&=\Big(\ \ShadedTableau[(1,0),(2,0)]{{1,2}}\ \Big|\
          \ShadedTableau[]{{10}}\ \Big|\
          \ShadedTableau[(1,0),(2,0)]{{7,8,9}}\ \Big|\
          \ShadedTableau[]{{3,4,5}}\ \Big|\
          \ShadedTableau[]{{6}}\
\Big),\\
\t_2&=\Big(\ \ShadedTableau[]{{1,2}}\ \Big|\
         \ShadedTableau[(1,0)]{{6}}\ \Big|\
         \ShadedTableau[]{{7,8,9}}\ \Big|\
         \ShadedTableau[]{{3,4,5}}\ \Big|\
         \ShadedTableau[(1,0)]{{10}}\
\Big),\\
\t_3&=\Big(\ \ShadedTableau[]{{1,2}}\ \Big|\
         \ShadedTableau[]{{6}}\ \Big|\
         \ShadedTableau[(1,0),(2,0),(3,0)]{{3,4,5}}\ \Big|\
         \ShadedTableau[(1,0),(2,0),(3,0)]{{7,8,9}}\ \Big|\
         \ShadedTableau[]{{10}}\
\Big),\\
\end{align*}
where at each step we have shaded the entries which moved. In the notation of
\autoref{T:main}, $a=0$, $b=2$, $c=3$ and $d=1$, so the theorem implies that there
exists a non-zero homomorphism $S^\blam\<3\>\to S^\bmu$. This map is induced by
the $\R$-module endomorphism $\thea$ of~$S^\bnu$ which is given by
$\thea(x)=xy_8y_9y_{10}e^\blam_\bmu$, for $x\in S^\bnu$,
where $e^\blam_\bmu=\sum_{\bi\in I^6}e(i_1,\dots,i_6,0,1,2,4)$.
\end{Example}

We now refer the reader to~\autoref{L:PushI}.  It is clear that the statement given in that lemma also holds for our more generalised setup.

As in~\autoref{Ex:GeneralCont}, there may be tableaux $\t_k$ for
$1 \leq k \leq z+1$ in which $\row_{\t_k}(n+f) > \row_{\t_k}(n+g)$
for some $1 \leq f<g \leq \gamma$.   As a consequence, certain words in
$\psi_1,\dots,\psi_{n+\gamma-1}$ which previously corresponded to reduced
expressions in $s_1,\ldots,s_{n+\gamma-1}$ no longer necessarily have this property.  As
we shall see, however, this will ultimately not cause us any difficulties.

\begin{Example} \label{Ex:GeneralCont2}
We continue~\autoref{Ex:General} (and~\autoref{Ex:GeneralCont}).  The braid diagrams below show
that $v_{\t_1} y_{10} = v_{\t_2}$.
$$
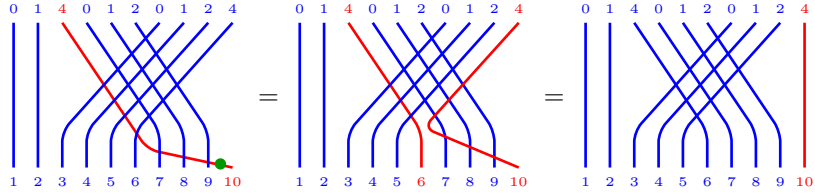
\begin{braid}\tikzset{xscale=0.8,yscale=1.2,baseline=24}
\draw(1,4)node[above ]{0}--(1,0)node[below]{1};
\draw(2,4)node[above ]{1}--(2,0)node[below]{2};
\draw(3,4)[red]node[above ]{4}--(6.5,0.5)--(10,0)node[below]{10};
\draw(4,4)node[above ]{0}--(7,1)--(7,0)node[below]{7};
\draw(5,4)node[above ]{1}--(8,1)--(8,0)node[below]{8};
\draw(6,4)node[above ]{2}--(9,1)--(9,0)node[below]{9};
\draw(7,4)node[above ]{0}--(3,1)--(3,0)node[below]{3};
\draw(8,4)node[above ]{1}--(4,1)--(4,0)node[below]{4};
\draw(9,4)node[above ]{2}--(5,1)--(5,0)node[below]{5};
\draw(10,4)node[above ]{4}--(6,1)--(6,0)node[below]{6};
\greendot(9.5,0.1);
\end{braid} =
\begin{braid}\tikzset{xscale=0.8,yscale=1.2,baseline=24}
\draw(1,4)node[above ]{0}--(1,0)node[below]{1};
\draw(2,4)node[above ]{1}--(2,0)node[below]{2};
\draw(3,4)[red]node[above ]{4}--(6,1)--(6,0)node[below]{6};
\draw(4,4)node[above ]{0}--(7,1)--(7,0)node[below]{7};
\draw(5,4)node[above ]{1}--(8,1)--(8,0)node[below]{8};
\draw(6,4)node[above ]{2}--(9,1)--(9,0)node[below]{9};
\draw(7,4)node[above ]{0}--(3,1)--(3,0)node[below]{3};
\draw(8,4)node[above ]{1}--(4,1)--(4,0)node[below]{4};
\draw(9,4)node[above ]{2}--(5,1)--(5,0)node[below]{5};
\draw(10,4)[red]node[above ]{4}--(6.1,1.1)--(10,0)node[below]{10};
\end{braid} =
\begin{braid}\tikzset{xscale=0.8,yscale=1.2,baseline=24}
\draw(1,4)node[above ]{0}--(1,0)node[below]{1};
\draw(2,4)node[above ]{1}--(2,0)node[below]{2};
\draw(3,4)node[above ]{4}--(6,1)--(6,0)node[below]{6};
\draw(4,4)node[above ]{0}--(7,1)--(7,0)node[below]{7};
\draw(5,4)node[above ]{1}--(8,1)--(8,0)node[below]{8};
\draw(6,4)node[above ]{2}--(9,1)--(9,0)node[below]{9};
\draw(7,4)node[above ]{0}--(3,1)--(3,0)node[below]{3};
\draw(8,4)node[above ]{1}--(4,1)--(4,0)node[below]{4};
\draw(9,4)node[above ]{2}--(5,1)--(5,0)node[below]{5};
\draw(10,4)[red]node[above ]{4}--(10,0)node[below]{10};
\end{braid}$$
\end{Example}

\begin{Lemma} \label{L:StringFamilies}
  Suppose that $\varpi$ is a word in $\Psin[n+\gamma]$
and let $B_\varpi$ be the braid diagram for $v_{\tnu} \varpi$.
Suppose that $B_\varpi$ has the following properties:
\begin{itemize}
\item If $1 \leq m <m' \leq n$, the $m$-string and the $m'$-string do not cross each other.
\item If $1 \leq m \leq n$ and $1 \leq g \leq \gamma$, the $m$-string and the $n+g$-string cross at most once.
\item If $1 \leq g,h \leq \gamma$ and the $n+g$-string and the $n+h$-string cross, then $i_g \neq i_h \pm 1$.
\end{itemize}
Let $s_{j_1}s_{j_2}\ldots s_{j_l}$ be a reduced expression for $\pi_{\varpi}$.  Then
$v_{\tnu} \varpi =v_{\tnu} \psi_{j_1}\psi_{j_2}\ldots\psi_{j_l}$.
\end{Lemma}

\begin{proof}
  Suppose $\varpi=\psi_{k_1}\psi_{k_2}\ldots\psi_{k_{l'}}$.
  As permutations, we may get from (the possibly non-reduced expression)
  $s_{k_1}s_{k_2} \ldots s_{k_{l'}}$ to any reduced expression for $\pi_{\varpi}$ by applying
  the braid relations and the relation $s_j^2=1$ for $1 \leq j<n+\gamma$.
The result follows, once we have checked that when making the corresponding manipulations of the diagram $B_\varpi$, we never get local diagrams of the form
$$\begin{braid}[7]
    \draw(1,4)node[above]{$i$}--(3,0);
    \draw(2,4)node[above]{$i{\pm}1$}--(3,2)--(2,0);
    \draw(3,4)node[above]{$i$}--(1,0);
  \end{braid},
\qquad \qquad
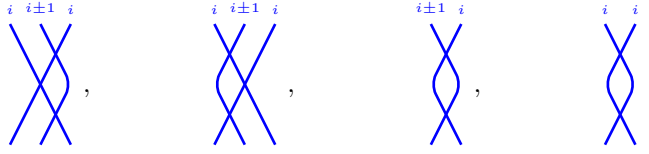
\begin{braid}[7]
    \draw(1,4)node[above]{$i$}--(3,0);
    \draw(2,4)node[above]{$i{\pm}1$}--(1,2)--(2,0);
    \draw(3,4)node[above]{$i$}--(1,0);
  \end{braid},
\qquad \qquad
\begin{braid}[7]
    \draw(2,4)node[above]{$i{\pm}1$}--(3,2)--(2,0);
    \draw(3,4)node[above]{$i$}--(2,2)--(3,0);
  \end{braid},
\qquad \qquad
\begin{braid}[7]
    \draw(2,4)node[above]{$i$}--(3,2)--(2,0);
    \draw(3,4)node[above]{$i$}--(2,2)--(3,0);
  \end{braid},
$$
but this follows easily from our assumptions on $B_\varpi$.
\end{proof}

As in \autoref{E:vtk}, set $v_{\t_k}=v_{\tnu_{\bsig_k}}\psi_{w_k}$, for
$0\le k\le z+1$, and consider the generalisation of \autoref{P:yngI} to the
current setting. We claim that if $0\le k\le z$ then
\begin{equation} \label{E:yngI}
v_{\t_k}y_{n+g_k}=v_{\t_{k+1}} +\sum_{\s\in\Std(\bnu)}a_\s v_\s,
\end{equation}
where $a_\s\ne0$ only if $\res(\s)=\res(\tnulam)$, $\row_\s(n+f_k)\ge \row_{\t_{k+1}}(n+f_k)$ and
either there exists $f_k \leq g\leq g_k$ with $\row_\s(n+g)>\row_{\t_{k+1}}(n+g)$
or $\t_{z+1}\not\gedom\s$ (or both).

To show that~\autoref{E:yngI} holds, we follow the proof of~\autoref{P:yngI} and write $v_{\t_k} y_{n+g_k}=\sum_{\u \in \Std(\bnu)}b_\u v_\u$ where $b_\u \neq 0$ only if $\res(\u)=\res(\tnulam)$ and $\row_\u(n+f_k)\geq \row_{t_{k+1}}(n+f_k)$.  Again, we must show that $b_{\t_{k+1}}=1$ and that if $\u \neq \t_{k+1}$ and $b_\u \neq 0$ then either there exists $f_k \leq g \leq g_k$ such that $\row_\u(n+g)>\row_{t_{k+1}}(n+g)$ or $\t_{z+1} \not \gedom \u$ (or both).

Using \autoref{L:GeneralSigma}, $v_{\tnu_{\bsig_k}} y_{n+g_k}$ is equal to the
sum of diagrams of the form below, where the sum is over all integer sequences
$\m=(m_{f_k},\ldots,m_{g_k})$ such that
$1 \leq m_{f_k}<m_{f_{k}+1}<\ldots<m_{g_k}\le n$,
$\res_{\tnu_{\bsig_k}}(m_h)=\res_{\tnu_{\bsig_k}}(n+h)$ and
$\row_{\tnu_{\bsig_k}}(m_h)>\row_{\tnu_{\bsig_k}}(n+h)$, for $f_k\le h\le g_k$.
$$
\begin{braid}\tikzset{xscale=0.8,yscale=1.2}
  \draw[brown](0,4.8)--(4.7,4.8)--(4.7,4)--(0,4);
  \draw[brown,dashed](-1,4.8)--(0,4.8);
  \draw[brown,dashed](-1,4)--(0,4);
  \foreach \k in {12,12.5,13} {
    \draw[blue!50](\k,4)--+(-5,-5);
  }
  \foreach \k in {0,4} {
    \draw[red](\k,4)--(\k+12.5,2.5)--(\k+9,-1);
    \draw[blue!50](\k+8.5,-1)--+(5,5);
    \draw[blue!50](\k+9.5,-1)--+(5,5);
  }
  \draw[red](9,-1)node[below]{$m_{f_k}$};
  \draw[dots](10.0,-0.9)--(11.3,-0.9);
  \draw[red](13,-1)node[below]{$m_{g_k}$};
  \draw(14,4)--(13.5,3.5)--(29,0)--(29,-1)node[below]{$n_{f_k}$};
  \draw[dots](29.3,-0.9)--(30.9,-0.9);
  \draw(18,4)--(17.2,3.2)--(31,0)--(31,-1)node[below]{$n_{g_k}$};
  \draw(0,3.9)node[above]{$i_{f_k}$};
  \draw[dots](0.9,4.4)--(3.3,4.4);
  \draw(4,3.9)node[above]{$i_{g_k}$};
  \draw[blue!50](12,-1)--+(5,5);
  \draw[blue!50](11.5,-1)--+(5,5);
  \draw[densely dashdotted,darkgreen!60](8,4)--(26,0)--(33,-1);
  \draw[densely dashdotted,darkgreen!60](15,4)--(32.5,0)--(33.5,-1);
  \draw[densely dashdotted,darkgreen!60](16,4)--(33.5,0)--(27.5,-1);
\end{braid}
$$
In all of these diagrams the $(n+h)$-strings, for $f_k\le h\le g_k$, correspond
to rows below $\row_{\tnu_{\bsig_k}}(n+f_k)$. Fix $l$ with $1\le l<f_k$ or
$g_k<l\le\gamma$. The dashed lines in the diagram indicate that
the $(n+l)$-string is in exactly the same position as it is in the diagram for
$v_{\tnu_{\bsig_k}}$. The $(n+l)$-string plays no role in the argument below
because, by assumption, the residue of this string is not equal to the residue
of an $(n+h)$-string, where $f_k\le h\le g_k$.

For each such $\m=(m_{f_k},\ldots,m_{g_k})$, define $\varpi_\m \in \W$ and $B_\m$ as in~\autoref{P:yngI}.
First consider the case when $m_h=m_h^{(k)}$ for $f_k \leq h \leq g_k$.
Then, considered as a word,
$$\varpi_\m=\varpi\prod_{h=f_k}^{g_k}
           \psi_{s(l^{(k)}_{f_k+g_k-h},m^{(k)}_{f_k+g_k-h})}$$
where $\varpi$ satisfies the assumptions of~\autoref{L:StringFamilies}.
Furthermore $\tnu \pi_{\varpi} = v_{\tnu_{\bsig_{k+1}}}$, so that $v_{\tnu}
\varpi = v_{\tnu_{\bsig_{k+1}}}$, regardless of the choice of reduced expression
that was used to define $v_{\tnu_{\bsig_{k+1}}}$.  Therefore, as in the proof
of~\autoref{P:yngI}, we obtain
\[v_{\tnu} \varpi_\m \psi_{w_k} = v_{\tnu_{\bsig_{k+1}}} \prod_{h=f_k}^{g_k}
    \psi_{s(l^{(k)}_{f_k+g_k-h},m^{(k)}_{f_k+g_k-h})} \psi_{w_k} = v_{\t_{k+1}}.\]

Now suppose that there exists $f_k \leq g \leq g_k$ such that $m_g>m^{(k)}_g$.
Exactly as in the proof of~\autoref{P:yngI}, $\tnu \varpi_m \psi_k$ is a linear
combination of terms $v_\u$ where $\u \in \Std(\bnu)$ is such that
$\row_\u(n+g)>\row_{\t_{k+1}}(n+g)$.

Finally suppose that $m_h \leq m_h^{(k)}$, for $1 \leq h \leq g_k$, and that
there exists an integer~$g$ such that $f_k \leq g \leq g_k$ and $m_g<m_g^{(k)}$.
As in the case where $m_h=m_h^{(k)}$ for all $h$, we
use~\autoref{L:StringFamilies} to write $v_{\tnu} \varpi_\m = v_{\tnu}
\psi_{d(\pi_{\varpi_\m})}$.  The argument used in the proof
of~\autoref{P:yngI} shows that $v_{\tnu}\varpi_\m \psi_{w_k}$ is a linear
combination of terms~$v_\u$ where $\u \in \Std(\bnu)$ and $\t_{z+1}\not\gedom\u$.
This completes the proof of~\autoref{E:yngI}.

This establishes the required generalisation of~\autoref{L:PushI}
and~\autoref{P:yngI}.  Now the argument used to prove \autoref{T:OneRow} can be
repeated, word for word, to complete the proof of \autoref{T:main}.

\subsection{Maps involving many cyclotomic Carter-Payne pairs} In this last
section we give two applications of \autoref{T:main} to situations which
involve more than one cyclotomic Carter-Payne pair. Again, we assume that
$e\in\{0,2,3,4,5,\dots\}$ is arbitrary.

If $(\blam,\bmu)$ is a cyclotomic Carter-Payne pair, define
$$\delta^\blam_\bmu=a-b+2c-d,$$
where $a$, $b$, $c$ and $d$ are as in \autoref{T:main}, and let
$\Theta^\blam_\bmu\in\Hom_{\R}(S^\blam\<\delta^\blam_\bmu\>,S^\bmu)$ be the
$\R$-module homomorphism induced by the map~$\theta^\blam_\bmu$ which was
defined in \autoref{E:thetaMain}.

The first result in this section is almost implicit in the proof of
\autoref{T:main}. In the statement of the next result we adopt the convention
that if $X$ and $Y$ are sets then $X\backslash Y=\set{x\in X|x\notin Y}$.

\begin{Theorem}
\label{T:composition}
  Suppose that $\blam,\bmu,\brho\in\Parts$ are multipartitions of~$n$ such that
  $\blam\backslash\brho=\bmu\backslash\brho$ and
  $(\blam,\brho)$ and $(\brho,\bmu)$ are cyclotomic Carter-Payne pairs.
  Then
  $(\blam,\bmu)$ is a cyclotomic Carter-Payne pair,
  $\delta^\blam_\bmu=\delta^\blam_\brho+\delta^\brho_\bmu\ge2$ and there exist maps
  $\Theta^\blam_{\brho}\in\Hom_{\R}(S^\blam\<\delta^\blam_\brho\>,S^\brho)$ and
  $\Theta^\brho_\bmu\in\Hom_{\R}(S^\brho\<\delta^\brho_\bmu\>,S^\bmu)$ such that
  $$0\ne\Theta^\blam_\bmu=\Theta^\brho_\bmu\circ\Theta^\blam_\brho
                 \in\Hom_{\R}(S^\blam\<\delta^\blam_\bmu\>,S^\bmu).$$
\end{Theorem}

Before we prove the theorem we give an example to illustrate it.

\begin{Example} \label{Ex:Composition}\shadedBaseline=-4mm
Suppose that $e=5$ and $\kappa = (4,1,0,0,2)$ and let
\[ \blam = (6,4,1|3,1|1|1|4,3), \quad \brho = (7,5,2|1,1|1|\emptyset|4,3) \quad\text{and}\quad
\bmu = (7,5,2|3,1|1|1|3,1).\]
Then $(\blam, \brho)$ and $(\brho,\bmu)$ are cyclotomic Carter-Payne pairs, so
\autoref{T:composition} implies that $(\blam,\bmu)$ is also a cyclotomic
Carter-Payne pair, which is easily checked. The residue diagram of $\bnu = \blam
\cup \bmu (= \blam \cup \brho = \brho \cup \bmu) = (7,5,2|3,1|1|1|4,3)$ is given
by
$$\Bigg(\ \ShadedTableau[]{{4,0,1,2,3,4,0},{3,4,0,1,2},{2,3}}\ \Bigg|\
          \ShadedTableau[(3,0),(2,0)]{{1,2,3},{0}}\ \Bigg|\
          \ShadedTableau[]{{0}}\ \Bigg|\
          \ShadedTableau[(1,0)]{{0}}\  \Bigg|\
          \ShadedTableau[]{{2,3,4,0},{1,2,3}}\
\Bigg),$$
where we have shaded the nodes in $\blam \backslash \brho = \bmu \backslash \brho$.  Let
\begin{align*}
\tnulam&= \Bigg(\ \ShadedTableau[]{{1,2,3,4,5,6,25},{7,8,9,10,26},{11,27}}\ \Bigg|\
          \ShadedTableau[]{{12,13,14},{15}}\ \Bigg|\
          \ShadedTableau[]{{16}}\ \Bigg|\
          \ShadedTableau[]{{17}}\  \Bigg|\
          \ShadedTableau[]{{18,19,20,21},{22,23,24}}\
\Bigg), \text{ and } \\
\tnu_\brho&= \Bigg(\ \ShadedTableau[]{{1,2,3,4,5,6,7},{8,9,10,11,12},{13,14}}\ \Bigg|\
          \ShadedTableau[]{{15,26,27},{16}}\ \Bigg|\
          \ShadedTableau[]{{17}}\ \Bigg|\
          \ShadedTableau[]{{25}}\  \Bigg|\
          \ShadedTableau[]{{18,19,20,21},{22,23,24}}\
\Bigg).
\end{align*}
Then the $\R$-module endomorphisms $\thea$, $\theta^\blam_\brho$ and $\theta^\brho_\bmu$ are given by
$$ \thea(x) = x y_{25}^4 y_{26} y_{27} e^\blam_\bmu, \quad
   \theta^\blam_\brho(x) = x y_{25}^3 y_{26} y_{27} e^\blam_\brho, \quad \text{ and } \quad
   \theta^\brho_\bmu(x) = x y_{25}y_{27} e^\brho_\bmu,$$
respectively, for all $x \in S^\bnu$ and where
$$e^{\blam}_{\bmu} = e^\blam_\brho = e^\brho_\bmu = \sum_{\bj \in I^{24}}e(j_1,\ldots,j_{24},0,2,3).$$
So, $\theta^\brho_\bmu\circ\theta^\blam_\brho=\thea\ne0$ and it induces a non-zero
$\R$-homomorphism $S^\blam\<7\> \rightarrow S^\bmu$ as claimed.
\end{Example}

\begin{proof}[Proof of \autoref{T:composition}] The condition
  $\blam\backslash\brho=\bmu\backslash\brho$ implies that
  $\bnu=\blam\cup\brho=\brho\cup\bmu=\blam\cup\bmu$ and, consequently, that
  $\res(\bnu\backslash\blam)=\res(\bnu\backslash\brho)=\res(\bnu\backslash\bmu)$
  since $(\blam,\brho)$ and $(\brho,\bmu)$ are cyclotomic Carter-Payne pairs.
  Let $\gamma=|\bnu\backslash\blam|$ and $J=\res(\bnu\backslash\blam)$ and note
  that if $j\le_J i$ then $\row\alpha^\bnu_\brho(i)=\row\alpha^\bnu_\brho(j)$
  and $\row\alpha^\bnu_\bmu(i)=\row\alpha^\bnu_\bmu(j)$. It
  follows that $(\blam,\bmu)$ is also a cyclotomic Carter-Payne pair. By
  \autoref{T:main} there exist non-zero maps
  $\Theta^\blam_\bmu$, $\Theta^\blam_\brho$ and $\Theta^\brho_\bmu$, as in the
  statement of the theorem, and it remains to show that
  $\Theta^\blam_\bmu=\Theta^\brho_\bmu\circ\Theta^\blam_\brho$.

  As in the proof of \autoref{T:main}, let $\tnulam$ be a standard
  $\bnu$-tableau such that $(\tnulam)_{\downarrow n}=\tlam$ and
  $\tnulam(\alpha^\bnu_\blam(i+1)) = \tnulam(\alpha^\bnu_\blam(i))+1$, for
  $i,i+1\in J$. Let
  $\t^\bnu_\brho$ be the unique standard $\bnu$-tableau such that $\t^\bnu_\brho
  \in \IniNu[\blam]$ and $(\t^\bnu_\brho)_{\downarrow n} = \t^\brho$ and define
  $\tnumu$ similarly.  Using these tableau, define maps $\thea$,
  $\theta^\blam_\brho$ and~$\theta^\brho_\bmu$ as in \autoref{E:thetaMain}. Then
  these three maps are $\R$-module endomorphisms of~$S^\bnu$ and by
  \autoref{L:ThetaIndependence} they depend only on the set of removable
  $j$-strips which occur between $\bnu\backslash\blam$ and~$\bnu\backslash\bmu$. Consequently,
  $\theta^\blam_\bmu=\theta^\brho_\bmu\circ\theta^\blam_\brho$ as endomorphisms
  of~$S^\bnu$. By the same reasoning, $\delta^\blam_\bmu = \delta^\blam_\brho +
  \delta^\brho_\bmu$.

  By the proof of \autoref{T:main}, the map $\theta^\brho_\bmu\circ\theta^\blam_\brho$ sends
  $\check S^\bnu_{\tnulam}$ into $\hat S^\bnu_{\tnumu}$. Therefore,
  $\theta^\brho_\bmu\circ\theta^\blam_\brho$ induces
  an $\R$-module homomorphism from $S^\blam\<\delta^\blam_\bmu\>$ to~$S^\bmu$ by
  \autoref{L:EllersMurray}.  As
  $\theta^\blam_\bmu=\theta^\brho_\bmu\circ\theta^\blam_\brho$, it follows that
  $\Theta^\blam_\bmu=\Theta^\brho_\bmu\circ\Theta^\blam_\brho$. In particular,
  the map $\Theta^\brho_\bmu\circ\Theta^\blam_\brho$ is a non-zero homomorphism in
  $\Hom_{\R}(S^\blam\<\delta^\blam_\bmu\>,S^\bmu)$.
\end{proof}

Before stating the second result in this section we introduce some new notation.
Suppose that $\blam=(\lambda^{(1)},\dots,\lambda^{(\ell)})\in\Parts$ and recall
that $\lambda^{(k)}=(\lambda^{(k)}_1,\lambda^{(k)}_2,\dots)$, for $1\le k\le\ell$.
Fix a row index $r=(l,a)$. Then $r$ determines the \textbf{row slice}
$\blam=\blam^{[1]}\sqcup\blam^{[2]}$ of~$\blam$, where
\begin{align*}
\blam^{[1]}&=\(\lambda^{(1)},\dots,\lambda^{(l-1)},
         (\lambda^{(l)}_1,\dots,\lambda^{(l)}_a)\)\\
\blam^{[2]}&=\((\lambda^{(l)}_{a+1},\lambda^{(l)}_{a+2},\dots),
         \lambda^{(l+1)},\dots,\lambda^{(\ell)}\).
\end{align*}
Define $n_1=|\blam^{[1]}|$,  $n_2=|\blam^{[2]}|$, $l_1=l$ and
$l_2=\ell-l+1$, for $b=1,2$. We consider the $l_b$-multipartition $\blam^{[b]}$ to be an element of
$\mathscr{P}^{\Lambda_b}_{n_b}$, where the weight $\Lambda_b\in P^+$ is determined by
the multisets of residues of the nodes in $\blam^{[b]}$, which we identify with a
subset of~$\blam$ in the obvious way, for $b=1,2$.

Iterating this construction, we  consider row slices
$\blam=\blam^{[1]}\sqcup\dots\sqcup\blam^{[s]}$. We emphasize that we consider
each slice $\blam^{[b]}$ as a subset of~$\blam$, so $\blam^{[b]}$ comes equipped
with a multiset of residues for $1\le b\le s$.

\begin{Theorem}
\label{T:Slice}
  Suppose that $e\in\{0,2,3,4,\dots\}$ and that $\blam,\bmu\in\Parts$ are
  multipartitions which admit row slices
  $\blam=\blam^{[1]}\sqcup\dots\sqcup\blam^{[s]}$ and
  $\bmu=\bmu^{[1]}\sqcup\dots\sqcup\bmu^{[s]}$ such that $\ell(\blam^{[b]})=\ell(\bmu^{[b]})$,
  $|\blam^{[b]}|=|\bmu^{[b]}|$ and $(\blam^{[b]},\bmu^{[b]})$ is
  a cyclotomic Carter-Payne pair, for $1\le b\le s$. Then
  $$\Hom_{\R}(S^\blam\<\delta^\blam_\bmu\>,S^\bmu)\ne0,$$
  where $\delta^\blam_\bmu=\delta^{\blam^{[1]}}_{\bmu^{[1]}}+\dots
  +\delta^{\blam^{[s]}}_{\bmu^{[s]}}\ge s$.
\end{Theorem}

If $\Lambda=\Lambda_0$ is a weight of level~$1$ then the main result of
\cite{LM:homs} implies that
$$\Hom_{\H}(\underline{S}^\blam,\underline{S}^\bmu)
\cong\bigotimes_{b=1}^s \Hom_{\mathscr{H}^{\Lambda_b}_{n_b}}(
          \underline{S}^{\blam^{[b]}},\underline{S}^{\bmu^{[b]}})$$
as vector spaces. Therefore, the $\H$-module analogue of \autoref{T:Slice} is a
corollary of \autoref{T:main} and \cite{LM:homs} when $\Lambda=\Lambda_0$.
\autoref{T:Slice} suggests that there is a higher level analogue of the
`row removal' results of \cite{LM:homs}.

\shadedBaseline=-4mm
\begin{Example}
\label{E:LastExample}
Let $e=3$ and let $\kappa=(0,2,0)$.  Consider
\begin{align*}
\blam&=(5,4,1|5,3|5)=(5,4,1|5,3)\sqcup(\emptyset|5),\quad\text{ and }\\
\bmu&=(6,5,1|3,3,2|3)=(6,5,1|3,3)\sqcup(2|3).
\end{align*}
Set $\bnu=\blam\cup\bmu=(6,5,1|5,3,2|5)$.  Then
$$\bnu= \Bigg(\ \ShadedTableau[(6,0),(5,-1)]{{0,1,2,0,1,2},{2,0,1,2,0},{1}}\ \Bigg|\
          \ShadedTableau[(4,0),(5,0),(1,-2),(2,-2)]{{2,0,1,2,0},{1,2,0},{0,1}}\ \Bigg|\
          \ShadedTableau[(4,0),(5,0)]{{0,1,2,0,1}}\ \Bigg),$$
where we have shaded all of the nodes in $\bnu\backslash (\blam \cap \bmu)$.
Using the obvious notation for the row slices of $\blam$ and $\bmu$ at row
$r=(2,2)$, we see that $(\blam^{[1]},\bmu^{[1]})$ and $(\blam^{[2]},\bmu^{[2]})$ are
both cyclotomic Carter-Payne pairs. Therefore,
$\Hom_{\R}(S^\blam\<5\>,S^\bmu)\ne0$ by \autoref{T:Slice}. To construct the
homomorphism $\thea$ in this example we would define
$$\tnulam=\Bigg(\ \ShadedTableau[(6,0),(5,-1)]{{1,2,3,4,5,24},{6,7,8,9,25},{10}}\ \Bigg|\
          \ShadedTableau[(1,-2),(2,-2)]{{11,12,13,14,15},{16,17,18},{26,27}}\ \Bigg|\
          \ShadedTableau[]{{19,20,21,22,23}}\ \Bigg)$$
and we define $\thea(x) =xy_{24}y_{25}y_{27}e^\blam_\bmu$.
\end{Example}

\begin{proof}[Sketch of the proof of \autoref{T:Slice}] We sketch the proof only
  in the special case when $s=2$. The general case is similar except that it
  requires more notation.

  Let $\bnu=\blam\cup\bmu=\blam^{[1]}\sqcup\blam^{[2]}\cup\bmu^{[1]}\sqcup\bmu^{[2]}$
  set $\gamma=|\bnu\backslash\blam|=\gamma_1+\gamma_2$, where
  $\gamma_1=|\blam^{[1]}\backslash\bmu^{[1]}|$ and
  $\gamma_2=|\blam^{[2]}\backslash\bmu^{[2]}|$. As in \autoref{E:LastExample},
  define $\tnulam$ to be the standard $\bnu$-tableau such that
  $(\tnulam)_{\downarrow n}=\tlam$, the numbers $n+1,\dots,n+\gamma_1$ appear
  in~$\bnu\backslash(\blam\cup\bmu^{[2]})$,
  $n+\gamma_1+1,\dots,n+\gamma_1+\gamma_2$ appear
  in~$\bnu\backslash(\blam\cup\bmu^{[1]})$, and $n+1,\dots,n+\gamma$ are
  entered into $\tnulam$ subject to the usual residue constraints.

  Using the tableau $\tnulam$ define $\R$-module endomorphisms $\theta^{[1]}$
  and $\theta^{[2]}$ of~$S^\bnu$ as in \autoref{E:thetaMain}. Then
  $\theta^{[1]}$  is right multiplication by a polynomial in
  $y_{n+1},\dots,y_{n+\gamma_1}$ times $e^\blam_\bmu$,  and $\theta^{[2]}$ is
  right multiplication by a polynomial in $y_{n+\gamma_1+1},\dots,y_{n+\gamma}$
  times $e^\blam_\bmu$. By construction, these two maps commute, so
  $\theta^{[1]}\circ\theta^{[2]}=\theta^{[2]}\circ\theta^{[1]}$ is an
  $\R$-module endomorphism of~$S^\bnu$. By \autoref{L:ThetaIndependence},
  $\theta^{[1]}$  depends only on~$\blam^{[1]}$ and $\bmu^{[1]}$ and,
  similarly,~$\theta^{[2]}$ depends only on~$\blam^{[2]}$ and $\bmu^{[2]}$. By
  \autoref{T:main}, the map~$\theta^{[1]}$ induces an $\R$-module homomorphism
  $\Theta^{[1]}$ from
  $S^{\blam^{[1]}\sqcup\blam^{[2]}}\<\delta^{\blam^{[1]}}_{\bmu^{[1]}}\>$ to
  $S^{\bmu^{[1]}\sqcup\blam^{[2]}}$. Similarly, because $\theta^{[2]}$ depends
  only on $\blam^{[2]}$ and $\bmu^{[2]}$, the map $\theta^{[2]}$ induces an
  $\R$-module homomorphism $\Theta^{[2]}$ from
  $S^{\bmu^{[1]}\sqcup\blam^{[2]}}\<\delta^{\blam^{[2]}}_{\bmu^{[2]}}\>$ to
  $S^{\bmu^{[1]}\sqcup\bmu^{[2]}}$. As in the proof of \autoref{T:composition},
  it follows that $\theta^{[2]}\circ\theta^{[1]}$ induces the $\R$-module
  homomorphism $\Theta^{[2]}\circ\Theta^{[1]}$ from
  $S^\blam\<\delta^\blam_\bmu\>$ to~$S^\bmu$.

  To complete the proof it remains to show that $\Theta^{[2]}\circ\Theta^{[1]}$
  is non-zero. By construction, $\theta^{[2]}$ only ever moves entries in the
  tableaux which are in the rows after $\ell(\blam^{[1]})=\ell(\bmu^{[1]})$. Let
  $\t^{[1]}_0,\dots,\t^{[1]}_{z_1+1}$ and $\t^{[2]}_0,\dots,\t^{[2]}_{z_2+1}$ be
  the tableaux constructed during the proof of \autoref{T:main} when showing
  that the induced maps~$\Theta^{[1]}$ and~$\Theta^{[2]}$, respectively, are
  well-defined and non-zero. Therefore, extending the slice notation to tableaux
  in the obvious way, and working modulo $\hat{S}^\bnu_\bmu$, the argument used
  to prove \autoref{E:dagger} shows that $\theta^{[2]}\circ\theta^{[1]}$ sends
  $v_{\tnulam}$ to $v_{\t^{[1]}_{z_1+1}\sqcup\t^{[2]}_{z_2+1}}$ plus a linear
  combination of terms~$v_\u$, where $\u\ne\t^{[1]}_{z_1+1}\sqcup\t^{[2]}_{z_2+1}$.
  It follows that $\Theta^{[2]}\circ\Theta^{[1]}\ne0$.
\end{proof}

\bibliographystyle{andrew}

\end{document}